\documentclass{article}

\usepackage[a4paper]{geometry}
\usepackage{indentfirst}
\usepackage[hidelinks]{hyperref}
\usepackage{enumitem}
\usepackage{amsmath,amsfonts,amssymb,amsthm,amsrefs}
\usepackage{mathtools}
\usepackage{mathrsfs}
\usepackage{stmaryrd}
\SetSymbolFont{stmry}{bold}{U}{stmry}{m}{n}
\usepackage{galois}
\usepackage[capitalise]{cleveref}
\usepackage{tikz-cd}
\usetikzlibrary{decorations.pathmorphing}

\newtheorem{thm}{Theorem}[section]
\newtheorem{lem}[thm]{Lemma}
\newtheorem{cor}[thm]{Corollary}
\newtheorem{prop}[thm]{Proposition}
\newtheorem{conj}[thm]{Conjecture}

\theoremstyle{definition}
\newtheorem{dfn}[thm]{Definition}
\newtheorem*{ack}{Acknowledgement}

\theoremstyle{remark}
\newtheorem{rmk}[thm]{Remark}
\newtheorem{eg}[thm]{Example}

\numberwithin{equation}{section}

\title{Weak admissibility of exponentially twisted cohomology associated with some nondegenerate functions}
\author{Peijiang Liu}
\date{January, 2025}

\begin{document}

\maketitle

\begin{abstract}
	\par In this article, we study the filtered $\Phi$-modules canonically attached to the exponentially twisted cohomology associated with some nondegenerate functions. Inspired by $p$-adic Hodge theory, we conjecture that those filtered $\Phi$-modules are weakly admissible. We show that this expectation is correct under some assumptions using the theory of Adolphson and Sperber.
\end{abstract}

\tableofcontents

\section*{Introduction}
\addcontentsline{toc}{section}{Introduction}

\par Let $p$ be a prime number, and let $K$ be a complete discrete valuation ring with mixed characteristic $(0,p)$. Let $O$ be the ring of integers of $K$, and let $k$ be the residue field of $O$. Let $X$ be a smooth scheme over $O$, with the special fiber $X_{k}$ over $k$, and the generic fiber $X_{K}$ over $K$. The \emph{Hodge filtration} on $H^{i}_{\mathrm{dR}}(X_{K})$ yields the \emph{Hodge polygon} (cf. \cite{fontaine1979modules}*{4.3.2}). If the Frobenius automorphism $a\mapsto{}a^{p}$ on $k$ is lifted to an automorphism on $K$, then the absolute Frobenius endomorphism on $X$ induces the Frobenius structure on $H^{i}(X_{k}/K)$, which yields the \emph{Newton polygon} (cf. \cite{fontaine1979modules}*{4.3.1}). By reformulating a conjecture of Katz in \cite{katz1971theorem}, Mazur conjectured in \cite{mazur1972frobenius} that if $X$ is proper and smooth over $O$, then the associated Newton polygon lies above the associated Hodge polygon. We call such property the \emph{Newton-above-Hodge} property for $X$. In \cite{mazur1973frobenius}, Mazur proved his conjecture under certain assumptions. Since $X$ is smooth and proper over $O$, the specialization map $H^{i}_{\mathrm{dR}}(X_{K})\rightarrow{}H^{i}_{\mathrm{rig}}(X_{k}/K)$ is an isomorphism of $K$-vector spaces. Using this isomorphism, we can canonically associate a filtered $\Phi$-module, namely a $K$-vector space equipped with a filtration and a Frobenius structure, to $X$. Then, we may regard the Newton-above-Hodge property for $X$ as a property for the associated filtered $\Phi$-modules. However, the Newton-above-Hodge property for filtered $\Phi$-modules is not a structural one. It is not stable under extensions and subquotients of filtered $\Phi$-modules, and the full subcategory consisting of filtered $\Phi$-modules that satisfies the Newton-above-Hodge property is not an abelian category. This deficit will be remediated by considering instead the \emph{weakly admissible} property, introduced by Fontaine in \cite{fontaine1979modules}*{4.1.4~D\'efinition}. For a filtered $\Phi$-module, being weakly admissible implies that it satisfies the Newton-above-Hodge property. Furthermore, the full subcategory consisting of weakly admissible filtered $\Phi$-module is an abelian category. One of the main theorems of $p$-adic Hodge theory, which states that the filtered $\Phi$-module associated with a proper and smooth scheme over $O$ is weakly admissible, reveals the hidden relationship between the Hodge filtration on $H^{i}_{\mathrm{dR}}(X_{K})$ and the Frobenius structure on $H^{i}_{\mathrm{rig}}(X_{k}/K)$.

\par Let $\mathbb{A}^{1}$ be the affine line over $O$, and let $\mathbb{T}^{n}$ be the $n$-dimensional torus over $O$. Now, we consider purely positive characteristic situation. Let $f:\mathbb{T}^{n}_{k}\rightarrow\mathbb{A}^{1}_{k}$ be a morphism. Following Dwork's philosophy in \cite{dwork1962zeta}, when $f$ satisfies certain good properties (i.e. $f$ is nondegenerate and $\dim\Delta(f)=n$, using notions that we introduce later), Adolphson and Sperber associated a Newton polygon and a Hodge polygon with $f$. In this introduction section, we call such polygons the \emph{AS-Newton polygon} and the \emph{AS-Hodge polygon} respectively. Adolphson and Sperber proved that the AS-Newton polygon lies above the AS-Hodge polygon (cf. \cite{adolphson1989exponential}*{Corollary~3.11} or \cite{adolphson1993twisted}*{Corollary~3.18}). However, giving a geometric interpretation to the mysterious combinatorially defined AS polygons is not straightforward. By a recent result of Li (cf. \cite{li2022exponential}*{Theorem~1.2}), the AS-Newton polygon is identified with the Newton polygon defined by the Frobenius structure on the exponentially twisted rigid cohomology $V_{\mathrm{rig}}$ over $\mathbb{T}^{n}_{k}$ associated with $f$. It is then natural to ask for a way to interpret the AS-Hodge polygon into the Hodge polygon coming from a geometric object on the generic fiber $\mathbb{T}^{n}_{K}$. For this, we consider the exponentially twisted de Rham cohomology $V_{\mathrm{dR}}$ over $\mathbb{T}^{n}_{K}$ associated with a morphism $\widehat{F}:\mathbb{T}^{n}_{K}\rightarrow\mathbb{A}^{1}_{K}$ determined by the Teichm\"uller lift of $f$. The idea is to find the relationship between the AS-Hodge polygon and the Hodge polygon defined by some Hodge filtration on $V_{\mathrm{dR}}$. The first huge obstacle to realizing this na\"ive idea is that the exponentially twisted de Rham cohomology is not in the realm of Hodge theory, and Hodge filtration a priori does not make sense. Recently, Sabbah and Yu developed irregular Hodge theory (cf. \cite{sabbah2018irregular}), generalizing Hodge theory. Irregular Hodge theory allows us to attach a canonical filtration, called the \emph{irregular Hodge filtration}, to $V_{\mathrm{dR}}$. Using a result of Yu in \cite{yu2014irregular}*{\S4}, we can show that the AS-Hodge polygon and the Hodge polygon defined by the irregular Hodge filtration on $V_{\mathrm{dR}}$ coincide. In parallel with the story about $p$-adic Hodge theory, it is now natural to ask if we can explain this Newton-above-Hodge phenomenon discovered by Adolphson and Sperber by showing weak admissibility of some filtered $\Phi$-module associated with the geometric information carried by $f$.

\par Let us explain our question in more detail. Let $\mathbb{F}_{q}$ denote the finite field with $q$ elements of characteristic $p$. Let $\mathrm{W}(\mathbb{F}_{q})$ denote the ring of Witt vectors over $\mathbb{F}_{q}$, and let $K_{0}$ be the field of fractions of $\mathrm{W}(\mathbb{F}_{q})$. Fix a primitive root of unity $\zeta_{p}$, and let $K_{1}=K_{0}(\zeta_{p})$. Let $O_{1}$ be the ring of integers of $K_{1}$. Fix a Dwork's uniformizer $\pi\in{}K_{1}$, namely the one satisfies $\pi^{p-1}+p=0$. Let $\mathcal{L}_{\pi}$ be the Dwork F-isocrystal over $\mathbb{A}^{1}_{\mathbb{F}_{q}}$. The exponentially twisted rigid cohomology, denoted by $V_{\mathrm{rig}}$ in the preceding paragraph, is the middle-degree rigid cohomology $H^{n}_{\mathrm{rig}}(\mathbb{T}^{n}_{\mathbb{F}_{q}}/K_{1},f^{\ast}\mathcal{L}_{\pi})$. The $p$-th power Frobenius automorphism on $\mathbb{F}_{q}$ is canonically lifted to an automorphism on $K_{1}$. Let $\phi_{f}$ be the Frobenius structure on $V_{\mathrm{rig}}$ induced by the absolute Frobenius endomorphism of $\mathbb{T}^{n}_{\mathbb{F}_{q}}$. Let $\hat{f}$ be the \emph{Teichm\"uller lift} of $f$ (cf. \cref{eq:lift}), and let $\widehat{F}=\pi\hat{f}$. Let $\nabla_{\widehat{F}}$ be the connection over $\mathbb{T}^{n}_{K_{1}}$ defined by $\nabla_{\widehat{F}}(1)=\pi\otimes{}d\hat{f}$. The exponentially twisted de Rham cohomology, denoted by $V_{\mathrm{dR}}$ in the preceding paragraph, is the middle-degree de Rham cohomology $H^{n}_{\mathrm{dR}}(\mathbb{T}^{n}_{K_{1}},\nabla_{\widehat{F}})$. Let $F^{\ast}_{\mathrm{irr}}$ denote the irregular Hodge filtration on $V_{\mathrm{dR}}$. There is a canonical morphism $\iota_{\widehat{F}}:V_{\mathrm{dR}}\rightarrow{}V_{\mathrm{rig}}$, which is often called the specialization map. Now, we are ready to state our question precisely in the following conjecture.

\begin{conj}\label{conj:main}
	\par If $f,\hat{f}$ are nondegenerate and $\dim\Delta(f)=n$, then $\iota_{\widehat{F}}$ is an isomorphism, and
	\begin{equation*}
		((V_{\mathrm{rig}},\phi_{f}),(V_{\mathrm{dR}},F^{\ast}_{\mathrm{irr}}),\iota_{\widehat{F}})\in\mathbf{M\overline{F}}{}^{\Phi}_{K_{1}}
	\end{equation*}
	is weakly admissible.
\end{conj}

\par The definition of $f,\hat{f}$ being nondegnerate can be found in \cref{dfn:nondegenerate}, and the definition of $\Delta(f)$ can be found at the beginning of \cref{section:np-module}. The tuple $((V_{\mathrm{rig}},\phi_{f}),(V_{\mathrm{dR}},F^{\ast}_{\mathrm{irr}}),\iota_{\widehat{F}})$ is not really a filtered $\Phi$-module in the sense of Fontaine in \cite{fontaine1979modules}*{1.2.1}. The filtration associated to a filtered $\Phi$-module is required to be indexed by integers, while $F^{\ast}_{\mathrm{irr}}$ is indexed by real numbers. Furthermore, the base field of the vector space endowed with a Frobenius structure associated with a filtered $\Phi$-module is $K_{0}$, while $V_{\mathrm{rig}}$ has base field $K_{1}$. Thus, in order to fit our situation, it is necessary to generalize Fontaine's definition of filtered $\Phi$-modules. From now on, an object obtained by such generalization will be called a \emph{filtered $\Phi$-module} over $K_{1}$, and a filtered $\Phi$-module in the sense of Fontaine in \cite{fontaine1979modules} will be called an \emph{integrally filtered $\Phi$-module} over $K_{0}$. The category consisting of filtered modules over $K_{1}$ is denoted by $\mathbf{M\overline{F}}{}^{\Phi}_{K_{1}}$. The definition of weakly admissible property for integrally filtered $\Phi$-module over $K_{0}$ can be naturally extended to filtered $\Phi$-modules over $K_{1}$. Moreover, the full subcategory consisting of weakly admissible filtered $\Phi$-modules over $K_{1}$ is an abelian category. The main result of this article is a proof of \cref{conj:main} under reasonable constraints. We now state the claim of our main result by the following Theorem.

\begin{thm}\label{thm:main}
	\par Assume that $p\neq{}2$. If $f,\hat{f}$ are nondegenerate and $\dim\Delta(f)=n$, then $\iota_{\widehat{F}}$ is an isomorphism, so that
	\begin{equation*}
		((V_{\mathrm{rig}},\phi_{f}),(V_{\mathrm{dR}},F^{\ast}_{\mathrm{irr}}),\iota_{\widehat{F}})\in\mathbf{M\overline{F}}{}^{\Phi}_{K_{1}}\rlap{.}
	\end{equation*}
	In addition, if $\ell_{\mathrm{HT}}(V_{\mathrm{dR}},F^{\ast}_{\mathrm{irr}})\leq{}p-2$, then $((V_{\mathrm{rig}},\phi_{f}),(V_{\mathrm{dR}},F^{\ast}_{\mathrm{irr}}),\iota_{\widehat{F}})$ is weakly admissible. Here, for a filtered module $(V,F^{\ast})$ over $K_{1}$, we denote by $\ell_{\mathrm{HT}}(V,F^{\ast})$ the difference between its maximal Hodge-Tate weight and its minimal Hodge-Tate weight (cf. \cref{dfn:hodge-tate}).
\end{thm}

\par The proof of \cref{thm:main} can be found at the end of \cref{section:proof}. We can actually compute $\ell_{\mathrm{HT}}(V_{\mathrm{dR}},F^{\ast}_{\mathrm{irr}})$ combinatorially from $f$ (cf. \cref{rmk:np-module}). It is worth noting that when $n=1$ and $p\neq{}2$, the requirement $\ell_{\mathrm{HT}}(V_{\mathrm{dR}},F^{\ast}_{\mathrm{irr}})\leq{}p-2$ is always satisfied, and in this situation, \cref{conj:main} is true without constraint (cf. \cref{eg:curve}). On the other hand, there are examples of $f$ where \cref{conj:main} is true but $\ell_{\mathrm{HT}}(V_{\mathrm{dR}},F^{\ast}_{\mathrm{irr}})>p-2$ (cf. \cref{eg:tensor}).

\par Next, we introduce the idea of the proof of \cref{thm:main}. Looking closely into Adolphson and Sperber's work in \cite{adolphson1989exponential} and \cite{adolphson1993twisted}, we construct a $K_{1}$-vector space $V_{\mathrm{NP}}$ and attach a Frobenius structure $\tilde{\phi}_{\mathrm{NP}}$ to $V_{\mathrm{NP}}$. The Newton polygon defined by this Frobenius structure coincides with the AS-Newton polygon. By Li's work in \cite{li2022exponential}, we obtain a canonical isomorphism $T_{\mathrm{rig}}:V_{\mathrm{NP}}\rightarrow{}V_{\mathrm{rig}}$ compatible with the Frobenius structures. Using an algorithm similar to the one Adolphson and Sperber developed in \cite{adolphson1989exponential} and \cite{adolphson1993twisted}, we decompose the specialization map $\iota_{\widehat{F}}:V_{\mathrm{dR}}\rightarrow{}V_{\mathrm{rig}}$ into $\iota_{\mathrm{dR}}:V_{\mathrm{dR}}\rightarrow{}V_{\mathrm{NP}}$ and $\iota_{\mathrm{rig}}:V_{\mathrm{NP}}\rightarrow{}V_{\mathrm{rig}}$. Using a result of Bourgeois in \cite{bourgeois1999annulation}, we can show that $\iota_{\mathrm{rig}}$ is an isomorphism. In general $\iota_{\mathrm{rig}}$ is not compatible with the Frobenius structures $\tilde{\phi}_{\mathrm{NP}}$ and $\phi_{f}$. Furthermore, using another result of Adolphson and Sperber in \cite{adolphson1997twisted}, we can show that $\iota_{\mathrm{dR}}$ is also an isomorphism. This proves the claim that the specialization map $\iota_{\widehat{F}}$ is an isomorphism.
\begin{equation*}
\begin{tikzcd}
	V_{\mathrm{dR}}
		\ar[r,"\iota_{\mathrm{dR}}"]
		\ar[dr,bend right,"\iota_{\widehat{F}}"']
		&
	V_{\mathrm{NP}}
		\ar[r,"T_{\mathrm{NP}}"]
		\ar[d,"\iota_{\mathrm{rig}}"']
		\ar[dr,"T_{\mathrm{rig}}"description]
		&
	V_{\mathrm{NP}}
		\ar[d,"\iota_{\mathrm{rig}}"]
		\\
		&
	V_{\mathrm{rig}}
		\ar[r,"\sim"]
		&
	V_{\mathrm{rig}}
\end{tikzcd}
\end{equation*}
Set $\hat{\phi}_{\mathrm{NP}}\coloneqq\iota_{\mathrm{rig}}^{-1}\comp\phi_{f}\comp\iota_{\mathrm{rig}}$, then it gives another Frobenius structure on $V_{\mathrm{NP}}$ that is different from $\tilde{\phi}_{\mathrm{NP}}$. Moreover, the work of Adolphson and Sperber in \cite{adolphson1989exponential} and \cite{adolphson1993twisted} yields a filtration $F^{\ast}_{\mathrm{NP}}$ on $V_{\mathrm{NP}}$. Using a result of Yu in \cite{yu2014irregular}*{\S4}, we can show that $\iota_{\widehat{F}}(F^{\ast}_{\mathrm{irr}})=F^{\ast}_{\mathrm{\mathrm{NP}}}$. Now, we obtain two filtered $\Phi$-modules $((V_{\mathrm{NP}},\tilde{\phi}_{\mathrm{NP}}),(V_{\mathrm{NP}},F^{\ast}_{\mathrm{NP}}),\mathrm{id})$ and $((V_{\mathrm{NP}},\hat{\phi}_{\mathrm{NP}}),(V_{\mathrm{NP}},F^{\ast}_{\mathrm{NP}}),\mathrm{id})$, where
\begin{equation*}
	((V_{\mathrm{rig}},\phi_{f}),(V_{\mathrm{dR}},F^{\ast}_{\mathrm{irr}}),\iota_{\widehat{F}})\cong((V_{\mathrm{NP}},\hat{\phi}_{\mathrm{NP}}),(V_{\mathrm{NP}},F^{\ast}_{\mathrm{NP}}),\mathrm{id})\rlap{.}
\end{equation*}
We may abbreviate their notations to $(V_{\mathrm{NP}},\tilde{\phi}_{\mathrm{NP}},F^{\ast}_{\mathrm{NP}})$ and $(V_{\mathrm{NP}},\hat{\phi}_{\mathrm{NP}},F^{\ast}_{\mathrm{NP}})$ respectively. Though there exists a natural automorphism $T_{\mathrm{NP}}\coloneqq\iota_{\mathrm{rig}}^{-1}\comp{}T_{\mathrm{rig}}\comp\iota_{\mathrm{rig}}$ on $V_{\mathrm{NP}}$ that is compatible with $\tilde{\phi}_{\mathrm{NP}}$ and $\hat{\phi}_{\mathrm{NP}}$, this automorphism is not compatible with the filtration $F^{\ast}_{\mathrm{NP}}$. In other words, there is no morphism between $(V_{\mathrm{NP}},\tilde{\phi}_{\mathrm{NP}},F^{\ast}_{\mathrm{NP}})$ and $(V_{\mathrm{NP}},\hat{\phi}_{\mathrm{NP}},F^{\ast}_{\mathrm{NP}})$ in general.

\par We introduce a property, called the \emph{NP-agreeable} property (cf. \cref{dfn:np-agreeable}), for filtered $\Phi$-modules. A benefit of considering this property is that, despite its being stronger than the weakly admissible property, it is much easier to check. By elaborating the arguments of Adolphson and Sperber in \cite{adolphson1989exponential} and \cite{adolphson1993twisted} which show the Newton-above-Hodge property, we can show that $(V_{\mathrm{NP}},\tilde{\phi}_{\mathrm{NP}},F^{\ast}_{\mathrm{NP}})$ is NP-agreeable, so that it is weakly admissible. Now, the idea is to show that replacing the Frobenius structure $\tilde{\phi}_{\mathrm{NP}}$ with $\hat{\phi}_{\mathrm{NP}}$ preserves the NP-agreeability. In order to do this, we need a good understanding of the automorphism $T_{\mathrm{NP}}$ on $V_{\mathrm{NP}}$. By some careful estimate, we find that $(V_{\mathrm{NP}},\tilde{\phi}_{\mathrm{NP}},F^{\ast}_{\mathrm{NP}})$ is agreeable if $\ell_{\mathrm{HT}}(V_{\mathrm{dR}},F^{\ast}_{\mathrm{irr}})\leq{}p-2$. This implies the claim of \cref{thm:main}.

\par Finally, we introduce the organization of this article. In \cref{section:definition}, we generalize the definition of integrally filtered $\Phi$-modules over $K_{0}$ by introducing filtered $\Phi$-modules over extensions of $K_{0}$, and extend the definition of weakly admissible property to these new objects. We also show that the full subcategory of weakly admissible filtered $\Phi$-modules is an abelian category.

\par \cref{section:cohomology} is a review on exponentially twisted cohomology. We describe the exponentially twisted rigid cohomology $V_{\mathrm{rig}}$ by a chain complex, and construct the chain map that induces the Frobenius structure $\phi_{f}$. We also describe the exponentially twisted de Rham cohomology $V_{\mathrm{dR}}$ by a chain complex, and construct the chain map which induces the specialization map.

\par \cref{section:np-module} serves as a preparation for the proof of \cref{thm:main}. We prove that the exponentially twisted cohomology associated with certain functions defines filtered $\Phi$-modules over $K_{1}$. The definition of $(V_{\mathrm{NP}},\tilde{\phi}_{\mathrm{NP}},F^{\ast}_{\mathrm{NP}})$ and $(V_{\mathrm{NP}},\hat{\phi}_{\mathrm{NP}},F^{\ast}_{\mathrm{NP}})$ can also be found in this section.

\par In \cref{section:proof}, we prove \cref{thm:main}. In the first half, we show that $(V_{\mathrm{NP}},\tilde{\phi}_{\mathrm{NP}},F^{\ast}_{\mathrm{NP}})$ is NP-agreeable and, in particular, weakly admissible. In the second half, we estimate $T_{\mathrm{NP}}$ carefully and complete the proof.

\par In \cref{section:example}, we give three examples and make some further discussions. The first example shows that \cref{conj:main} is true without extra constraints if $n=1$ and $p\neq{}2$. The second example clarifies that $T_{\mathrm{NP}}$ is not compatible with the filtration $F^{\ast}_{\mathrm{NP}}$ in general. The third example demonstrates that the condition $\ell_{\mathrm{HT}}(V_{\mathrm{dR}},F^{\ast}_{\mathrm{irr}})$ is not indispensable. After these examples, we raise some questions about how \cref{conj:main} could be upgraded to fit more general setups.

\begin{ack}
	\par I would like to thank Jeng-Daw Yu, Takahiro Saito, and Yichen Qin for their helpful discussions on irregular Hodge theory. I also want to thank Takeshi Tsuji (Tsuji-Sensei) for his kind discussion about $p$-adic Hodge theory, and I really appreciate his inspiring advice. Let me express my thanks to Joe Kramer-Miller for his valuable discussion, before which, \cref{lem:valuation.uniformizer} was only claiming a lower bound for the valuation. Without his pointing out that the valuation is actually computable, we were not even able to expect the constraint for the Hodge-Tate length in \cref{thm:main}. Finally, please allow me to express my heartful gratitude to my advisor Tomoyuki Abe (Abe-Sensei). I am very grateful for the inspiring seminars with Abe-Sensei, and I am indebted to him for his generous advice and warm encouragement. It is such an honor for me to be a student of him.
\end{ack}

\section{Filtered \texorpdfstring{$\Phi$}{Phi}-modules}
\label{section:definition}

\par The goal of this section is to extend the weak admissible property to a generalized version of filtered $\Phi$-modules. Originally, a filtered $\Phi$-module in the sense of Fontiane in \cite{fontaine1979modules}*{1.2.1} is a vector space over a certain base field equipped with two extra structures: a filtration that is integrally indexed, and a bijective semilinear endomorphism. We need to generalize this definition because of two reasons: one is that irregular Hodge filtrations are indexed by real numbers, and the other is that exponentially twisted cohomology groups are vector spaces over a larger base field. The weakly admissible property defined by Fontaine in \cite{fontaine1979modules}*{4.1.4~D\'efinition} is naturally extended to these new objects.

\par Let us fix some basic notations. Let $\mathbb{F}_{p}$ be the finite field with $p$ elements of characteristic $p$. Let $k$ be an algebraic extension of $\mathbb{F}_{p}$, and let $\mathrm{W}(k)$ denote the ring of Witt vectors over $k$. Let $K_{0}$ be the field of fractions of $\mathrm{W}(k)$. For $a\in{}k$, we denote by $\hat{a}\in{}K$ its Teichm\"uller lift. The isomorphism $\mathrm{Gal}(K_{0}/\mathbb{Q}_{p})\cong\mathrm{Gal}(k/\mathbb{F}_{p})$ allows us to uniquely lift the absolute Frobenius automorphism $a\mapsto{}a^{p}$ on $k$ to an automorphism on $K_{0}$, which we call the \emph{Frobenius automorphism}. Fix an algebraic closure $\overline{\mathbb{Q}}_{p}$ of $\mathbb{Q}_{p}$, and let $\mathbb{C}_{p}$ be the completion of $\overline{\mathbb{Q}}_{p}$. Let $\operatorname{ord}:\mathbb{C}_{p}\rightarrow\mathbb{R}\cup\{\infty\}$ be the additive $p$-adic valuation on $\mathbb{C}_{p}$ normalized by $\operatorname{ord}p=1$, and let $\vert\lambda\rvert_{p}=p^{-\operatorname{ord}\lambda}$ denote the $p$-adic norm of $\lambda\in\mathbb{C}_{p}$ corresponding to this valuation.

\begin{dfn}\label{dfn:filtration}
	\par Let $K$ be a field. A \emph{filtration} $F^{\ast}$ on a $K$-vector space $V$ is a collection $F^{\ast}V=\{F^{i}V\}_{i\in\mathbb{R}}$ of $K$-subspaces of $V$ such that $F^{j}V\subseteq{}F^{i}V$ for all $j>i$ in $\mathbb{R}$. Such a filtration is said to be \emph{exhaustive} if $\bigcup_{i\in\mathbb{R}}F^{i}V=V$, and \emph{separated} if $\bigcap_{i\in\mathbb{R}}F^{i}V=0$.
	\begin{enumerate}
		\item Let $\mathbf{Mod}_{K}$ be the category of finite-dimensional $K$-vector spaces.
		\item A \emph{filtered module} over $K$ is a pair $(V,F^{\ast})$, where $V\in\mathbf{Mod}_{K}$, and $F^{\ast}$ is a filtration on $V$ that is exhaustive and separated.
		\item Let $(V,F^{\ast})$ and $(V',F^{\ast})$ be filtered modules over $K$. A morphism $T:V'\rightarrow{}V$ in $\mathbf{Mod}_{K}$ is said to be \emph{filtration-compatible} if $T(F^{i}V')\subseteq{}F^{i}V$ for all $i\in\mathbb{R}$.
	\end{enumerate}
	Let $\mathbf{M\overline{F}}_{K}$ be the category consisting of filtered modules over $K$ and filtration-compatible morphisms. For an object $(V,F^{\ast})\in\mathbf{M\overline{F}}_{K}$, we sometimes omit $F^{\ast}$ and denote it by $V$ instead, if it causes no confusion.
\end{dfn}

\begin{dfn}\label{dfn:exact}
	\begin{enumerate}
		\item Let $T:(V',F^{\ast})\rightarrow(V,F^{\ast})$ be a morphism in $\mathbf{M\overline{F}}_{K}$. We say that $T$ is \emph{strict} if $T(F^{i}V')=T(V')\cap{}F^{i}V$ for all $i\in\mathbb{R}$.
		\item Let $(V',F^{\ast})\rightarrow(V,F^{\ast})\rightarrow(V'',F^{\ast})$ be a sequence in $\mathbf{M\overline{F}}_{K}$. We say that this sequence is \emph{exact} if the morphisms are strict, and the associated sequence $V'\rightarrow{}V\rightarrow{}V''$ in $\mathbf{Mod}_{K}$ is exact.
	\end{enumerate}
\end{dfn}

\begin{dfn}\label{dfn:hodge-tate}
	\par Let $(V,F^{\ast})\in\mathbf{M\overline{F}}_{K}$. For $i\in\mathbb{R}$, set $\mathrm{gr}^{i}(V,F^{\ast})=F^{i}V/F^{>i}V$, where $F^{>i}V=\bigcup_{j>i}F^{j}V$. We define the \emph{associated graded module} of $(V,F^{\ast})$ to be $\mathrm{gr}^{\ast}(V,F^{\ast})=\bigoplus_{i\in\mathbb{R}}\mathrm{gr}^{i}(V,F^{\ast})$.
	\begin{enumerate}
		\item We say that $i\in\mathbb{R}$ is a \emph{Hodge-Tate weight} of $V$ if $\dim_{K}\mathrm{gr}^{i}(V)\neq{}0$. Let $\mathfrak{W}_{\mathrm{HT}}(V,F^{\ast})$ be the set of all Hodge-Tate weights of $(V,F^{\ast})$.
		\item We define the \emph{Hodge-Tate length} of $(V,F^{\ast})$ to be
			\begin{equation*}
				\ell_{\mathrm{HT}}(V,F^{\ast})=\max\mathfrak{W}_{\mathrm{HT}}(V,F^{\ast})-\min\mathfrak{W}_{\mathrm{HT}}(V,F^{\ast})\rlap{.}
			\end{equation*}
	\end{enumerate}
\end{dfn}

\par A filtered module over $K$ in the sense of Fontaine in \cite{fontaine1979modules} is an object $V\in\mathbf{M\overline{F}}_{K}$ such that $\mathfrak{W}_{\mathrm{HT}}(V)\subseteq\mathbb{Z}$. We call such an object an \emph{integrally filtered module} over $K$.

\begin{dfn}\label{dfn:hodge.number}
	\par Let $V\in\mathbf{M\overline{F}}_{K}$. We define the \emph{Hodge number} of $V$ to be
	\begin{equation*}
		t_{\mathrm{H}}(V)=\sum_{i\in\mathfrak{W}_{\mathrm{HT}}(V)}i\cdot\dim_{K}\mathrm{gr}^{i}(V)\rlap{.}
	\end{equation*}
\end{dfn}

\par For $V\in\mathbf{M\overline{F}}_{K}$, the tensor product filtration on $V^{\otimes\dim_{K}V}$ induces a filtration on $\det{}V$. We note that $t_{\mathrm{H}}(V)=t_{\mathrm{H}}(\det{}V)$. In \cref{dfn:filtration.basis} we define a class of basis of and object in $\mathbf{M\overline{F}}_{K}$ that is useful for the proof of \cref{thm:main}. In \cref{dfn:hodge-tate.weight}, we fix a notation for later convenience. \cref{lem:hodge.inequality} serves as a preparation for the proof of \cref{thm:category}.

\begin{dfn}\label{dfn:filtration.basis}
	\par Let $(V,F^{\ast})\in\mathbf{M\overline{F}}_{K}$, and let $d=\dim_{K}V$. A basis $\{v_{i}\}_{i=1}^{d}$ of $V$ is said to be \emph{filtration-generating} if $F^{j}V=\langle{}v_{1},\dots,v_{d_{j}}\rangle_{K}$ for all $j\in\mathfrak{W}_{\mathrm{HT}}(V,F^{\ast})$, where $d_{i}=\dim_{K}F^{i}V$.
\end{dfn}

\begin{dfn}\label{dfn:hodge-tate.weight}
	\par Let $(V,F^{\ast})\in\mathbf{M\overline{F}}_{K}$. We define the \emph{Hodge-Tate weight} of $v\in{}V$ to be
	\begin{equation*}
		w_{\mathrm{HT}}(v)=\max\{i\in\mathfrak{W}_{\mathrm{HT}}(V,F^{\ast})\mid{}v\in{}F^{i}V\}\rlap{.}
	\end{equation*}
\end{dfn}

\begin{lem}\label{lem:hodge.inequality}
	\par Let $T:V'\rightarrow{}V$ be a morphism in $\mathbf{M\overline{F}}_{K}$. If $T$ is an isomorphism in $\mathbf{Mod}_{K}$, then $t_{\mathrm{H}}(V')\leq{}t_{\mathrm{H}}(V)$. This inequality is an equality if and only if $T$ is an isomorphism in $\mathbf{M\overline{F}}_{K}$.
\end{lem}
\begin{proof}
	\par We show the assertion by induction on $d=\dim_{K}V$. The assertion is obvious if $d=1$. Assuming the assertion for $1\leq{}d\leq{}r-1$, we prove it for $d=r$. Pick $v'\in{}V\setminus\{0\}$, and let $v=T(v')$. Then, we have the following commutative diagram of short exact sequences in $\mathbf{M\overline{F}}_{K}$.
	\begin{equation*}
	\begin{tikzcd}
		0
			\ar[r]
			&
		\langle{}v'\rangle_{K}
			\ar[r]
			\ar[d,"T"']
			&
		V'
			\ar[r]
			\ar[d,"T"']
			&
		V'/\langle{}v'\rangle_{K}
			\ar[r]
			\ar[d,"\overline{T}"']
			&
		0
			\\
		0
			\ar[r]
			&
		\langle{}v\rangle_{K}
			\ar[r]
			&
		V
			\ar[r]
			&
		V/\langle{}v\rangle_{K}
			\ar[r]
			&
		0
	\end{tikzcd}
	\end{equation*}
	Here $\langle{}v\rangle_{K}$ and $\langle{}v'\rangle_{K}$ are equipped with the subspace filtrations, while $V/\langle{}v\rangle_{K}$ and $V'/\langle{}v'\rangle_{K}$ are equipped with the quotient filtrations. Note that the vertical arrows are isomorphisms in $\mathbf{Mod}_{K}$. By the assertion for $d=1$, we get $t_{\mathrm{H}}(\langle{}v'\rangle_{K})\leq{}t_{\mathrm{H}}(\langle{}v\rangle_{K})$, and by the assertion for $d=r-1$, we get $t_{\mathrm{H}}(V'/\langle{}v'\rangle_{K})\leq{}t_{\mathrm{H}}(V/\langle{}v\rangle_{K})$. Those inequalities are equalities if and only if both $T:\langle{}v'\rangle_{K}\rightarrow\langle{}v\rangle_{K}$ and $\overline{T}:V'/\langle{}v'\rangle_{K}\rightarrow{}V/\langle{}v\rangle_{K}$ are isomorphisms in $\mathbf{M\overline{F}}_{K}$. Thus
	\begin{equation*}
		t_{\mathrm{H}}(V')=t_{\mathrm{H}}(\langle{}v'\rangle_{K})+t_{\mathrm{H}}(V'/\langle{}v'\rangle_{K})\leq{}t_{\mathrm{H}}(\langle{}v\rangle_{K})+t_{\mathrm{H}}(V/\langle{}v\rangle_{K})=t_{\mathrm{H}}(V)\rlap{,}
	\end{equation*}
	with equality if and only if $T:V'\rightarrow{}V$ is an isomorphism in $\mathbf{M\overline{F}}_{K}$. Hence, the assertion is true for $d=r$. By mathematical induction, we conclude that the assertion is true for $d\in\mathbb{Z}_{\geq{}1}$.
\end{proof}

\par To define $\Phi$-modules, we need a Frobenius automorphism on the base field. Since not all the extension of $K_{0}$ carries a Frobenius automorphism, we restrict our attention to a certain class of extensions of $K_{0}$ specified by the following definition.

\begin{dfn}\label{dfn:extension}
	\par Let $K_{\ast}$ be a finite extension of $K_{0}$ with residue field $k$. If the Frobenius automorphism on $K_{0}$ extends to an automorphism $\sigma$ on $K_{\ast}$, then we say that $K_{\ast}$ is a \emph{Frobenius extension} of $K_{0}$ with respect to $\sigma$.
\end{dfn}

\begin{dfn}\label{dfn:phi}
	\par Let $K_{\ast}$ be a Frobenius extension of $K_{0}$ with respect to $\sigma$. An endomorphism $\phi$ on a $K_{\ast}$-vector space is said to be $\sigma$-semilinear if $\phi(\lambda\cdot{}v)=\sigma(\lambda)\cdot\phi(v)$ for all $\lambda\in{}K_{\ast}$ and $v\in{}V$.
	\begin{enumerate}
		\item A \emph{$\Phi$-module} over $K_{\ast}$ is a pair $(V,\phi)$, where $V\in\mathbf{Mod}_{K_{\ast}}$, and $\phi$ is a $\sigma$-semilinear endomorphism on $V$ that is bijective.
		\item Let $(V,\phi)$ and $(V',\phi')$ be $\Phi$-modules over $K_{\ast}$. A morphism $T:V'\rightarrow{}V$ in $\mathbf{Mod}_{K_{\ast}}$ is said to be \emph{Frobenius-compatible} if $T\comp\phi'=\phi\comp{}T$.
	\end{enumerate}
	Let $\mathbf{Mod}{}^{\Phi}_{K_{\ast}}$ be the category consisting of $\Phi$-modules over $K_{\ast}$ and Frobenius-compatible morphisms. For an object $(V,\phi)\in\mathbf{Mod}{}^{\Phi}_{K_{\ast}}$, we sometimes omit $\phi$ and denote it by $V$ instead, if it causes no confusion.
\end{dfn}

\par A $\Phi$-module in the sense of Fontaine in \cite{fontaine1979modules} is an object in $\mathbf{Mod}{}^{\Phi}_{K_{0}}$.

\begin{dfn}\label{dfn:newton.number}
	\par Let $(V,\phi)\in\mathbf{Mod}{}^{\Phi}_{K_{\ast}}$
	\begin{enumerate}
		\item If $\dim_{K}V=1$, pick $v\in{}V\setminus\{0\}$. Let $\lambda\in{}K_{\ast}$ be characterized by $\phi(v)=\lambda\cdot{}v$. We define the \emph{Newton number} of $(V,\phi)$ to be $t_{\mathrm{N}}(V,\phi)=\operatorname{ord}\lambda$, which is independent of the choice of such $v$.
		\item We define the \emph{Newton number} of $(V,\phi)$ to be $t_{\mathrm{N}}(V,\phi)=t_{\mathrm{N}}(\det{}V,\det\phi)$.
	\end{enumerate}
\end{dfn}

\begin{dfn}\label{dfn:category}
	\par Let $K$ be a complete discrete valuation field extending $K_{\ast}$ with a perfect residue field of characteristic $p$. We define the category of \emph{filtered $\Phi$-modules} over $K/K_{\ast}$ to be
	\begin{equation*}
		\mathbf{M\overline{F}}{}^{\Phi}_{K/K_{\ast}}=\mathbf{Mod}{}^{\Phi}_{K_{\ast}}\times_{\mathbf{Mod}_{K}}\mathbf{M\overline{F}}_{K}\rlap{.}
	\end{equation*}
\end{dfn}

\par An object in $\mathbf{M\overline{F}}{}^{\Phi}_{K/K_{\ast}}$ consists of three data: a $\Phi$-module $(V,\phi)\in\mathbf{Mod}{}^{\Phi}_{K_{\ast}}$, a filtered module $(V_{K},F^{\ast})\in\mathbf{M\overline{F}}_{K}$, and an isomorphism $\iota:V_{K}\rightarrow{}K\otimes_{K_{\ast}}V$ in $\mathbf{Mod}_{K}$. We denote such an object by $((V,\phi),(V_{K},F^{\ast}),\iota)$. We also denote this object abstractly by $\boldsymbol{V}$. The \emph{dimension} of $\boldsymbol{V}$ is defined to be $\dim\boldsymbol{V}=\dim_{K_{\ast}}V=\dim_{K}V_{K}$. Furthermore, the \emph{Hodge number} of $\boldsymbol{V}$ is defined to be $t_{\mathrm{H}}(\boldsymbol{V})=t_{\boldsymbol{H}}(V_{K},F^{\ast})$, and the \emph{Newton number} of $\boldsymbol{V}$ is defined to be $t_{\mathrm{N}}(\boldsymbol{V})=t_{\mathrm{N}}(V,\phi)$. For objects $\boldsymbol{V}=((V,\phi),(V_{K},F^{\ast}),\iota)$ and $\boldsymbol{V}{}'=((V',\phi'),(V'_{K},F^{\ast}),\iota')$ in $\mathbf{M\overline{F}}{}^{\Phi}_{K/K_{\ast}}$, a morphism $\boldsymbol{T}:\boldsymbol{V}{}'\rightarrow\boldsymbol{V}$ in $\mathbf{M\overline{F}}{}^{\Phi}_{K/K_{\ast}}$ is a pair $(T,T_{K})$ such that $T_{K}=\iota^{-1}\comp(1\otimes{}T)\comp\iota'$, where $T:V'\rightarrow{}V$ gives a morphism in $\mathbf{Mod}{}^{\Phi}_{K_{\ast}}$, and $T_{K}:V'_{K}\rightarrow{}V_{K}$ gives a morphism in $\mathbf{M\overline{F}}_{K}$. We say that $\boldsymbol{T}$ is \emph{strict} if $T_{K}$ is strict. Moreover, for an object $((V,\phi),(V_{K},F^{\ast}),\iota)\in\mathbf{M\overline{F}}{}^{\Phi}_{K/K_{\ast}}$, if $V_{K}=K\otimes_{K_{\ast}}V$, namely $\iota$ is the identity map, then we denote it by $(V,\phi,F^{\ast})$. The full subcategory of $\mathbf{M\overline{F}}{}^{\Phi}_{K/K_{\ast}}$ consisting of such objects, denoted by $\mathbf{M\overline{F}}{}^{\Phi,\mathrm{id}}_{K/K_{\ast}}$, is equivalent to $\mathbf{M\overline{F}}{}^{\Phi}_{K/K_{\ast}}$. Natural objects associated with exponentially twisted cohomology are in $\mathbf{M\overline{F}}{}^{\Phi}_{K/K_{\ast}}$, but it is more convenient to use objects in $\mathbf{M\overline{F}}{}^{\Phi,\mathrm{id}}_{K/K_{\ast}}$ for certain proof, which is the reason why we introduce this full subcategory here. Furthermore, let $\mathbf{MF}{}^{\Phi,\mathrm{id}}_{K/K_{\ast}}$ be the full subcategory of $\mathbf{M\overline{F}}{}^{\Phi,\mathrm{id}}_{K/K_{\ast}}$ consisting of objects $(V,\phi,F^{\ast})$ such that $(V_{K},F^{\ast})$ is a integrally filtered. We note that $\mathbf{MF}{}^{\Phi,\mathrm{id}}_{K/K_{\ast}}$ coincides with the category of filtered $\Phi$-modules defined by Fontaine in \cite{fontaine1979modules}*{1.2.2}.

\begin{dfn}\label{dfn:weakly-admissible}
	\par Let $\boldsymbol{V}\in\mathbf{M\overline{F}}{}^{\Phi}_{K/K_{\ast}}$. We say that $\boldsymbol{V}$ is \emph{weakly admissible} if $t_{\mathrm{N}}(\boldsymbol{V}{}')\geq{}t_{\mathrm{H}}(\boldsymbol{V}{}')$ for all subobject $\boldsymbol{V}{}'$ of $\boldsymbol{V}$ in $\mathbf{M\overline{F}}{}^{\Phi}_{K/K_{\ast}}$, with equality when $\boldsymbol{V}{}'=\boldsymbol{V}$. Let $\mathbf{M\overline{F}}{}^{\Phi,\mathbf{w.a.}}_{K/K_{\ast}}$ denote the full subcategory of $\mathbf{M\overline{F}}{}^{\Phi}_{K/K_{\ast}}$ consisting of weakly admissible objects.
\end{dfn}

\par \cref{dfn:weakly-admissible} extends the weakly admissible property for objects in $\mathbf{MF}{}^{\Phi,\mathrm{id}}_{K/K_{\ast}}$ given by Fontaine in \cite{fontaine1979modules}*{4.1.4~D\'efinition} to objects in $\mathbf{M\overline{F}}{}^{\Phi}_{K/K_{\ast}}$. The following theorem is a generalized version of \cite{fontaine1979modules}*{4.2.1~Proposition}. Even though the proof is parallel, we include it for the sake of completeness.

\begin{thm}\label{thm:category}
	\begin{enumerate}
		\item Let $0\rightarrow\boldsymbol{V}{}'\rightarrow\boldsymbol{V}\rightarrow\boldsymbol{V}{}''\rightarrow{}0$ be a short exact sequence in $\mathbf{M\overline{F}}{}^{\Phi}_{K/K_{\ast}}$. If two of the three terms are weakly admissible, then so is the third one.
		\item Let $\boldsymbol{T}:\boldsymbol{V}{}'\rightarrow\boldsymbol{V}$ be a morphism in $\mathbf{M\overline{F}}{}^{\Phi}_{K/K_{\ast}}$. If $\boldsymbol{V}{}'$ and $\boldsymbol{V}$ are weakly admissible, then $\boldsymbol{T}$ is strict, and $\ker\boldsymbol{T}$ and $\operatorname{coker}\boldsymbol{T}$ are weakly admissible.
		\item The category $\mathbf{M\overline{F}}{}^{\Phi,\mathbf{w.a.}}_{K/K_{\ast}}$ is an abelian category.
	\end{enumerate}
\end{thm}
\begin{proof}
	\par We start from the first assertion. It is straightforward to verify that $\boldsymbol{V}$ and $\boldsymbol{V}{}'$ being weakly admissible implies that $\boldsymbol{V}{}''$ is weakly admissible, and that $\boldsymbol{V}$ and $\boldsymbol{V}{}''$ being weakly admissible implies that $\boldsymbol{V}{}'$ is weakly admissible. Therefore, it remains to show that if $\boldsymbol{V}{}'$ and $\boldsymbol{V}{}''$ being weakly admissible, then $\boldsymbol{V}$ is weakly admissible. Let $\boldsymbol{W}$ be a subobject of $\boldsymbol{V}$. Let $\boldsymbol{W}{}'=\boldsymbol{W}\cap\boldsymbol{V}{}'$, and let $\boldsymbol{\overline{W}}{}''=\boldsymbol{W}/\boldsymbol{W}{}'$. Then, we get a short exact sequence
	\begin{equation*}
		0\rightarrow\boldsymbol{W}{}'\rightarrow\boldsymbol{W}\rightarrow\boldsymbol{\overline{W}}{}''\rightarrow{}0
	\end{equation*}
	in $\mathbf{M\overline{F}}{}^{\Phi}_{K/K_{\ast}}$. Write $\boldsymbol{\overline{W}}{}''=((W'',\phi''),(W''_{K},\overline{F}{}^{\ast}),\iota'')$, where $\overline{F}{}^{\ast}$ denotes the quotient filtration. Replace $\overline{F}{}^{\ast}$ with the subspace filtration $F^{\ast}$ and let $\boldsymbol{W}{}''=((W'',\phi''),(W''_{K},F^{\ast}),\iota'')$, then $\boldsymbol{W}{}''$ is a subobject of $\boldsymbol{V}''$ in $\mathbf{M\overline{F}}{}^{\Phi}_{K/K_{\ast}}$. Note that the identity map $W''_{K}\rightarrow{}W''_{K}$ gives a morphism $(W''_{K},\overline{F}{}^{\ast})\rightarrow(W''_{K},F^{\ast})$ in $\mathbf{M\overline{F}}_{K}$. By \cref{lem:hodge.inequality}, we get $t_{\mathrm{H}}(W''_{K},\overline{F}{}^{\ast})\leq{}t_{\mathrm{H}}(W''_{K},F^{\ast})$. Thus
	\begin{align*}
		t_{\mathrm{H}}(\boldsymbol{W})=t_{\mathrm{H}}(\boldsymbol{W}{}')+t_{\mathrm{H}}(W''_{K},\overline{F}{}^{\ast})\leq{}&t_{\mathrm{H}}(\boldsymbol{W}{}')+t_{\mathrm{H}}(W''_{K},F^{\ast})\\
		\leq{}&t_{\mathrm{N}}(\boldsymbol{W}{}')+t_{\mathrm{N}}(W'',\phi)=t_{\mathrm{N}}(\boldsymbol{W})\rlap{.}
	\end{align*}
	When $\boldsymbol{W}=\boldsymbol{V}$, we have $t_{\mathrm{H}}(\boldsymbol{V})=t_{\mathrm{H}}(\boldsymbol{V}{}')+t_{\mathrm{H}}(\boldsymbol{V}{}'')=t_{\mathrm{N}}(\boldsymbol{V}{}')+t_{\mathrm{N}}(\boldsymbol{V}{}'')=t_{\mathrm{N}}(\boldsymbol{V})$. Thus $\boldsymbol{V}$ is weakly admissible.
	
	\par Note that the third assertion follows straightforwardly from the second one, so we only need to consider the second assertion. To show that $\boldsymbol{T}$ is strict, it suffices to show $\operatorname{coim}\boldsymbol{T}\cong\operatorname{im}\boldsymbol{T}$ in $\mathbf{M\overline{F}}{}^{\Phi}_{K/K_{\ast}}$. Here $\operatorname{coim}\boldsymbol{T}=((\operatorname{coim}T,\phi'),(\operatorname{coim}T_{K},\overline{F}{}^{\ast}),\iota')$ and $\operatorname{im}\boldsymbol{T}=((\operatorname{im}T,\phi),(\operatorname{im}T_{K},F^{\ast}),\iota)$, where $\overline{F}{}^{\ast}$ denotes the quotient filtration, and $F^{\ast}$ denotes the subspace filtration. Note that the canonical isomorphisms $\operatorname{coim}T\rightarrow\operatorname{im}T$ and $\operatorname{coim}T_{K}\rightarrow\operatorname{im}T_{K}$ in of vector spaces give a morphism $\operatorname{coim}\boldsymbol{T}\rightarrow\operatorname{im}\boldsymbol{T}$ in $\mathbf{M\overline{F}}{}^{\Phi}_{K}$. By \cref{lem:hodge.inequality}, we get $t_{\mathrm{H}}(\operatorname{coim}\boldsymbol{T})\leq{}t_{\mathrm{H}}(\operatorname{im}\boldsymbol{T})$. This inequality is an equality if and only if $\operatorname{coim}\boldsymbol{T}\rightarrow\operatorname{im}\boldsymbol{T}$ is an isomorphism in $\mathbf{M\overline{F}}{}^{\Phi}_{K}$. Since $\boldsymbol{V}{}'$ is weakly admissible, we have $t_{\mathrm{N}}(\operatorname{coim}\boldsymbol{T})\leq{}t_{\mathrm{H}}(\operatorname{coim}\boldsymbol{T})$. Since $\boldsymbol{V}$ is weakly admissible, we have $t_{\mathrm{N}}(\operatorname{im}\boldsymbol{T})\geq{}t_{\mathrm{H}}(\operatorname{im}\boldsymbol{T})$. Combining those inequalities, we get
	\begin{equation}\label{eq:inequality}
		t_{\mathrm{N}}(\operatorname{coim}\boldsymbol{T})\leq{}t_{\mathrm{H}}(\operatorname{coim}\boldsymbol{T})\leq{}t_{\mathrm{H}}(\operatorname{im}\boldsymbol{T})\leq{}t_{\mathrm{N}}(\operatorname{im}\boldsymbol{T})
	\end{equation}
	At the same time, note that the canonical isomorphism $\operatorname{coim}T\rightarrow\operatorname{im}T$ gives an isomorphism $(\operatorname{coim}T,\phi')\rightarrow(\operatorname{im}T,\phi)$ in $\mathbf{Mod}{}^{\Phi}_{K_{\ast}}$. Thus $t_{\mathrm{N}}(\operatorname{coim}\boldsymbol{T})=t_{\mathrm{N}}(\operatorname{im}\boldsymbol{T})$, which implies that all inequalities in \cref{eq:inequality} are equalities. In particular, we get $t_{\mathrm{H}}(\operatorname{coim}\boldsymbol{T})=t_{\mathrm{H}}(\operatorname{im}\boldsymbol{T})$, which implies $\operatorname{coim}\boldsymbol{T}\cong\operatorname{im}\boldsymbol{T}$ in $\mathbf{M\overline{F}}{}^{\Phi}_{K/K_{\ast}}$. Thus $\boldsymbol{T}$ is strict. Now, by the first assertion, to show that $\ker\boldsymbol{T}$ and $\operatorname{coker}\boldsymbol{T}$ are weakly admissible, it suffices to show that $\operatorname{im}\boldsymbol{T}$ is weakly admissible. Since $\operatorname{im}\boldsymbol{T}$ is a subobject of a weakly admissible object in $\mathbf{M\overline{F}}{}^{\Phi}_{K/K_{\ast}}$, we only need to prove $t_{\mathrm{N}}(\operatorname{im}\boldsymbol{T})=t_{\mathrm{H}}(\operatorname{im}\boldsymbol{T})$, which we have already shown.
\end{proof}

\par In the following definition, we define a class of basis of an object in $\mathbf{M\overline{F}}{}^{\Phi}_{K/K_{\ast}}$ that is useful for the proof of \cref{thm:main}.

\begin{dfn}\label{dfn:agreeable}
	\par Let $\boldsymbol{V}=((V,\phi),(V_{K},F^{\ast}),\iota)\in\mathbf{M\overline{F}}{}^{\Phi}_{K/K_{\ast}}$, and let $d=\dim\boldsymbol{V}$. By saying a basis of $\boldsymbol{V}$, we mean a basis of $V$.
	\begin{enumerate}
		\item A basis $\{v_{i}\}_{i=1}^{d}$ of $\boldsymbol{V}$ is said to be \emph{filtration-generating} if $\{v_{K,i}\}_{i=1}^{d}$ is a filtration-generating basis of $V_{K}$, where $v_{K,i}=\iota^{-1}(1\otimes{}v_{i})$.
		\item Let $\{v_{i}\}_{i=1}^{d}$ be a filtration generating basis of $\boldsymbol{V}$. We say that $\{v_{i}\}_{i=1}^{d}$ is \emph{agreeable} if $\operatorname{ord}A_{\phi}(v_{i},v_{j})\geq{}w_{\mathrm{HT}}(v_{K,j})$ for all $i,j=1,\dots,d$, where $A_{\phi}(v_{i},v_{j})\in{}K_{\ast}$ is characterised by
			\begin{equation*}
				\phi(v_{i})=\sum_{j=1}^{d}A_{\phi}(v_{i},v_{j})\cdot{}v_{j}\rlap{.}
			\end{equation*}
	\end{enumerate}
\end{dfn}

\begin{prop}\label{prop:agreeable}
	\par If $\boldsymbol{V}\in\mathbf{M\overline{F}}{}^{\Phi}_{K/K_{\ast}}$ admits an agreeable basis, then $t_{\mathrm{N}}(\boldsymbol{V})\geq{}t_{\mathrm{H}}(\boldsymbol{V})$.
\end{prop}
\begin{proof}
	\par Let $d=\dim\boldsymbol{V}$, and let $\{v_{i}\}_{i=1}^{d}$ be an agreeable basis of $\boldsymbol{V}$. Let $\boldsymbol{A}_{\phi}$ denote the $d\times{}d$ matrix whose $(i,j)$ entry is equal to $A_{\phi}(v_{i},v_{j})$ for $i,j=1,\dots,d$. Note that $t_{\mathrm{N}}(\boldsymbol{V})=\operatorname{ord}\det\boldsymbol{A}_{\phi}$, where
	\begin{equation*}
		\operatorname{ord}\det\boldsymbol{A}_{\phi}\geq\min\left\{\sum_{j=1}^{d}\operatorname{ord}A_{\phi}(v_{i_{j}},v_{j})\;\middle\vert\;{}i_{1},\dots,i_{d}=1,\dots,d\right\}\geq\sum_{j=1}^{d}w_{\mathrm{HT}}(v_{K,j})\rlap{.}
	\end{equation*}
	The first inequality follows from the definition of determinants, and the second inequality follows from the assumption that $\{v_{i}\}_{i=1}^{d}$ is agreeable. Since $\{v_{i}\}_{i=1}^{d}$ is also filtration-generating, we have $t_{\mathrm{H}}(\boldsymbol{V})=\sum_{j=1}^{d}w_{\mathrm{HT}}(v_{K,j})$. Therefore, we get
	\begin{equation*}
		t_{\mathrm{N}}(\boldsymbol{V})=\operatorname{ord}\det\boldsymbol{A}_{\phi}\geq\sum_{i=1}^{d}w_{\mathrm{HT}}(v_{K,i})=t_{\mathrm{H}}(\boldsymbol{V})\rlap{.}\qedhere
	\end{equation*}
\end{proof}

\section{Exponentially twisted cohomology}
\label{section:cohomology}

\par This section is a preliminary to the arguments in \cref{section:np-module}. The goal of this section is to recollect some knowledge of exponentially twisted cohomology, and clarify how it relates to the categories defined in \cref{section:definition}. We recall firstly the exponentially twisted rigid cohomology, and secondly the exponentially twisted de Rham cohomology. We also consider a canonical map between them, which is often called the specialization map.

\par From this section, we fix $k=\mathbb{F}_{q}$. Fix a primitive $p$-th root of unity $\zeta_{p}$, and let $K_{1}=K_{0}(\zeta_{p})$. The isomorphism $\mathrm{Gal}(K_{1}/\mathbb{Q}_{p}(\zeta_{p}))\cong\operatorname{Gal}(K_{0}/\mathbb{Q}_{p})$ implies that the Frobenius automorphism on $K_{0}$ is uniquely extended to an automorphism $\sigma$ on $K_{1}$ such that $\sigma(\zeta_{p})=\zeta_{p}$. Then $K_{1}$ is a Frobenius extension of $K_{0}$ with respect to $\sigma$. Let $O_{1}$ be the ring of integers of $K_{1}$.

\par We begin with the exponentially twisted rigid cohomology. First, we review some basic theory of overconvergent F-isocrystals over affine smooth varieties. Let $\Gamma_{0}$ be a finitely generated smooth $O_{1}$-algebra presented by $\Gamma_{0}=O_{1}[\boldsymbol{y}]/I_{0}$, where $\boldsymbol{y}=(y_{1},\dots,y_{r})$, and $I_{0}$ is an ideal of $O_{1}[\boldsymbol{y}]$. For $\boldsymbol{v}=(v_{1},\dots,v_{r})\in\mathbb{Z}^{r}_{\geq{}0}$, let $\lvert\boldsymbol{v}\rvert=\sum_{i=1}^{r}v_{i}$, and let $\boldsymbol{y}^{\boldsymbol{v}}=y_{1}^{v_{1}}\dots{}y_{r}^{v_{r}}$. Let
\begin{equation*}
	O_{1}[\boldsymbol{y}]^{\dagger}=\bigcup_{\lambda>1}\left\{\sum_{\boldsymbol{v}\in\mathbb{Z}^{r}_{\geq{}0}}A_{\boldsymbol{v}}\boldsymbol{y}^{\boldsymbol{v}}\in{}O_{1}[[\boldsymbol{y}]]\;\middle\vert\;\lvert{}A_{\boldsymbol{v}}\rvert_{p}\lambda^{\lvert\boldsymbol{v}\rvert}\rightarrow{}0\mbox{ as }\lvert\boldsymbol{v}\rvert\rightarrow\infty\right\}\rlap{,}
\end{equation*}
and let $\Gamma^{\dagger}_{0}=O_{1}[\boldsymbol{y}]^{\dagger}/I_{0}\cdot{}O_{1}[\boldsymbol{y}]^{\dagger}$. Let $\Gamma=K_{1}\otimes_{O_{1}}\Gamma_{0}$, and let $\Gamma^{\dagger}=K_{1}\otimes_{O_{1}}\Gamma^{\dagger}_{0}$. Let $\overline{\Gamma}=k\otimes_{O_{1}}\Gamma_{0}$, and let $X_{k}=\operatorname{Spec}\overline{\Gamma}$. Let $\varphi:\Gamma^{\dagger}\rightarrow\Gamma^{\dagger}$ denote the lift of the absolute Frobenius on $\overline{\Gamma}$. An \emph{overconvergent F-isocrystal} over $X_{k}$ consists of three data: a projective $\Gamma^{\dagger}$-module $\mathcal{E}$ of finite type, an integrable connection $\nabla:\mathcal{E}\rightarrow\mathcal{E}\otimes_{\Gamma^{\dagger}}\Omega^{1}_{\Gamma^{\dagger}/K_{1}}$, and an isomorphism $\varPhi:\mathcal{E}^{\varphi}\rightarrow\mathcal{E}$ of $\Gamma^{\dagger}$-modules such that $\varPhi\comp\nabla^{\varphi}=\nabla\comp\varPhi$. Here $\mathcal{E}^{\varphi}$ denotes the extension of scalars of $\mathcal{E}$ along $\varphi$, and $\nabla^{\varphi}$ is the connection on $\mathcal{E}^{\varphi}$ induced by $\nabla$. We denote such an overconvergent F-isocrystal by $(\mathcal{E},\nabla,\varPhi)$. We often omit $\nabla$ and $\varPhi$, and write $\mathcal{E}$ instead, if it causes no confusion. Let $\mathrm{DR}^{\ast}(\mathcal{E},\nabla)$ denote the de Rham complex of the connection $\nabla$ on $\mathcal{E}$. For each $i\in\mathbb{Z}$, the $i$-th degree \emph{rigid cohomology} of the overconvergent F-isocrystal $\mathcal{E}$ over $X_{k}$ is given by
\begin{equation*}
	H^{i}_{\mathrm{rig}}(X_{k}/K_{1},\mathcal{E})=H^{i}(\mathrm{DR}^{\ast}(\mathcal{E},\nabla))\rlap{.}
\end{equation*}
At the same time, note that $\varPhi$ yields a chain map $\phi:\mathrm{DR}^{\ast}(\mathcal{E},\nabla)\rightarrow\mathrm{DR}^{\ast}(\mathcal{E},\nabla)$, which induces a bijective $\sigma$-semilinear endomorphism on $H^{i}_{\mathrm{rig}}(X_{k}/K_{1},\mathcal{E})$. Denoting this endomorphism still by $\phi$, we get
\begin{equation*}
	(H^{i}_{\mathrm{rig}}(X_{k}/K_{1},\mathcal{E}),\phi)\in\mathbf{Mod}{}^{\Phi}_{K_{1}}\rlap{.}
\end{equation*}
The terminology exponentially twisted rigid cohomology means the rigid cohomology of certain type of overconvergent F-isocrystals: the pullbacks of the Dwork F-isocrystal over the affine line. We recollect some fact of the Dwork F-isocrystal. By \cite{dwork1962zeta}*{Lemma~4.1}, there exists a unique $\pi\in\mathbb{Q}_{p}(\zeta_{p})$ satisfying $\pi+\frac{\pi^{p}}{p}=0$, such that
\begin{equation*}
	\pi\equiv\zeta_{p}-1\bmod(\zeta_{p}-1)^{2}\rlap{.}
\end{equation*}
Note that $\pi$ is a uniformizer of $K_{1}$ and $\operatorname{ord}\pi=\frac{1}{p-1}$. Let $\mathbb{A}^{1}=\operatorname{Spec}O_{1}[t]$ be the affine line over $O_{1}$. The \emph{Dwork F-isocrystal} $\mathcal{L}_{\pi}$ over $\mathbb{A}^{1}_{k}$ associated with $\pi$ is the overconvergent F-isocrystal over $\mathbb{A}^{1}_{k}$ presented by $(\mathcal{L}_{\pi},\nabla_{\pi{}t},\varPhi_{t})$, where $\mathcal{L}_{\pi}=K_{1}[t]^{\dagger}$ is the free $K_{1}[t]^{\dagger}$-module of rank $1$, the connection $\nabla_{\pi{}t}$ on $\mathcal{L}_{\pi}$ is defined by $\nabla_{\pi{}t}(1)=\pi\otimes{}dt$, and $\varPhi_{t}:\mathcal{L}^{\varphi}_{\pi}\rightarrow\mathcal{L}_{\pi}$ is defined by $\varPhi_{t}(1)=\exp(\pi(t^{p}-t))$. We explain more in detail about this formulation. Firstly $\varPhi_{t}$ is well-defined, because by \cite{dwork1982lectures}*{21.1~Proposition}, we have $\exp(\pi(t^{p}-t))\in{}K_{1}[t]^{\dagger}$. Secondly, the connection $\nabla^{\varphi}_{\pi{}t}$ on $\mathcal{L}^{\varphi}_{\pi}$ is defined by $\nabla^{\varphi}_{\pi{}t}(1)=\pi{}pt^{p-1}\otimes{}dt$, and it is straightforward to verify $\varPhi_{t}\comp\nabla^{\varphi}_{\pi{}t}=\nabla_{\pi{}t}\comp\varPhi_{t}$. Now, we are ready to recall the exponentially twisted rigid cohomology. Let $A_{0}=O_{1}[x_{1},\dots,x_{n},(x_{1}\dots{}x_{n})^{-1}]$, and let $\mathbb{T}^{n}=\operatorname{Spec}A_{0}$ be the $n$-dimensional torus over $O_{1}$. For $\boldsymbol{u}=(u_{1},\dots,u_{n})\in\mathbb{Z}^{n}$, let $\lvert\boldsymbol{u}\rvert=\sum_{i=1}^{n}\lvert{}u_{i}\rvert$, and let $\boldsymbol{x}^{\boldsymbol{u}}=x_{1}^{u_{1}}\dots{}x_{n}^{u_{n}}$. Note that
\begin{equation*}
	A^{\dagger}=\bigcup_{\lambda>1}\left\{\sum_{\boldsymbol{u}\in\mathbb{Z}^{n}}A_{\boldsymbol{u}}\boldsymbol{x}^{\boldsymbol{u}}\;\middle\vert\;A_{\boldsymbol{u}}\in{}K_{1}\mbox{, }\lvert{}A_{\boldsymbol{u}}\rvert_{p}\lambda^{\lvert\boldsymbol{u}\rvert}\rightarrow{}0\mbox{ as }\lvert\boldsymbol{u}\rvert\rightarrow\infty\right\}\rlap{.}
\end{equation*}
Let $f:\mathbb{T}^{n}_{k}\rightarrow\mathbb{A}^{1}_{k}$ be the morphism defined by $t\mapsto{}f(\boldsymbol{x})$, where
\begin{equation}\label{eq:function}
	f(\boldsymbol{x})=\sum_{\boldsymbol{u}\in\mathbb{Z}^{n}}\alpha_{\boldsymbol{u}}\boldsymbol{x}^{\boldsymbol{u}}\in\overline{A}\coloneqq{}k\otimes_{O_{1}}A_{0}
\end{equation}
is a Laurent polynomial. We define the \emph{Teichm\"uller lift} of $f(\boldsymbol{x})$ to be
\begin{equation}\label{eq:lift}
	\hat{f}(\boldsymbol{x})=\sum_{\boldsymbol{u}\in\mathbb{Z}^{n}}\hat{\alpha}_{\boldsymbol{u}}\boldsymbol{x}^{\boldsymbol{u}}\in{}A\coloneqq{}K_{1}\otimes_{O_{1}}A_{0}\rlap{,}
\end{equation}
where $\hat{\alpha}_{\boldsymbol{u}}\in{}K_{0}$ denotes the Teichm\"uller lift of $\alpha_{\boldsymbol{u}}\in{}k$. Let $\widehat{F}=\pi\hat{f}$. We present the overconvergent F-isocrystal $f^{\ast}\mathcal{L}_{\pi}$ over $\mathbb{T}^{n}_{k}$ by the triple $(\mathcal{E}_{\widehat{F}},\nabla_{\widehat{F}},\varPhi_{f})$, where $\mathcal{E}_{\widehat{F}}=A^{\dagger}$ is the free $A^{\dagger}$-module of rank $1$, the connection $\nabla_{\widehat{F}}$ on $\mathcal{E}_{\widehat{F}}$ is defined by $\nabla_{\widehat{F}}(1)=\pi\otimes{}d\hat{f}$, and $\varPhi_{f}:\mathcal{E}^{\varphi}_{\widehat{F}}\rightarrow\mathcal{E}_{\widehat{F}}$ is defined by $\varPhi_{f}(1)=\exp(\pi(\varphi(\hat{f}(\boldsymbol{x}))-\hat{f}(\boldsymbol{x})))$. Here, we note that $\varphi$ is defined by the endomorphism of $A^{\dagger}$ given by
\begin{equation*}
	\sum_{\boldsymbol{u}\in\mathbb{Z}^{n}}A_{\boldsymbol{u}}\boldsymbol{x}^{\boldsymbol{u}}\mapsto\sum_{\boldsymbol{u}\in\mathbb{Z}^{n}}\sigma(A_{\boldsymbol{u}})\boldsymbol{x}^{p\boldsymbol{u}}\rlap{.}
\end{equation*}
Note that $\varPhi_{f}$ is well-defined, because we have $\exp(\pi(\varphi(\hat{f}(\boldsymbol{x}))-\hat{f}(\boldsymbol{x})))\in{}A^{\dagger}$ by the proof of \cite{dwork1962zeta}*{Lemma~4.1}. The connection $\nabla^{\varphi}_{\widehat{F}}$ on $\mathcal{E}^{\varphi}_{\widehat{F}}$ is defined by $\nabla^{\varphi}_{\widehat{F}}(1)=\pi\otimes{}d\varphi(\hat{f}(\boldsymbol{x}))$, and it is straightforward to verify that $\varPhi_{f}\comp\nabla^{\varphi}_{\widehat{F}}=\nabla_{\widehat{F}}\comp\varPhi_{f}$. For each $i\in\mathbb{Z}$, the $i$-th degree rigid cohomology of $f^{\ast}\mathcal{L}_{\pi}$ over $\mathbb{T}^{k}$ is given by
\begin{equation*}
	H^{i}_{\mathrm{rig}}(\mathbb{T}^{n}_{k}/K_{1},f^{\ast}\mathcal{L}_{\pi})=H^{i}(\mathrm{DR}^{\ast}(\mathcal{E}_{\widehat{F}},\nabla_{\widehat{F}}))\rlap{.}
\end{equation*}
More explicitly, the de Rham complex $\mathrm{DR}^{\ast}(\mathcal{E}_{\widehat{F}},\nabla_{\widehat{F}})$ is the chain complex defined by setting
\begin{equation*}
	\mathrm{DR}^{l}(\mathcal{E}_{\widehat{F}},\nabla_{\widehat{F}})=\bigoplus_{1\leq{}i_{1}<\dots<i_{l}\leq{}n}A^{\dagger}dx_{i_{1}}\wedge\dots\wedge{}dx_{i_{l}}
\end{equation*}
for $l\in\mathbb{Z}$, with the differential $\mathrm{DR}^{l}(\mathcal{E}_{\widehat{F}},\nabla_{\widehat{F}})\rightarrow\mathrm{DR}^{l+1}(\mathcal{E}_{\widehat{F}},\nabla_{\widehat{F}})$ given by
\begin{equation*}
	\xi{}dx_{i_{1}}\wedge\dots\wedge{}dx_{i_{l}}\mapsto\sum_{i=1}^{n}(\partial_{i}\xi+\pi\xi\partial_{i}\hat{f})dx_{i}\wedge{}dx_{1}\wedge\dots\wedge{}dx_{i_{l}}
\end{equation*}
for $\xi\in{}A^{\dagger}$ and $1\leq{}i_{1}<\dots<i_{l}\leq{}n$. Here, we set $\partial_{i}=\frac{\partial}{\partial{}x_{i}}$ for $i=1,\dots,n$. Furthermore, the isomorphism $\varPhi_{f}$ is represented by the map $\phi_{f}:A^{\dagger}\rightarrow{}A^{\dagger}$ given by $\xi(\boldsymbol{x})\mapsto\exp(\pi(\varphi(\hat{f}(\boldsymbol{x}))-\hat{f}(\boldsymbol{x})))\cdot\varphi(\xi(\boldsymbol{x}))$. We observe that
\begin{equation*}
	\nabla_{\widehat{F}}\comp\phi_{f}=p\cdot(\phi_{f}\comp\nabla_{\widehat{F}})\rlap{,}
\end{equation*}
which implies that $\phi_{f}$ induces a chain map $\phi_{f}:\mathrm{DR}^{\ast}(\mathcal{E}_{\widehat{F}},\nabla_{\widehat{F}})\rightarrow\mathrm{DR}^{\ast}(\mathcal{E}_{\widehat{F}},\nabla_{\widehat{F}})$. We denote the induced bijective $\sigma$-semilinear endomorphsm on $H^{i}(\mathrm{DR}^{\ast}(\mathcal{E}_{\widehat{F}},\nabla_{\widehat{F}}))$ still by $\phi_{f}$. Then, we get
\begin{equation*}
	(H^{i}_{\mathrm{rig}}(\mathbb{T}^{n}_{k}/K_{1},f^{\ast}\mathcal{L}_{\pi}),\phi_{f})\in\mathbf{Mod}{}^{\Phi}_{K_{1}}\rlap{.}
\end{equation*}

\par Next, we move on to the exponentially twisted de Rham cohomology. We first go back to the setup considered at the beginning of this section and review some general theory. Recall that $X_{K_{1}}=\operatorname{Spec}\Gamma$ is an affine smooth variety over $K_{1}$. Let $E$ be a projective $\Gamma$-module of finite type, and let $\nabla:E\rightarrow{}E\otimes_{\Gamma}\Omega^{1}_{\Gamma/K_{1}}$ be an integrable connection. Let $\mathrm{DR}^{\ast}(E,\nabla)$ denote the de Rham complex of the connection $\nabla$ on $E$. For each $i\in\mathbb{Z}$, the $i$-th degree \emph{de Rham cohomology} of $\nabla$ over $X_{K_{1}}$ is given by
\begin{equation*}
	H^{i}_{\mathrm{dR}}(X_{K_{1}},\nabla)=H^{i}(\mathrm{DR}^{\ast}(E,\nabla))\rlap{.}
\end{equation*}
If $\nabla$ is not regular, then it is not clear how to attach a Hodge theoretic filtration on $H^{i}_{\mathrm{dR}}(X_{K_{1}},\nabla)$. In order to attach a canonical filtration to $H^{i}_{\mathrm{dR}}(X_{K_{1}},\nabla)$, a generalized version of the Hodge filtrations is expected. This issue is first raised by Deligne in \cite{deligne5theorie}, where he attached some sensible filtrations on the de Rham cohomology of exponentially twisted connections over curvers. Later, irregular Hodge theory is established by Sabbah and Yu (cf. \cite{sabbah2018irregular}). More explicitly, Sabbah defined the category $\mathbf{irrMHM}$ of \emph{irregular mixed Hodge modules} containing the category $\mathbf{MHM}$ of \emph{mixed Hodge modules} as a full subcategory (cf. \cite{sabbah2018irregular}*{Theorem~0.2~(1)}). An object in $\mathbf{irrMHM}$ is equipped with a filtration indexed by real numbers, called the \emph{irregular Hodge filtration} (cf. \cite{sabbah2018irregular}*{Definition~2.22}). The irregular Hodge filtration generalizes the Hodge filtration, in the sense that the irregular Hodge filtration associated with an object in $\mathbf{MHM}$ coincides with the associated Hodge filtration (cf. \cite{sabbah2018irregular}*{Theorem~0.3~(1)}).

\par Now, we explain how irregular Hodge theory allows us to attach a canonical filtration on the exponentially twisted de Rham cohomology. Assume that $X_{K_{1}}$ is smooth and quasi-projective, and let $F$ be a global regular function on $X_{K_{1}}$. For a pair $(X_{K_{1}},F)$ who admits a \emph{good compactification} (cf. \cite{yu2014irregular}*{\S1}), Yu provided a way in \cite{yu2014irregular} to construct the irregular Hodge filtration on the exponentially twisted de Rham cohomology $H^{i}_{\mathrm{dR}}(X_{K_{1}},\nabla_{F})$. Here $\nabla_{F}$ denotes the exponentially twisted connection over $X_{K_{1}}$ associated with $F$, namely the connection over $X_{K_{1}}$ defined by $\nabla_{F}(1)=1\otimes{}dF$. Note that Yu was assuming the base field to be $\mathbb{C}$, but his construction of the irregular Hodge filtration (cf. \cite{yu2014irregular}*{p.~110 Definition}) can be applied without change in our situation. Denoting the irregular Hodge filtration by $F^{\ast}_{\mathrm{irr}}$, we get
\begin{equation*}
	(H^{i}_{\mathrm{dR}}(X_{K_{1}},\nabla_{F}),F^{\ast}_{\mathrm{irr}})\in\mathbf{M\overline{F}}_{K_{1}}\rlap{.}
\end{equation*}
If $F$ is a regular function on $\mathbb{T}^{n}_{K_{1}}$ that is nondegenerate (cf. \cref{dfn:nondegenerate}), then the pair $(\mathbb{T}^{n}_{k},F)$ admits a good compactification (cf. \cite{yu2014irregular}*{\S4}), and the irregular Hodge filtration $F^{\ast}_{\mathrm{irr}}$ is in particular defined on $H^{i}_{\mathrm{dR}}(\mathbb{T}^{n}_{K_{1}},\nabla_{F})$.

\par Finally, we consider the relationship between the exponentially twisted rigid cohomology and the exponentially twisted de Rham cohomology. If $(\mathcal{E},\nabla)=(\Gamma^{\dagger}\otimes_{\Gamma}E,1\otimes\nabla)$, then the inclusion $\Gamma\hookrightarrow\Gamma^{\dagger}$ induces a chain map $\mathrm{DR}^{\ast}(E,\nabla)\rightarrow\mathrm{DR}^{\ast}(\mathcal{E},\nabla)$. For each $i$, this chain map induces the specialization map $H^{i}_{\mathrm{dR}}(X_{K_{1}},\nabla)\rightarrow{}H^{i}_{\mathrm{rig}}(X_{k}/K_{1},\mathcal{E})$, which is a morphism in $\mathbf{Mod}_{K_{1}}$. In particular, considering the exponentially twisted cohomology over the $n$-dimensional torus, the inclusion $A\hookrightarrow{}A^{\dagger}$ induces the specialization map
\begin{equation}\label{eq:specialization}
	\iota_{\widehat{F}}:H^{i}_{\mathrm{dR}}(\mathbb{T}^{n}_{K_{1}},\nabla_{\widehat{F}})\rightarrow{}H^{i}_{\mathrm{rig}}(\mathbb{T}^{n}_{k}/K_{1},f^{\ast}\mathcal{L}_{\pi})\rlap{.}
\end{equation}
In the following section, we prove that $\iota_{\widehat{F}}$ is an isomorphism under certain conditions for $f$.

\section{Newton polyhedron modules}
\label{section:np-module}

\par This section is a preparation for the proof of \cref{thm:main}. We first prove that the morphism in \cref{eq:specialization} is an isomorphism under certain conditions, so that we can associate a filtered $\Phi$-module to $f$. Secondly, we construct two filtered $\Phi$-modules, one of them can be directly compared with the exponentially twisted cohomology, and the other one is closely related to the work of Adolphson and Sperber in \cite{adolphson1989exponential} and \cite{adolphson1993twisted}. Those filtered $\Phi$-modules are useful for the proof of \cref{thm:main}.

\par Let $\mathfrak{K}$ be a field. For a Laurent polynomial
\begin{equation*}
	\mathfrak{f}(\boldsymbol{x})=\sum_{\boldsymbol{u}\in\mathbb{Z}^{n}}\mathfrak{a}_{\boldsymbol{u}}\boldsymbol{x}^{\boldsymbol{u}}\in\mathfrak{K}[x_{1},\dots,x_{n},(x_{1}\dots{}x_{n})^{-1}]\rlap{,}
\end{equation*}
we define the \emph{support} of $\mathfrak{f}$ to be the finite set $\operatorname{supp}\mathfrak{f}=\{\boldsymbol{u}\in\mathbb{Z}^{n}\mid\mathfrak{a}_{\boldsymbol{u}}\neq{}0\}$. Let $\Delta(\mathfrak{f})$ denote the convex hull of $\operatorname{supp}\mathfrak{f}\cup\{\boldsymbol{0}\}$ in $\mathbb{R}^{n}$, which we call the \emph{Newton polyhedron} of $\mathfrak{f}$. Let $\dim\Delta(\mathfrak{f})$ denote the dimension of the smallest linear subspace of $\mathbb{R}^{n}$ that contains $\Delta(\mathfrak{f})$. Let $\mathrm{vol}(\mathfrak{f})$ denote the $\dim\Delta(\mathfrak{f})$-dimensional volume of $\Delta(\mathfrak{f})$ relative to the induced lattice. For a face $\tau\in\Delta(\mathfrak{f})$, let
\begin{equation*}
	\mathfrak{f}_{\tau}(\boldsymbol{x})=\sum_{\boldsymbol{u}\in\tau\cap\operatorname{supp}\mathfrak{f}}\mathfrak{a}_{\boldsymbol{u}}\boldsymbol{x}^{\boldsymbol{u}}\rlap{.}
\end{equation*}

\begin{dfn}\label{dfn:nondegenerate}
	\par Let $\mathfrak{K}$ be a field, and let $\overline{\mathfrak{K}}$ be an algebraic closure of $\mathfrak{K}$. A Laurent polynomial $\mathfrak{f}(\boldsymbol{x})\in\mathfrak{K}[x_{1},\dots,x_{n},(x_{1}\dots{}x_{n})^{-1}]$ is said to be \emph{nondegenerate}, if for every face $\tau$ of $\Delta(\mathfrak{f})$ that does not contain $\boldsymbol{0}$, the Laurent polynomials $\partial_{1}\mathfrak{f}_{\tau}(\boldsymbol{x}),\dots,\partial_{n}\mathfrak{f}_{\tau}(\boldsymbol{x})$ have no common zero in $(\overline{\mathfrak{K}}{}^{\times})^{n}$.
\end{dfn}

\par We retain the notations from \cref{section:cohomology}. Recall that $f(\boldsymbol{x})$ is the Laurent polynomial defined in \cref{eq:function}, and $\hat{f}$ is its Teichm\"uller lift. Recall that $\widehat{F}=\pi\hat{f}$.

\begin{thm}\label{thm:object.cohomology}
	\par Assume that $p\neq{}2$. If $f,\hat{f}$ are nondegenerate and $\dim\Delta(f)=n$, then the specialization map $\iota_{\widehat{F}}$ in \cref{eq:specialization} is an isomorphism, so that we get an object
	\begin{equation*}
		((H^{n}_{\mathrm{rig}}(\mathbb{T}^{n}_{k}/K_{1},f^{\ast}\mathcal{L}_{\pi}),\phi_{f}),(H^{n}_{\mathrm{dR}}(\mathbb{T}^{n}_{K_{1}},\nabla_{\widehat{F}}),F^{\ast}_{\mathrm{irr}}),\iota_{\widehat{F}})\in\mathbf{M\overline{F}}{}^{\Phi}_{K_{1}/K_{1}}\rlap{.}
	\end{equation*}
\end{thm}

\begin{rmk}\label{rmk:object.cohomology}
	\par Assume that $p\neq{}2$. Let $F(\boldsymbol{x})\in{}A_{0}=O_{1}[x_{1},\dots,x_{n},(x_{1}\dots{}x_{n})^{-1}]$ be a nondegenerate Laurent polynomial. If $\pi^{-1}F\in{}A_{0}$ and the reduction of $\pi^{-1}F$ by $\pi$ coincides with $f$, then we still have a canonical map $\iota_{F}:H^{n}_{\mathrm{dR}}(\mathbb{T}^{n}_{K_{1}},\nabla_{F})\rightarrow{}H^{n}_{\mathrm{rig}}(\mathbb{T}^{n}_{k}/K_{1},f^{\ast}\mathcal{L}_{\pi})$. If we assume in addition that $\Delta(F)=\Delta(f)$, and the $p$-adic distance between $\pi^{-1}F$ and $\hat{f}$ is less than $p^{-\frac{1}{p-1}}$, then by a strategy similar to that in the proof of \cref{thm:object.cohomology}, we can show that $\iota_{F}$ is an isomorphism, so that we obtain an object $((H^{n}_{\mathrm{rig}}(\mathbb{T}^{n}_{k}/K_{1},f^{\ast}\mathcal{L}_{\pi}),\phi_{f}),(H^{n}_{\mathrm{dR}}(\mathbb{T}^{n}_{K_{1}},\nabla_{F}),F^{\ast}_{\mathrm{irr}}),\iota_{F})$ in $\mathbf{M\overline{F}}{}^{\Phi}_{K_{1}/K_{1}}$
\end{rmk}

\par The proof of \cref{thm:object.cohomology} can be found right after \cref{cor:decomposition}. Now, we make some preparation for the proof. Let $\operatorname{cone}f$ be the conical hull of $\operatorname{supp}f$ in $\mathbb{R}^{n}$. Define the \emph{weight} of $\boldsymbol{u}\in\operatorname{cone}f$ to be
\begin{equation*}
	w(\boldsymbol{u})=\inf\{w\in\mathbb{R}_{\geq{}0}\mid\boldsymbol{u}\in{}w\cdot\Delta(f)\}\rlap{.}
\end{equation*}
Let $\mathrm{M}(f)=\mathbb{Z}^{n}\cap\operatorname{cone}f$, and let $\boldsymbol{x}^{\mathrm{M}(f)}$ denote the multiplicative monoid $\{\boldsymbol{x}^{\boldsymbol{u}}\mid\boldsymbol{u}\in\mathrm{M}(f)\}$. Let $R=k[\boldsymbol{x}^{\mathrm{M}(f)}]$ be the monoid algebra spanned by $\boldsymbol{x}^{\mathrm{M}(f)}$ over $k$. By \cite{adolphson1993twisted}*{Lemma~1.13.(c)}, there exists a positive integer $m$, such that $w(\mathrm{M}(f))\subseteq{}m^{-1}\mathbb{Z}$. For $i\in{}m^{-1}\mathbb{Z}$, set
\begin{equation*}
	R^{i}=\left\{\sum_{\boldsymbol{u}\in\mathrm{M}(f)}a_{\boldsymbol{u}}\boldsymbol{x}^{\boldsymbol{u}}\in{}R\;\middle\vert\;w(\boldsymbol{u})\leq{}-i\mbox{ if }a_{\boldsymbol{u}}\neq{}0\right\}\rlap{.}
\end{equation*}
Let $\overline{R}=\bigoplus_{i\in{}m^{-1}\mathbb{Z}}R^{(i)}$, where $R^{(i)}=R^{i}/R^{i+\frac{1}{m}}$. For $i=1,\dots,n$, let $\bar{f}_{i}(\boldsymbol{x})\in{}R^{(-1)}$ denote the image of $x_{i}\partial_{i}f(\boldsymbol{x})\in{}R^{-1}$ in $\overline{R}$, and let $D_{\bar{f},i}=x_{i}\partial_{i}+\bar{f}_{i}$. Define a chain complex $C^{\ast}(\overline{R},\bar{f})$ by setting
\begin{equation*}
	C^{l}(\overline{R},\bar{f})=\bigoplus_{1\leq{}i_{1}<\dots<i_{l}\leq{}n}\overline{R}\frac{dx_{i_{1}}}{x_{i_{1}}}\wedge\dots\wedge\frac{dx_{i_{l}}}{x_{i_{l}}}
\end{equation*}
for $l\in\mathbb{Z}$, with the differential $C^{l}(\overline{R},\bar{f})\rightarrow{}C^{l+1}(\overline{R},\bar{f})$ given by
\begin{equation*}
	\bar{\xi}\frac{dx_{i_{1}}}{x_{i_{1}}}\wedge\dots\wedge\frac{dx_{i_{l}}}{x_{i_{l}}}\mapsto\sum_{i=1}^{n}D_{\bar{f},i}\bar{\xi}\frac{dx_{i}}{x_{i}}\wedge\frac{dx_{i_{1}}}{x_{i_{1}}}\wedge\dots\wedge\frac{dx_{i_{l}}}{x_{i_{l}}}
\end{equation*}
for $\bar{\xi}\in\overline{R}$ and $1\leq{}i_{1}<\dots<i_{l}\leq{}n$.

\begin{thm}[\cite{adolphson1989exponential}*{Theorem~2.18}]\label{thm:decomposition.reduction}
	\par Assume that $f$ is nondegenerate and $\dim\Delta(f)=n$.
	\begin{enumerate}
		\item If $i\neq{}n$, then $\dim_{k}H^{i}(C^{\ast}(\overline{R},\bar{f}))=0$.
		\item There exists a finite subset $\mathrm{M}_{\mathrm{NP}}\subseteq\mathrm{M}(f)$, such that
			\begin{equation*}
				\overline{R}=\overline{V}_{\mathrm{NP}}\oplus\sum_{i=1}^{n}D_{\bar{f},i}\overline{R}\rlap{,}
			\end{equation*}
			where $\overline{V}_{\mathrm{NP}}=\langle\boldsymbol{x}^{\boldsymbol{u}}\mid\boldsymbol{u}\in\mathrm{M}_{\mathrm{NP}}\rangle_{k}$. Furthermore, we have $\lvert\mathrm{M}_{\mathrm{NP}}\rvert=n!\mathrm{vol}(f)$.\qed
	\end{enumerate}
\end{thm}

\begin{rmk}\label{rmk:decomposition.reduction}
	\par Assume that $f$ is nondegenerate and $\dim\Delta(f)=n$. Let $\overline{I}_{\bar{f}}$ denote the ideal of $\overline{R}$ generated by $\bar{f}_{1},\dots,\bar{f}_{n}$. For $i\in{}m^{-1}\mathbb{Z}$, there exists a finite subset
	\begin{equation*}
		\mathrm{M}^{(i)}_{\mathrm{NP}}\subseteq\{\boldsymbol{u}\in\mathrm{M}(f)\mid{}w(\boldsymbol{u})=i\}\rlap{,}
	\end{equation*}
	such that $\overline{V}{}^{(i)}_{\mathrm{NP}}=\langle\boldsymbol{x}^{\boldsymbol{u}}\mid\boldsymbol{u}\in\mathrm{M}^{(i)}_{\mathrm{NP}}\rangle_{k}$ is complimentary to $R^{(i)}\cap\overline{I}_{\bar{f}}$ in $R^{(i)}$, namely
	\begin{equation*}
		R^{(i)}=\overline{V}{}^{(i)}_{\mathrm{NP}}\oplus(R^{(i)}\cap\overline{I}_{\bar{f}})\rlap{.}
	\end{equation*}
	Setting $\mathrm{M}_{\mathrm{NP}}=\bigcup_{i\in{}m^{-1}\mathbb{Z}}\mathrm{M}^{(i)}_{\mathrm{NP}}$, we get $\overline{V}_{\mathrm{NP}}=\bigoplus_{i\in{}m^{-1}\mathbb{Z}}\overline{V}{}^{(i)}_{\mathrm{NP}}$ and $\overline{R}=\overline{V}_{\mathrm{NP}}\oplus\sum_{i=1}^{n}D_{\bar{f},i}\overline{R}$.
\end{rmk}

\par Fix a $m$-th root $\pi_{m}$ of $\pi$, and let $K_{m}=K_{1}(\pi_{m})$. Note that by setting $\sigma(\pi_{m})=\pi_{m}$, we obtain an extension of $\sigma$ to $K_{m}$. Let $O_{m}$ be the ring of integers of $K_{m}$. For $b\in\mathbb{R}_{\geq{}0}$ and $c\in\mathbb{R}$, let
\begin{equation*}
	L(b,c)=\left\{\sum_{\boldsymbol{u}\in\mathrm{M}(f)}A_{\boldsymbol{u}}\boldsymbol{x}^{\boldsymbol{u}}\;\middle\vert\;A_{\boldsymbol{u}}\in{}K_{m}\mbox{, }\operatorname{ord}A_{\boldsymbol{u}}\geq{}b\cdot{}w(\boldsymbol{u})+c\right\}\rlap{.}
\end{equation*}
For $b\in\mathbb{R}_{\geq{}0}$, let $L(b)=\bigcup_{c\in\mathbb{R}}L(b,c)$. Let $\varpi=p^{\frac{1}{p-1}}\in\mathbb{R}$. Let
\begin{align*}
	B={}&\left\{\sum_{\boldsymbol{u}\in\mathrm{M}(f)}A_{\boldsymbol{u}}\boldsymbol{x}^{\boldsymbol{u}}\;\middle\vert\;A_{\boldsymbol{u}}\in{}K_{m}\mbox{, }\lvert{}A_{\boldsymbol{u}}\rvert_{p}\varpi^{w(\boldsymbol{u})}\rightarrow{}0\mbox{ as }w(\boldsymbol{u})\rightarrow\infty\right\}\rlap{,}\\
	B_{0}={}&\left\{\sum_{\boldsymbol{u}\in\mathrm{M}(f)}A_{\boldsymbol{u}}\boldsymbol{x}^{\boldsymbol{u}}\;\middle\vert\;A_{\boldsymbol{u}}\in\pi^{w(\boldsymbol{u})}O_{m}\mbox{, }\lvert{}A_{\boldsymbol{u}}\rvert_{p}\varpi^{w(\boldsymbol{u})}\rightarrow{}0\mbox{ as }w(\boldsymbol{u})\rightarrow\infty\right\}\rlap{.}
\end{align*}
For $b>\frac{1}{p-1}$, we have $A\subseteq{}L(b)\subseteq{}B\subseteq{}L(\frac{1}{p-1})\subseteq{}A^{\dagger}_{m}$. In addition, if $c\geq{}0$, then $L(b,c)_{m}\subseteq{}B_{0}$. Here $A^{\dagger}_{m}=K_{m}\otimes_{K_{1}}A^{\dagger}$. Let $\mathfrak{A}=\{L(b)\mid{}b>0\}\cup\{A^{\dagger}_{m},B\}$. For $b>0$, set
\begin{equation*}
	\mathfrak{A}^{b}=\{\Lambda\in\mathfrak{A}\mid{}L(b)\subseteq\Lambda\}\rlap{.}
\end{equation*}
By \cite{dwork1962zeta}*{Lemma~4.1}, there exists a unique $\gamma\in\mathbb{Q}_{p}(\zeta_{p})$ satisfying $\sum_{i=1}^{\infty}\frac{\gamma^{p^{i}}}{p^{i}}=0$, such that
\begin{equation*}
	\gamma\equiv\zeta_{p}-1\bmod(\zeta_{p}-1)^{2}\rlap{.}
\end{equation*}
Note that $\operatorname{ord}\gamma=\frac{1}{p-1}$. For $l\in\mathbb{Z}_{\geq{}0}$, set $\gamma_{l}=\sum_{i=1}^{l}\frac{\gamma^{p^{i}}}{p^{i}}$. Let $\theta(t)=\sum_{l=0}^{\infty}\gamma_{l}t^{p^{l}}$, and let
\begin{equation*}
	\widetilde{F}(\boldsymbol{x})=\sum_{\boldsymbol{u}\in\operatorname{supp}f}\theta(\hat{\alpha}_{\boldsymbol{u}}\boldsymbol{x}^{\boldsymbol{u}})\rlap{.}
\end{equation*}
For $i=1,\dots,n$ and $F\in\{\widehat{F},\widetilde{F}\}$, we define an operator $D_{F,i}=x_{i}\partial_{i}+x_{i}\partial_{i}F$. For $\Lambda\in\mathfrak{A}^{\frac{p}{p-1}}\cup\{A^{\dagger},B_{0}\}$ and $F\in\{\widehat{F},\widetilde{F}\}$, we define a chain complex $C^{\ast}(\Lambda,F)$ by setting
\begin{equation*}
	C^{l}(\Lambda,F)=\bigoplus_{1\leq{}i_{1}<\dots<i_{l}\leq{}n}\Lambda\frac{dx_{i_{1}}}{x_{i_{1}}}\wedge\dots\wedge\frac{dx_{i_{l}}}{x_{i_{l}}}
\end{equation*}
for $l\in\mathbb{Z}$, with differential $C^{l}(\Lambda,F)\rightarrow{}C^{l+1}(\Lambda,F)$ given by $\omega\mapsto\sum_{i=1}^{n}D_{F,i}\omega\wedge\frac{dx_{i}}{x_{i}}$. Note that $C^{\ast}(\Lambda,F)$ is well-defined, because by \cref{lem:valuation.differential}, we have $x_{i}\partial_{i}F(\boldsymbol{x})\in{}L(\frac{p}{p-1},-1)$ for $i=1,\dots,n$. Moreover, the inclusion $\Lambda\hookrightarrow{}A^{\dagger}_{m}$ induces a chain map $\iota_{\Lambda}:C^{\ast}(\Lambda,F)\rightarrow{}C^{\ast}(A^{\dagger}_{m},F)$. Here $\iota_{\Lambda}$ is written as $\hat{\iota}_{\Lambda}$ if $F=\widehat{F}$, and as $\tilde{\iota}_{\Lambda}$ if $F=\widetilde{F}$.  For $i\in\mathbb{Z}$, this chain map induces a morphism
\begin{equation}\label{eq:inclusion}
	\iota_{\Lambda}:H^{i}(C^{\ast}(\Lambda,F))\rightarrow{}H^{i}(C^{\ast}(A^{\dagger}_{m},F))\in\mathbf{Mod}_{K_{m}}\rlap{.}
\end{equation}
For $l\in\mathbb{Z}$, let $\hat{\imath}:C^{l}(A^{\dagger},\widehat{F})\rightarrow\mathrm{DR}^{l}(\mathcal{E}_{\widehat{F}},\nabla_{\widehat{F}})$ be the map given by
\begin{equation*}
	\xi(\boldsymbol{x})\frac{dx_{i_{1}}}{x_{i_{1}}}\wedge\dots\wedge\frac{dx_{i_{l}}}{x_{i_{l}}}\mapsto{}x_{1}^{-1}\dots{}x_{n}^{-1}\xi(\boldsymbol{x})dx_{i_{1}}\wedge\dots\wedge{}dx_{i_{l}}\rlap{.}
\end{equation*}
It is straightforward to verify that $\hat{\imath}$ is an isomorphism, and $\hat{\imath}\comp\nabla_{\widehat{F}}=\nabla_{\widehat{F}}\comp\hat{\imath}$. Thus $\hat{\imath}$ induces an isomorphism $\hat{\imath}:C^{\ast}(A^{\dagger},\widehat{F})\rightarrow\mathrm{DR}^{\ast}(\mathcal{E}_{\widehat{F}},\nabla_{\widehat{F}})$ of chain complexes. For $i\in\mathbb{Z}$, this chain map induces an isomorphism
\begin{equation}\label{eq:logarithm}
	\hat{\imath}:H^{i}(C^{\ast}(A^{\dagger},\widehat{F}))\rightarrow{}H^{i}_{\mathrm{rig}}(\mathbb{T}^{n}_{k}/K_{1},f^{\ast}\mathcal{L}_{\pi})\in\mathbf{Mod}_{K_{1}}\rlap{.}
\end{equation}

\begin{lem}\label{lem:valuation.differential}
	\par Let $F\in\{\widehat{F},\widetilde{F}\}$. For $i=1,\dots,n$, we have $x_{i}\partial_{i}F(\boldsymbol{x})\in{}L(\frac{p}{p-1},-1)$.
\end{lem}
\begin{proof}
	\par We first consider the case where $F=\widehat{F}$. In this situation, for $i=1,\dots,n$ we have
	\begin{equation*}
		x_{i}\partial_{i}\widehat{F}(\boldsymbol{x})=\sum_{\boldsymbol{u}\in\operatorname{supp}f}\pi{}u_{i}\hat{\alpha}_{\boldsymbol{u}}\boldsymbol{x}^{\boldsymbol{u}}\rlap{,}
	\end{equation*}
	where $u_{i}$ is characterized by $\boldsymbol{u}=(u_{1},\dots,u_{n})$. Note that $w(\boldsymbol{x})\leq{}1$ if $\boldsymbol{u}\in\operatorname{supp}f$. Thus, we have
	\begin{equation*}
		\operatorname{ord}(\pi{}u_{i}\hat{\alpha}_{\boldsymbol{u}})\geq\operatorname{ord}\pi=\frac{1}{p-1}=\frac{p}{p-1}-1\geq\frac{p}{p-1}\cdot{}w(\boldsymbol{u})-1
	\end{equation*}
	for $\boldsymbol{u}\in\operatorname{supp}f$, which shows the assertion. Next, we consider the case where $F=\widetilde{F}$ and show that $x_{i}\partial_{i}\widetilde{F}(\boldsymbol{x})\in{}L(\frac{p}{p-1},-1)$ for $i=1,\dots,n$. In this situation, we have
	\begin{equation*}
		x_{i}\partial_{i}\widetilde{F}(\boldsymbol{x})=\sum_{\boldsymbol{u}\in\operatorname{supp}f}\sum_{l=0}^{\infty}\gamma_{l}p^{l}u_{i}\hat{\alpha}_{\boldsymbol{u}}^{p^{l}}\boldsymbol{x}^{p^{l}\boldsymbol{u}}\rlap{.}
	\end{equation*}
	Note that for $l\in\mathbb{Z}_{\geq{}0}$, we have $\operatorname{ord}\gamma_{l}=\operatorname{ord}\frac{\gamma^{p^{l+1}}}{p^{l+1}}=\frac{p^{l+1}}{p-1}-l-1$. Therefore, we have
	\begin{equation*}
		\operatorname{ord}(\gamma_{l}p^{l}u_{i}\hat{\alpha}_{\boldsymbol{u}}^{p^{l}})\geq\operatorname{ord}(\gamma_{l}p^{l})=\frac{p^{l+1}}{p-1}-1\geq\frac{p}{p-1}\cdot{}w(p^{l}\boldsymbol{u})-1
	\end{equation*}
	for $l\in\mathbb{Z}_{\geq{}0}$ and $\boldsymbol{u}\in\operatorname{supp}f$, which shows the assertion.
\end{proof}

\begin{lem}\label{lem:reduction}
	\par The map $\rho:B_{0}\rightarrow\overline{R}$ given by
	\begin{equation*}
		\sum_{\boldsymbol{u}\in\mathrm{M}(f)}A_{\boldsymbol{u}}\boldsymbol{x}^{\boldsymbol{u}}\mapsto\sum_{\boldsymbol{u}\in\mathrm{M}(f)}a_{\boldsymbol{u}}\boldsymbol{x}^{\boldsymbol{u}}\rlap{,}
	\end{equation*}
	is a ring homomorphism, where $a_{\boldsymbol{u}}\in{}k$ is the reduction of $\pi^{-w(\boldsymbol{u})}A_{\boldsymbol{u}}\in{}O_{m}$ modulo $\pi$.
\end{lem}
\begin{proof}
	\par The proof is parallel to that in \cite{adolphson1989exponential}*{Lemma~2.10}. We say that $\boldsymbol{u}$ and $\boldsymbol{u}'$ in $\mathrm{M}(f)$ are \emph{cofacial} if $w(\boldsymbol{u})^{-1}\boldsymbol{u}$ and $w(\boldsymbol{u}')^{-1}\boldsymbol{u}'$ lie on the same closed face of $\Delta(f)$. By \cite{adolphson1989exponential}*{Lemma~1.9.(c)}, we have $w(\boldsymbol{u}+\boldsymbol{u}')\leq{}w(\boldsymbol{u})+w(\boldsymbol{u}')$ for $\boldsymbol{u},\boldsymbol{u}'\in\mathrm{M}(f)$. This inequality is an equality if and only if $\boldsymbol{u}$ and $\boldsymbol{u}'$ are cofacial. On the other hand, by \cite{adolphson1989exponential}*{(1.11)}, the multiplication in $\overline{R}$ is given by
	\begin{equation*}
		\boldsymbol{x}^{\boldsymbol{u}}\cdot\boldsymbol{x}^{\boldsymbol{u}'}=\begin{cases}
			\boldsymbol{x}^{\boldsymbol{u}+\boldsymbol{u}'}
				&\mbox{if $\boldsymbol{u}$ and $\boldsymbol{u}'$ are cofacial}\rlap{,}\\
			0
				&\mbox{otherwise}\rlap{.}
		\end{cases}
	\end{equation*}
	Then, the assertion follows from the observation for $\pi\in\{\pi,\gamma\}$ and $\boldsymbol{u},\boldsymbol{u}'\in\mathrm{M}(f)$ that
	\begin{equation*}
		\pi^{w(\boldsymbol{u})}\boldsymbol{x}^{\boldsymbol{u}}\cdot\pi^{w(\boldsymbol{u}')}\boldsymbol{x}^{\boldsymbol{u}'}=\pi^{w(\boldsymbol{u})+w(\boldsymbol{u}')-w(\boldsymbol{u}+\boldsymbol{u}')}\cdot\pi^{w(\boldsymbol{u}+\boldsymbol{u}')}\boldsymbol{x}^{\boldsymbol{u}+\boldsymbol{u}'}\rlap{.}\qedhere
	\end{equation*}
\end{proof}

\begin{rmk}\label{rmk:reduction}
	\par Note that $R$ and $\overline{R}$ are identical as $k$-modules, but not as $k$-algebras. Regarding $\rho$ as a map from $B_{0}$ to $R$ does not give a homomorphism of rings.
\end{rmk}

\begin{prop}\label{prop:decomposition.lift}
	\par Assume that $f$ is nondegenerate and $\dim\Delta(f)=n$. Let $\mathrm{M}_{\mathrm{NP}}\subseteq\mathrm{M}(f)$ be the finite subset given by \cref{thm:decomposition.reduction}. Let $F\in\{\widehat{F},\widetilde{F}\}$.
	\begin{enumerate}
		\item If $i\neq{}n$, then $H^{i}(C^{\ast}(B,F))=0$.
		\item Set $V_{\mathrm{NP},m}=\langle\boldsymbol{x}^{\boldsymbol{u}}\mid\boldsymbol{u}\in\mathrm{M}_{\mathrm{NP}}\rangle_{K_{m}}$. Then, we have $B=V_{\mathrm{NP},m}\oplus\sum_{i=1}^{n}D_{F,i}B$.
	\end{enumerate}
\end{prop}
\begin{proof}
	\par If $F=\widetilde{F}$, then the assertion is given by the proof of \cite{adolphson1989exponential}*{Theorem~2.9}. Next, we consider the case where $F=\widehat{F}$. Note that for $i=1,\dots,n$, we have
	\begin{equation*}
		D_{\bar{f},i}\comp\rho=\rho\comp{}D_{\widehat{F},i}\rlap{,}
	\end{equation*}
	which yields a chain map $C^{\ast}(B_{0},\widehat{F})\rightarrow{}C^{\ast}(\overline{R},\bar{f})$. Note that $B_{0}$ is a \emph{flat separated complete} $O_{1}$-module in the sense of Monsky in \cite{monsky1970p}*{p.~91}, with respect to the norm on $B_{0}$ given by
	\begin{equation*}
		\lVert\xi\rVert=\sup\{\lvert{}A_{\boldsymbol{u}}\rvert_{p}\mid\boldsymbol{u}\in\mathrm{M}(f)\}
	\end{equation*}
	for $\xi(\boldsymbol{x})=\sum_{\boldsymbol{u}\in\mathrm{M}(f)}A_{\boldsymbol{u}}\boldsymbol{x}^{\boldsymbol{u}}\in{}B_{0}$. By \cref{thm:decomposition.reduction}.1 and \cite{monsky1970p}*{Theorem~8.5~(1)}, for $i\neq{}n$, we have $H^{i}(C^{\ast}(B,F))=0$. By \cref{thm:decomposition.reduction}.2 and the proof of \cite{adolphson1989exponential}*{Theorem~A.1}, we get
	\begin{equation*}
		B=V_{\mathrm{NP},m}\oplus\sum_{i=1}^{n}D_{F,i}B\rlap{.}\qedhere
	\end{equation*}
\end{proof}

\begin{prop}\label{prop:decomposition.growth}
	\par Assume that $f$ is nondegenerate and $\dim\Delta(f)=n$. Let $\frac{1}{p-1}\leq{}b\leq\frac{p}{p-1}$ and $c\in\mathbb{R}$. Set $V(b,c)=V_{\mathrm{NP},m}\cap{}L(b,c)$, and let $e=b-\frac{1}{p-1}$. Let $F\in\{\widehat{F},\widetilde{F}\}$. Then, we have
	\begin{equation*}
		L(b,c)=V(b,c)+\sum_{i=1}^{n}D_{F,i}L(b,c+e)\rlap{.}
	\end{equation*}
\end{prop}
\begin{proof}
	\par If $F=\widetilde{F}$, then the assertion is the same as \cite{adolphson1989exponential}*{Proposition~3.6}. It remains to show the assertion for $F=\widehat{F}$. We first consider the case where $b>\frac{1}{p-1}$. For $\xi(\boldsymbol{x})\in{}L(b,c)$, by \cite{adolphson1989exponential}*{Proposition~3.2}, there exists $v_{0}(\boldsymbol{x})\in{}V(b,c)$ and $\eta_{1,1}(\boldsymbol{x}),\dots,\eta_{n,1}(\boldsymbol{x})\in{}L(b,c+e)$, such that
	\begin{equation*}
		\xi(\boldsymbol{x})=v_{0}(\boldsymbol{x})+\sum_{i=1}^{n}\pi\hat{f}_{i}(\boldsymbol{x})\cdot\eta_{i,1}(\boldsymbol{x})\rlap{.}
	\end{equation*}
	Setting $\xi_{1}(\boldsymbol{x})=-\sum_{i=1}^{n}x_{i}\partial_{i}\eta_{i,1}(\boldsymbol{x})$, we have $\xi_{1}(\boldsymbol{x})\in{}L(b,c+e)$ and
	\begin{equation*}
		\xi(\boldsymbol{x})=v_{0}(\boldsymbol{x})+\xi_{1}(\boldsymbol{x})+\sum_{i=1}^{n}D_{\widehat{F},i}\eta_{i,1}(\boldsymbol{x})\rlap{.}
	\end{equation*}
	Applying the same argument for $\xi(\boldsymbol{x})$ as above to $\xi_{1}(\boldsymbol{x})$ and proceeding recursively, we obtain a sequence $\{(v_{l}(\boldsymbol{x}),\xi_{l}(\boldsymbol{x}),\eta_{1,l}(\boldsymbol{x}),\dots,\eta_{n,l}(\boldsymbol{x}))\}_{l=1}^{\infty}$ in $V(b,c+l\cdot{}e)\times{}L(b,c+l\cdot{}e)^{n+1}$ satisfying
	\begin{equation*}
		\xi_{l}(\boldsymbol{x})=v_{l}(\boldsymbol{x})+\xi_{l+1}(\boldsymbol{x})+\sum_{i=1}^{n}D_{\widehat{F},i}\eta_{i,l+1}(\boldsymbol{x})\rlap{.}
	\end{equation*}
	Since $b>\frac{1}{p-1}$, we have $e=b-\frac{1}{p-1}>0$, which implies that $\sum_{l=0}^{\infty}v_{l}(\boldsymbol{x})$ has a limit in $V(b,c)$, and that $\sum_{l=1}^{\infty}\eta_{i,l}(\boldsymbol{x})$ has a limit in $L(b,c)$ for $i=1,\dots,n$. let $v(\boldsymbol{x})=\sum_{l=0}^{\infty}v_{l}(\boldsymbol{x})$, and for $i=1,\dots,n$, let $\eta_{i}(\boldsymbol{x})=\sum_{l=1}^{\infty}\eta_{i,l}(\boldsymbol{x})$. Then, we have
	\begin{equation*}
		\xi(\boldsymbol{x})=v(\boldsymbol{x})+\sum_{i=1}^{n}D_{\widehat{F},i}\eta_{i}(\boldsymbol{x})\in{}V(b,c)+\sum_{i=1}^{n}D_{\widehat{F},i}L(b,c+e)\rlap{.}
	\end{equation*}
	
	\par Now, we consider the case where $b=\frac{1}{p-1}$. In this situation, we have $e=0$, so the argument above does not apply. We may assume $c=0$, since $L(b,c)=\varpi^{c}\cdot{}L(b,0)$, where we recall that $\varpi=p^{\frac{1}{p-1}}$. For $\xi(\boldsymbol{x})=\sum_{\boldsymbol{u}\in\mathrm{M}(f)}A_{\boldsymbol{u}}\boldsymbol{x}^{\boldsymbol{u}}\in{}L(b,0)$ and $l\in\mathbb{Z}_{\geq{}0}$, let
	\begin{equation*}
		\xi^{(l)}(\boldsymbol{x})=\sum_{w(\boldsymbol{u})\leq{}l}A_{\boldsymbol{u}}\boldsymbol{x}^{\boldsymbol{u}}\rlap{.}
	\end{equation*}
	Note that $\xi^{(l)}(\boldsymbol{x})\in{}B_{0}$ for $l\in\mathbb{Z}_{\geq{}0}$, so by the proof of \cref{prop:decomposition.lift}, there exist $v^{(l)}(\boldsymbol{x})\in{}V(\frac{1}{p-1},0)$ and $\eta^{(l)}_{1}(\boldsymbol{x}),\dots,\eta^{(l)}_{1}(\boldsymbol{x})\in{}B_{0}$, such that
	\begin{equation*}
		\xi^{(l)}(\boldsymbol{x})=v^{(l)}(\boldsymbol{x})+\sum_{i=1}^{n}D_{\widehat{F},i}\eta^{(l)}_{i}(\boldsymbol{x})\rlap{.}
	\end{equation*}
	Now, we get a sequence $\{(v^{(l)}(\boldsymbol{x}),\eta^{(l)}_{1}(\boldsymbol{x}),\dots,\eta^{(l)}_{1}(\boldsymbol{x}))\}_{l=0}^{\infty}$ in the space $V(\frac{1}{p-1},0)\times{}L(\frac{1}{p-1},0)^{n}$, which is compact in the topology of coefficient-wise convergence. Therefore, this sequence has a limit point $(v(\boldsymbol{x}),\eta_{1}(\boldsymbol{x}),\dots,\eta_{n}(\boldsymbol{x}))$ in $V(\frac{1}{p-1},0)\times{}L(\frac{1}{p-1},0)^{n}$ that satisfies
	\begin{equation*}
		\xi(\boldsymbol{x})=v(\boldsymbol{x})+\sum_{i=1}^{n}D_{\widehat{F},i}\eta_{i}(\boldsymbol{x})\rlap{.}\qedhere
	\end{equation*}
\end{proof}

\begin{cor}\label{cor:decomposition}
	\par Assume that $f$ is nondegenerate and $\dim\Delta(f)=n$. Let $\Lambda\in\mathfrak{A}^{\frac{p}{p-1}}$.
	\begin{enumerate}
		\item We have $\Lambda=V_{\mathrm{NP},m}\oplus\sum_{i=1}^{n}D_{\widetilde{F},i}\Lambda$.
		\item If $p\neq{}2$, then we have $\Lambda=V_{\mathrm{NP},m}\oplus\sum_{i=1}^{n}D_{\widehat{F},i}\Lambda$.
	\end{enumerate}
\end{cor}
\begin{proof}
	\par By \cref{prop:decomposition.lift}, the assertion is true if $\Lambda=B$. \cref{prop:decomposition.lift} also implies that the inclusion $V_{\mathrm{NP},m}\hookrightarrow{}B$ induces an isomorphism $V_{\mathrm{NP},m}\rightarrow{}H^{n}(C^{\ast}(B,F))$ for $F\in\{\widehat{F},\widetilde{F}\}$. We first consider the first assertion. By the proof of \cite{li2022exponential}*{Theorem~3.2}, we know that the map $\tilde{\iota}_{B}:H^{n}(C^{\ast}(B,\widetilde{F}))\rightarrow{}H^{n}(C^{\ast}(A^{\dagger}_{m},\widetilde{F}))$ induced by the inclusion $B\hookrightarrow{}A^{\dagger}_{m}$ is an isomorphism. Therefore, the inclusion $V_{\mathrm{NP},m}\hookrightarrow{}A^{\dagger}_{m}$ induces an isomorphism
	\begin{equation*}
		\tilde{\iota}_{\mathrm{NP},m}:V_{\mathrm{NP},m}\rightarrow{}H^{n}(C^{\ast}(A^{\dagger}_{m},\widetilde{F}))\rlap{.}
	\end{equation*}
	Hence, the first assertion is true for $\Lambda=A^{\dagger}_{m}$. For $B\subseteq\Lambda\subseteq{}A^{\dagger}_{m}$, namely $\Lambda\in\mathfrak{A}^{\frac{1}{p-1}}\cup\{B\}$, by an argument similar to that in the proof of \cite{li2022exponential}*{Proposition~3.1}, we can prove that the inclusion $B\hookrightarrow\Lambda$ induces a surjection $H^{n}(C^{\ast}(B,\widetilde{F}))\rightarrow{}H^{n}(C^{\ast}(\Lambda,\widetilde{F}))$, which implies the first assertion. For $\frac{1}{p-1}<b\leq\frac{p}{p-1}$, we have $\sum_{i=1}^{n}D_{\widetilde{F},i}L(b)\subseteq\sum_{i=1}^{n}D_{\widetilde{F},i}B$, since $L(b)\subseteq{}B$. This gives $V_{\mathrm{NP},m}\cap\sum_{i=1}^{n}D_{\widetilde{F},i}L(b)\subseteq{}V_{\mathrm{NP},m}\cap\sum_{i=1}^{n}D_{\widetilde{F},i}L(b)$. Since the first assertion is true $\Lambda=B$, we have $V_{\mathrm{NP},m}\cap\sum_{i=1}^{n}D_{\widetilde{F},i}L(b)=0$, which implies $V_{\mathrm{NP},m}\cap\sum_{i=1}^{n}D_{\widetilde{F},i}L(b)=0$. Therefore, \cref{prop:decomposition.growth} gives the first assertion for $\Lambda=L(b)$ with $\frac{1}{p-1}<b\leq\frac{p}{p-1}$. The argument above shows the first assertion.
	
	\par Next, we deal with the second assertion. By \cite{bourgeois1999annulation}*{Th\'eor\`em~1.4.(2)}, if $p\neq{}2$, then the map $\hat{\iota}_{B}:H^{n}(C^{\ast}(B,\widehat{F}))\rightarrow{}H^{n}(C^{\ast}(A^{\dagger}_{m},\widehat{F}))$ induced by the inclusion $B\hookrightarrow{}A^{\dagger}_{m}$ is an isomorphism. Therefore, if $p\neq{}2$, then the inclusion $V_{\mathrm{NP},m}\hookrightarrow{}A^{\dagger}_{m}$ induces an isomorphism
	\begin{equation*}
		\hat{\iota}_{\mathrm{NP},m}:V_{\mathrm{NP},m}\rightarrow{}H^{n}(C^{\ast}(A^{\dagger}_{m},\widehat{F}))\rlap{.}
	\end{equation*}
	Hence, the second assertion is true for $\Lambda=A^{\dagger}_{m}$. For $B\subseteq\Lambda\subseteq{}A^{\dagger}_{m}$, namely $\Lambda\in\mathfrak{A}^{\frac{1}{p-1}}\cup\{B\}$, by an argument similar to that in the proof of \cite{bourgeois1999annulation}*{Corollary~1.3}, we can prove that the inclusion $B\hookrightarrow\Lambda$ induces an surjection $H^{n}(C^{\ast}(B,\widehat{F}))\rightarrow{}H^{n}(C^{\ast}(\Lambda,\widehat{F}))$, which implies the first assertion. For $\frac{1}{p-1}<b\leq\frac{p}{p-1}$, we have $\sum_{i=1}^{n}D_{\widehat{F},i}L(b)\subseteq\sum_{i=1}^{n}D_{\widehat{F},i}B$, since $L(b)\subseteq{}B$. This gives $V_{\mathrm{NP},m}\cap\sum_{i=1}^{n}D_{\widehat{F},i}L(b)\subseteq{}V_{\mathrm{NP},m}\cap\sum_{i=1}^{n}D_{\widehat{F},i}L(b)$. The second assertion for $\Lambda=B$ gives $V_{\mathrm{NP},m}\cap\sum_{i=1}^{n}D_{\widehat{F},i}L(b)=0$, so $V_{\mathrm{NP},m}\cap\sum_{i=1}^{n}D_{\widehat{F},i}L(b)=0$. Then, by \cref{prop:decomposition.growth}, the second assertion is true for $\Lambda=L(b)$ with $\frac{1}{p-1}<b\leq\frac{p}{p-1}$. The argument above shows the second assertion.
\end{proof}

\begin{proof}[Proof of \cref{thm:object.cohomology}]
	\par Let $V_{\mathrm{NP}}=\langle\boldsymbol{x}^{\boldsymbol{u}}\mid\boldsymbol{u}\in\mathrm{M}_{\mathrm{NP}}\rangle_{K_{1}}$, then we have $V_{\mathrm{NP},m}=K_{m}\otimes_{K_{1}}V_{\mathrm{NP}}$. Since $p\neq{}2$, by \cref{cor:decomposition}.2, we have $A^{\dagger}_{m}=V_{\mathrm{NP},m}\oplus\sum_{i=1}^{n}D_{\widehat{F},i}A^{\dagger}_{m}$. It is straightforward to verify that $D_{\widehat{F},i}A^{\dagger}_{m}=K_{m}\otimes_{K_{1}}D_{\widehat{F},i}A^{\dagger}$ for $i=1,\dots,n$. Thus $A^{\dagger}=V_{\mathrm{NP}}\oplus\sum_{i=1}^{n}D_{\widehat{F},i}A^{\dagger}$. This implies that the inclusion $V_{\mathrm{NP}}\hookrightarrow{}A^{\dagger}$ induces an isomorphism
	\begin{equation*}
		\hat{\iota}_{\mathrm{NP}}:V_{\mathrm{NP}}\rightarrow{}H^{n}(C^{\ast}(A^{\dagger},\widehat{F}))\rlap{.}
	\end{equation*}
	For $l\in\mathbb{Z}$, we consider a map $\mathrm{DR}^{l}(A,\nabla_{\widehat{F}})\rightarrow{}C^{l}(A^{\dagger},\widehat{F})$ given by
	\begin{equation*}
		\xi(\boldsymbol{x})dx_{i_{1}}\wedge\dots\wedge{}dx_{i_{l}}\mapsto{}x_{i_{1}}\dots{}x_{i_{l}}\xi(\boldsymbol{x})\frac{dx_{i_{1}}}{x_{i_{1}}}\wedge\dots\wedge\frac{dx_{i_{i_{l}}}}{x_{i_{l}}}
	\end{equation*}
	for $\xi(\boldsymbol{x})\in{}A$ and $1\leq{}i_{1}<\dots<i_{l}\leq{}n$. It is straightforward to verify that this map is compatible with the differentials, so that we obtain a chain map $\mathrm{DR}^{\ast}(A,\nabla_{\widehat{F}})\rightarrow{}C^{\ast}(A^{\dagger},\widehat{F})$. This chain map induces a morphism $H^{n}(\mathrm{DR}^{\ast}(A,\nabla_{\widehat{F}}))\rightarrow{}H^{n}(C^{\ast}(A^{\dagger},\widehat{F}))$ in $\mathbf{Mod}_{K_{1}}$. Note that the specialization map $\iota_{\widehat{F}}$ is the composition of $\hat{\imath}$ and this morphism, where $\hat{\imath}:H^{n}(C^{\ast}(A^{\dagger},\widehat{F}))\rightarrow{}H^{n}_{\mathrm{rig}}(\mathbb{T}^{n}_{k}/K_{1},f^{\ast}\mathcal{L}_{\pi})$ is the isomorphism given by \cref{eq:logarithm}. To show that $\iota_{\widehat{F}}$ is an isomorphism, it suffices to show that the morphism $H^{n}(\mathrm{DR}^{\ast}(A,\nabla_{\widehat{F}}))\rightarrow{}H^{n}(C^{\ast}(A^{\dagger},\widehat{F}))$ is an isomorphism. Since $V_{\mathrm{NP}}\subseteq{}A$, the morphism is surjective. On the other hand, by \cite{adolphson1997twisted}*{Theorem~1.4}, we have $\dim_{K_{1}}H^{n}(\mathrm{DR}^{\ast}(A,\nabla_{\widehat{F}}))=n!\mathrm{vol}(\widehat{F})$. Since $\Delta(\widehat{F})=\Delta(f)$, we have $\mathrm{vol}(\widehat{F})=\mathrm{vol}(f)$. Therefore, by \cref{thm:decomposition.reduction}.2 and \cref{cor:decomposition}, we get
	\begin{equation*}
		\dim_{K_{1}}H^{n}(\mathrm{DR}^{\ast}(A,\nabla_{\widehat{F}}))=\dim_{K_{1}}H^{n}(C^{\ast}(A^{\dagger},\widehat{F}))\rlap{.}
	\end{equation*}
	This implies that the surjective morphism $H^{n}(\mathrm{DR}^{\ast}(A,\nabla_{\widehat{F}}))\rightarrow{}H^{n}(C^{\ast}(A^{\dagger},\widehat{F}))$ is an isomorphism, and hence the $\iota_{\widehat{F}}$ is an isomorphism.
\end{proof}

\begin{rmk}\label{rmk:inclusion}
	\par Assume that $f$ is nondegenerate and $\dim\Delta(f)=n$. By \cref{cor:decomposition}.1, we have $A^{\dagger}_{m}=V_{\mathrm{NP},m}\oplus\sum_{i=1}^{n}D_{\widetilde{F},i}A^{\dagger}_{m}$. Note that for $i=1,\dots,n$, we have $D_{\widetilde{F},i}A^{\dagger}_{m}=K_{m}\otimes_{K_{1}}D_{\widetilde{F},i}A^{\dagger}$. Thus $A^{\dagger}=V_{\mathrm{NP}}\oplus\sum_{i=1}^{n}D_{\widetilde{F},i}A^{\dagger}$, which implies that the inclusion $V_{\mathrm{NP}}\hookrightarrow{}A^{\dagger}$ induces an isomorphism
	\begin{equation*}
		\tilde{\iota}_{\mathrm{NP}}:V_{\mathrm{NP}}\rightarrow{}H^{n}(C^{\ast}(A^{\dagger},\widetilde{F}))\rlap{.}
	\end{equation*}
\end{rmk}

\par Next, we construct two filtered $\Phi$-modules which are useful for the proof of \cref{thm:main}. For $F\in\{\widehat{F},\widetilde{F}\}$, let $\phi_{F}:A^{\dagger}\rightarrow{}A^{\dagger}$ be the endomorphism given by
\begin{equation*}
	\xi(\boldsymbol{x})\mapsto\exp(\varphi(F(\boldsymbol{x}))-F(\boldsymbol{x}))\cdot\varphi(\xi(\boldsymbol{x}))\rlap{.}
\end{equation*}
It is straightforward to verify that $D_{F,i}\comp\phi_{F}=p\cdot(\phi_{F}\comp{}D_{F,i})$ for $F\in\{\widehat{F},\widetilde{F}\}$ and $i=1,\dots,n$. This yields a chain map $\phi_{F}:C^{\ast}(A^{\dagger},F)\rightarrow{}C^{\ast}(A^{\dagger},F)$. For $i\in\mathbb{Z}$, let $\phi_{F}:H^{i}(C^{\ast}(A^{\dagger},F))\rightarrow{}H^{i}(C^{\ast}(A^{\dagger},F))$ denote the induced endomorphism. We note that $\phi_{F}$ is $\sigma$-semilinear and bijective.

\begin{dfn}\label{dfn:np-module}
	\par Assume that $f$ is nondegenerate and $\dim\Delta(f)=n$. Let $F\in\{\widehat{F},\widetilde{F}\}$.
	\begin{enumerate}
		\item Let $\phi_{\mathrm{NP}}=\iota^{-1}_{\mathrm{NP}}\comp\phi_{F}\comp\iota_{\mathrm{NP}}$, which is a bijective $\sigma$-semilinear endomorphism on $V_{\mathrm{NP}}$. Here, we write $\phi_{\mathrm{NP}},\iota_{\mathrm{NP}}$ as $\hat{\phi}_{\mathrm{NP}},\hat{\iota}_{\mathrm{NP}}$ if $F=\widehat{F}$, and as $\tilde{\phi}_{\mathrm{NP}},\tilde{\iota}_{\mathrm{NP}}$ if $F=\widetilde{F}$.
		\item Define an exhaustive separated filtration $F^{\ast}_{\mathrm{NP}}$, which we call the \emph{Newton polyhedron filtration}, on $V_{\mathrm{NP}}$ by setting $F^{i}_{\mathrm{NP}}V_{\mathrm{NP}}=\langle\boldsymbol{x}^{\boldsymbol{u}}\mid\boldsymbol{u}\in\mathrm{M}_{\mathrm{NP}}\mbox{, }w(\boldsymbol{u})\leq{}-i\rangle_{K_{1}}$ for $i\in\mathbb{R}$.
	\end{enumerate}
	We call the filtered $\Phi$-module $(V_{\mathrm{NP}},\phi_{\mathrm{NP}},F^{\ast}_{\mathrm{NP}})\in\mathbf{M\overline{F}}{}^{\Phi}_{K_{1}/K_{1}}$ the \emph{Newton polyhedron module} associated with $F$.
\end{dfn}

\begin{prop}\label{prop:np-module}
	\par Assume that $p\neq{}2$ and $f,\hat{f}$ are nondegenerate with $\dim\Delta(f)=n$. Then
	\begin{equation*}
		((H^{n}_{\mathrm{rig}}(\mathbb{T}^{n}_{k}/K_{1},f^{\ast}\mathcal{L}_{\pi}),\phi_{f}),(H^{n}_{\mathrm{dR}}(\mathbb{T}^{n}_{K_{1}},\nabla_{\widehat{F}}),F^{\ast}_{\mathrm{irr}}),\iota_{\widehat{F}})\cong(V_{\mathrm{NP}},\hat{\phi}_{\mathrm{NP}},F^{\ast}_{\mathrm{NP}})\in\mathbf{M\overline{F}}{}^{\Phi}_{K_{1}/K_{1}}\rlap{.}
	\end{equation*}
\end{prop}
\begin{proof}
	\par Let $\iota_{\mathrm{rig}}=\hat{\imath}\comp\hat{\iota}_{\mathrm{NP}}$. Note that $\iota_{\mathrm{rig}}:V_{\mathrm{NP}}\rightarrow{}H^{n}_{\mathrm{rig}}(\mathbb{T}^{n}_{k}/K_{1},f^{\ast}\mathcal{L}_{\pi})$ is an isomorphism. It is straightforward to verify that $\hat{\imath}\comp\phi_{\widehat{F}}=\phi_{f}\comp\hat{\imath}$, which implies $\iota_{\mathrm{rig}}\comp\hat{\phi}_{\mathrm{NP}}=\phi_{f}\comp\iota_{\mathrm{rig}}$. Thus $\iota_{\mathrm{rig}}$ gives an isomorphism
	\begin{equation*}
		\iota_{\mathrm{rig}}:(V_{\mathrm{NP}},\hat{\phi}_{\mathrm{NP}})\rightarrow(H^{n}_{\mathrm{rig}}(\mathbb{T}^{n}_{k}/K_{1},f^{\ast}\mathcal{L}_{\pi}),\phi_{f})\in\mathbf{Mod}{}^{\Phi}_{K_{1}}\rlap{.}
	\end{equation*}
	On the other hand, let $\iota_{\mathrm{dR}}=\iota_{\mathrm{rig}}^{-1}\comp\iota_{\widehat{F}}$. Note that $\iota_{\mathrm{dR}}:H^{n}_{\mathrm{dR}}(\mathbb{T}^{n}_{K_{1}},\nabla_{\widehat{F}})\rightarrow{}V_{\mathrm{NP}}$ is an isomorphism. By \cite{yu2014irregular}*{p.~126~footnote}, we have $\iota_{\mathrm{dR}}(F^{\ast}_{\mathrm{irr}})=F^{\ast}_{\mathrm{NP}}$. This implies that $\iota_{\mathrm{dR}}$ gives an isomorphism
	\begin{equation*}
		\iota_{\mathrm{dR}}:(H^{n}_{\mathrm{dR}}(\mathbb{T}^{n}_{K_{1}},\nabla_{\widehat{F}}),F^{\ast}_{\mathrm{irr}})\rightarrow(V_{\mathrm{NP}},F^{\ast}_{\mathrm{NP}})\in\mathbf{M\overline{F}}_{K_{1}}\rlap{.}
	\end{equation*}
	Thus $\boldsymbol{\iota}:((H^{n}_{\mathrm{rig}}(\mathbb{T}^{n}_{k}/K_{1},f^{\ast}\mathcal{L}_{\pi}),\phi_{f}),(H^{n}_{\mathrm{dR}}(\mathbb{T}^{n}_{K_{1}},\nabla_{\widehat{F}}),F^{\ast}_{\mathrm{irr}}),\iota_{\widehat{F}})\rightarrow(V_{\mathrm{NP}},\hat{\phi}_{\mathrm{NP}},F^{\ast}_{\mathrm{NP}})$ is an isomorphism in $\mathbf{M\overline{F}}{}^{\Phi}_{K_{1}/K_{1}}$, where $\boldsymbol{\iota}=(\iota_{\mathrm{rig}}^{-1},\iota_{\mathrm{dR}})$.
\end{proof}

\begin{rmk}\label{rmk:np-module}
	\par Assume that $p\neq{}2$ and $f,\hat{f}$ are nondegenerate with $\dim\Delta(f)=n$. By \cref{prop:np-module}, we have $\ell_{\mathrm{HT}}(H^{n}_{\mathrm{dR}}(\mathbb{T}^{n}_{K_{1}},\nabla_{\widehat{F}}),F^{\ast}_{\mathrm{irr}})=\ell_{\mathrm{HT}}(V_{\mathrm{NP}},F^{\ast}_{\mathrm{NP}})$. By \cref{dfn:np-module}, we have
	\begin{equation*}
		\ell_{\mathrm{HT}}(V_{\mathrm{NP}},F^{\ast}_{\mathrm{NP}})=\max\{w(\boldsymbol{u})\mid\boldsymbol{u}\in\mathrm{M}_{\mathrm{NP}}\}\rlap{,}
	\end{equation*}
	since $\min\{w(\boldsymbol{u})\mid\boldsymbol{u}\in\mathrm{M}_{\mathrm{NP}}\}=0$. Therefore, \cref{rmk:decomposition.reduction} provides an algorighm to compute the Hodge-Tate length $\ell_{\mathrm{HT}}(H^{n}_{\mathrm{dR}}(\mathbb{T}^{n}_{K_{1}},\nabla_{\widehat{F}}),F^{\ast}_{\mathrm{irr}})$ combinatorially.
\end{rmk}

\par We consider the relationship between $(V_{\mathrm{NP}},\hat{\phi}_{\mathrm{NP}},F^{\ast}_{\mathrm{NP}})$ and $(V_{\mathrm{NP}},\tilde{\phi}_{\mathrm{NP}},F^{\ast}_{\mathrm{NP}})$. For $\Lambda\in\mathfrak{A}^{\frac{p-1}{p}}$, let $T_{f}:\Lambda\rightarrow\Lambda$ be the endomorphism given by
\begin{equation*}
	\xi(\boldsymbol{x})\mapsto\exp(\widetilde{F}(\boldsymbol{x})-\widehat{F}(\boldsymbol{x}))\cdot\xi(\boldsymbol{x})\rlap{,}
\end{equation*}
which is well-defined by \cref{lem:valuation.transformation}. It is straightforward to verify that $T_{f}$ is an isomorphism. For $i=1,\dots,n$, Note that $D_{\widehat{F},i}\comp{}T_{f}=T_{f}\comp{}D_{\widetilde{F},i}$. This yields an isomorphism $T_{f}:C^{\ast}(\Lambda,\widetilde{F})\rightarrow{}C^{\ast}(\Lambda,\widehat{F})$ of chain complexes. For $i\in\mathbb{Z}$, we denote by $T_{f}:H^{i}(C^{\ast}(\Lambda,\widetilde{F}))\rightarrow{}H^{i}(C^{\ast}(\Lambda,\widehat{F}))$ the induced map, which is an isomorphism in $\mathbf{Mod}_{K_{1}}$. Furthermore, it is straightforward to verify that $\phi_{\widehat{F}}\comp{}T_{f}=T_{f}\comp\phi_{\widetilde{F}}$, which implies that $T_{f}$ gives an isomorphism
\begin{equation*}
	T_{f}:(H^{i}(C^{\ast}(\Lambda,\widetilde{F})),\phi_{\widetilde{F}})\rightarrow(H^{i}(C^{\ast}(\Lambda,\widehat{F})),\phi_{\widehat{F}})\in\mathbf{Mod}{}^{\Phi}_{K_{1}}\rlap{.}
\end{equation*}

\begin{lem}\label{lem:valuation.uniformizer}
	\par We have $\operatorname{ord}(\gamma-\pi)=p-1+\frac{1}{p-1}$.
\end{lem}
\begin{proof}
	\par Let $\delta=\frac{\gamma}{\pi}$. It suffice to show that $\operatorname{ord}(\delta-1)=p-1$, since $\operatorname{ord}\pi=\frac{1}{p-1}$. Recall that
	\begin{equation}\label{eq:congruence}
		\gamma\equiv\pi\equiv\zeta_{p}-1\bmod(\zeta_{p}-1)^{2}\rlap{,}
	\end{equation}
	which implies that $\operatorname{ord}\gamma=\frac{1}{p-1}$. Thus $\operatorname{ord}\delta=0$. \cref{eq:congruence} also implies that $\operatorname{ord}(\delta-1)\geq\frac{1}{p-1}$, since it gives $\gamma-\pi\equiv{}0\bmod(\zeta_{p}-1)^{2}$. Recall that $\pi+\frac{\pi^{p}}{p}=0$ and $\sum_{i=0}^{\infty}\frac{\gamma^{p^{i}}}{p^{i}}=0$. Thus
	\begin{equation*}
		1-\delta^{p-1}=1+\frac{\gamma^{p-1}}{p}=\gamma^{-1}\left(\gamma+\frac{\gamma^{p}}{p}\right)=\gamma^{-1}\left(-\sum_{i=2}^{\infty}\frac{\gamma^{p^{i}}}{p^{i}}\right)\rlap{.}
	\end{equation*}
	Note that $\operatorname{ord}(\sum_{i=2}^{\infty}\frac{\gamma^{p^{i}}}{p^{i}})=\operatorname{ord}\frac{\gamma^{p^{2}}}{p^{2}}=\frac{p^{2}}{p-1}-2$. Therefore, we have
	\begin{equation*}
		\operatorname{ord}(1-\delta^{p-1})=-\operatorname{ord}\gamma+\operatorname{ord}\left(\sum_{i=2}^{\infty}\frac{\gamma^{p^{i}}}{p^{i}}\right)=-\frac{1}{p-1}+\frac{p^{2}}{p-1}-2=p-1\rlap{.}
	\end{equation*}
	If $p=2$, then this shows $\operatorname{ord}(1-\delta)=p-1$. Next, we deal with the case where $p\neq{}2$. Note that
	\begin{equation*}
		1-\delta^{p-1}=(1-\delta)(1+\dots+\delta^{p-2})\rlap{.}
	\end{equation*}
	By \cref{eq:congruence}, we have $\delta\equiv{}1\bmod\frac{(\zeta_{p}-1)^{2}}{\pi}$, which gives $1+\dots+\delta^{p-2}\equiv{}p-1\equiv{}-1\bmod\frac{(\zeta_{p}-1)^{2}}{\pi}$. Thus $\operatorname{ord}(1+\dots+\delta^{p-2})=0$, which implies that $\operatorname{ord}(1-\delta)=\operatorname{ord}(1-\delta^{p-1})=p-1$.
\end{proof}

\begin{lem}\label{lem:valuation.transformation}
	\par We have $\exp(\widetilde{F}(\boldsymbol{x})-\widehat{F}(\boldsymbol{x}))\in{}L(\frac{p-1}{p},0)$.
\end{lem}
\begin{proof}
	\par Note that $\exp(\widetilde{F}(\boldsymbol{x})-\widehat{F}(\boldsymbol{x}))=\prod_{\boldsymbol{u}\in\operatorname{supp}f}\exp((\gamma-\pi)\hat{\alpha}_{\boldsymbol{u}}\boldsymbol{x}^{\boldsymbol{u}})\cdot\exp(\sum_{l=1}^{\infty}\gamma_{l}\hat{\alpha}_{\boldsymbol{u}}^{p^{l}}\boldsymbol{x}^{p^{l}\boldsymbol{u}})$. By the proof of \cite{li2022exponential}*{Proposition~2.1}, we have $\exp(\sum_{l=1}^{\infty}\gamma_{l}\hat{\alpha}_{\boldsymbol{u}}^{p^{l}}\boldsymbol{x}^{p^{l}\boldsymbol{u}})\in{}L(\frac{p-1}{p},0)$ for $\boldsymbol{u}\in\operatorname{supp}f$. Next, we show that $\exp((\gamma-\pi)\hat{f}(\boldsymbol{x}))\in{}L(p-1,0)$. For $\boldsymbol{u}\in\operatorname{supp}f$, we have
	\begin{equation*}
		\exp((\gamma-\pi)\hat{\alpha}_{\boldsymbol{u}}\boldsymbol{x}^{\boldsymbol{u}})=\sum_{i=0}^{\infty}\frac{(\gamma-\pi)^{i}\hat{\alpha}_{\boldsymbol{u}}^{i}}{i!}\cdot\boldsymbol{x}^{i\boldsymbol{u}}\rlap{.}
	\end{equation*}
	Note that $\operatorname{ord}(i!)=\sum_{j=1}^{\infty}\left\lfloor\frac{i}{p^{j}}\right\rfloor\leq{}\sum_{j=1}^{\infty}\frac{i}{p^{j}}=\frac{i}{p-1}$, where $\lfloor\lambda\rfloor$ denotes the maximal integer that is less or equal to $\lambda\in\mathbb{R}$. Therefore, for $\boldsymbol{u}\in\operatorname{supp}f$ and $i\in\mathbb{Z}_{\geq{}0}$, we have
	\begin{align*}
		\operatorname{ord}\left(\frac{(\gamma-\pi)^{i}\hat{\alpha}_{\boldsymbol{u}}^{i}}{i!}\right)={}&i\operatorname{ord}(\gamma-\pi)-\operatorname{ord}(i!)\geq{}i\cdot\left(p-1+\frac{1}{p-1}\right)-\frac{i}{p-1}\\
		={}&i\cdot(p-1)\geq(p-1)\cdot{}w(i\boldsymbol{u})\rlap{.}
	\end{align*}
	This shows that $\exp((\gamma-\pi)\hat{\alpha}_{\boldsymbol{u}}\boldsymbol{x}^{\boldsymbol{u}})\in{}L(p-1,0)$ for $\boldsymbol{u}\in\operatorname{supp}f$. Therefore, the assertion follows, since we have $\frac{p-1}{p}<1\leq{}p-1$.
\end{proof}

\begin{prop}\label{prop:transformation}
	\par Assume that $p\neq{}2$ and $f$ is nondegenerate with $\dim\Delta(f)=n$. Then
	\begin{equation*}
		(V_{\mathrm{NP}},\tilde{\phi}_{\mathrm{NP}})\cong(V_{\mathrm{NP}},\hat{\phi}_{\mathrm{NP}})\in\mathbf{Mod}{}^{\Phi}_{K_{1}}\rlap{.}
	\end{equation*}
\end{prop}
\begin{proof}
	\par Let $T_{\mathrm{NP}}=\hat{\iota}^{-1}_{\mathrm{NP}}\comp{}T_{f}\comp\tilde{\iota}_{\mathrm{NP}}$, which is an automorphism on $V_{\mathrm{NP}}$. It is straightforward to verify that $\hat{\phi}_{\mathrm{NP}}\comp{}T_{\mathrm{NP}}=T_{\mathrm{NP}}\comp\tilde{\phi}_{\mathrm{NP}}$, which implies that $T_{\mathrm{NP}}$ gives an isomorphism
	\begin{equation*}
		T_{\mathrm{NP}}:(V_{\mathrm{NP}},\tilde{\phi}_{\mathrm{NP}})\rightarrow(V_{\mathrm{NP}},\hat{\phi}_{\mathrm{NP}})\in\mathbf{Mod}{}^{\Phi}_{K_{1}}\rlap{.}\qedhere
	\end{equation*}
\end{proof}

\begin{rmk}\label{rmk:transformation}
	\par Assume that $p\neq{}2$ and $f$ is nondegenerate with $\dim\Delta(f)=n$. We note that $T_{\mathrm{NP}}$ cannot be extended to an isomorphism in $\mathbf{M\overline{F}}{}^{\Phi}_{K_{1}/K_{1}}$ from $(V_{\mathrm{NP}},\tilde{\phi}_{\mathrm{NP}},F^{\ast}_{\mathrm{NP}})$ to $(V_{\mathrm{NP}},\hat{\phi}_{\mathrm{NP}},F^{\ast}_{\mathrm{NP}})$ in general. If this were true, we would have \cref{conj:main} right away. See \cref{eg:counterexample} for an example where $T_{\mathrm{NP}}$ is not even filtration-compatible with respect to $F^{\ast}_{\mathrm{NP}}$.
\end{rmk}

\section{Weak admissibility}
\label{section:proof}

\par We retain the notations from \cref{section:np-module}. In this section, we prove \cref{thm:main}. There are two steps: the first step is to prove that $(V_{\mathrm{NP}},\tilde{\phi}_{\mathrm{NP}},F^{\ast}_{\mathrm{NP}})$ is weakly admissible (cf. \cref{thm:weak-admissibility}), and the second step is to show under certain conditions that the automorphism $T_{\mathrm{NP}}:V_{\mathrm{NP}}\rightarrow{}V_{\mathrm{NP}}$ preserves being weakly admissible.

\par In this section, assume that $f$ is nondegenerate and $\dim\Delta(f)$. For $v\in{}V_{\mathrm{NP},m}$, we may write
\begin{equation*}
	v=\sum_{\boldsymbol{u}\in\mathrm{M}_{\mathrm{NP}}}A_{\boldsymbol{u}}(v)\cdot\pi^{pw(\boldsymbol{u})}\boldsymbol{x}^{\boldsymbol{u}}
\end{equation*}
with unique $A_{\boldsymbol{u}}(v)\in{}K_{m}$. We define the \emph{order} of $v$ to be $\operatorname{ord}v=\min\{\operatorname{ord}A_{\boldsymbol{u}}(v)\mid\boldsymbol{u}\in\mathrm{M}_{\mathrm{NP}}\}$. Define the \emph{weight} of $v$ to be $w(v)=\max\{w(\boldsymbol{u})\mid\boldsymbol{u}\in\mathrm{M}_{\mathrm{NP}}\mbox{, }A_{\boldsymbol{u}}(v)\neq{}0\}$. We note that
\begin{equation*}
	w(v)=-w_{\mathrm{HT}}(v)\rlap{.}
\end{equation*}
Here $w_{\mathrm{HT}}(v)$ is the Hodge-Tate weight of $v$ with respect to $F^{\ast}_{\mathrm{NP}}$ defined in \cref{dfn:hodge-tate.weight}. Let $d_{\mathrm{NP}}=\dim_{K_{1}}V_{\mathrm{NP}}$. Sort $\mathrm{M}_{\mathrm{NP}}$ and write $\mathrm{M}_{\mathrm{NP}}=\{\boldsymbol{u}_{1},\dots,\boldsymbol{u}_{d_{\mathrm{NP}}}\}$ so that
\begin{equation*}
	w(\boldsymbol{u}_{1})\geq\dots\geq{}w(\boldsymbol{u}_{d_{\mathrm{NP}}})
\end{equation*}
Note that $\ell_{\mathrm{HT}}(V_{\mathrm{NP}},F^{\ast}_{\mathrm{NP}})=w(\boldsymbol{u}_{1})$. Define an order on $\mathrm{M}_{\mathrm{NP}}$ by setting $\boldsymbol{u}_{d_{\mathrm{NP}}}\prec\dots\prec\boldsymbol{u}_{1}$. Then, for $v\in{}V_{\mathrm{NP}}$, we define its \emph{leading power} to be
\begin{equation*}
	\boldsymbol{\mu}(v)=\max\{\boldsymbol{u}\mid\boldsymbol{u}\in\mathrm{M}_{\mathrm{NP}}\mbox{, }\operatorname{ord}v=\operatorname{ord}A_{\boldsymbol{u}}(v)\}\rlap{.}
\end{equation*}
Note that $w(v)\geq{}w(\boldsymbol{\mu}(v))$.

\begin{dfn}\label{dfn:np-basis}
	\par Assume that $f$ is nondegenerate and $\dim\Delta(f)$. Let $V$ be a $K_{m}$-subspace of $V_{\mathrm{NP},m}$, and let $d=\dim_{K_{m}}V$. A basis $\{v_{i}\}_{i=1}^{d}$ of $V$ is called a \emph{quasi-NP basis}, if it is filtration-generating with respect to $F^{\ast}_{\mathrm{NP}}$, and satisfies
	\begin{equation*}
		A_{\boldsymbol{\mu}(v_{i})}(v_{j})=\begin{cases}
			1
				&i=j\rlap{,}\\
			0
				&i\neq{}j\mbox{ and }w(v_{i})\leq{}w(v_{j})\rlap{,}\\
			\mbox{anything}
				&\mbox{otherwise}\rlap{.}
		\end{cases}
	\end{equation*}
\end{dfn}

\begin{lem}\label{lem:np-basis}
	\par Assume that $f$ is nondegenerate and $\dim\Delta(f)$. Any subobject of $(V_{\mathrm{NP},m},F^{\ast}_{\mathrm{NP}})$ in $\mathbf{M\overline{F}}_{K_{m}}$ admits a quasi-NP basis.
\end{lem}
\begin{proof}
	\par Let $V\subseteq{}V_{\mathrm{NP},m}$ be a subobject of $(V_{\mathrm{NP},m},F^{\ast}_{\mathrm{NP}})$, and let $d=\dim_{K_{m}}V$. We show by induction on $d$ that $V$ admits a quasi-NP basis. If $d=1$, pick $v\in{}V\setminus\{0\}$, and let
	\begin{equation*}
		v_{1}=A_{\boldsymbol{\mu}(v)}(v)^{-1}\cdot{}v\rlap{.}
	\end{equation*}
	Note that $\boldsymbol{\mu}(v_{1})=\boldsymbol{\mu}(v)$, so that $A_{\boldsymbol{\mu}(v_{1})}(v_{1})=A_{\boldsymbol{\mu}(v)}(v)^{-1}\cdot{}A_{\boldsymbol{\mu}(v_{1})}(v)=1$. Thus $\{v_{1}\}$ is a quasi-NP basis of $V$.
	
	\par Next, Assuming the assertion for $1\leq{}d\leq{}r-1$, we show that the assertion is true for $d=r$. Sort $\mathfrak{W}_{\mathrm{HT}}(V,F^{\ast}_{\mathrm{NP}})$ and write $\mathfrak{W}_{\mathrm{HT}}(V,F^{\ast}_{\mathrm{NP}})=\{i_{1},\dots,i_{l}\}$ such that $i_{1}>\dots>i_{l}$. For $h=1,\dots,l$, let $d_{h}=\dim_{K_{m}}F^{i_{h}}V$. Pick an $(r-1)$-dimensional $K_{m}$-subspace $W$ of $V$ such that $F^{>i_{l}}V\subseteq{}W$. By the inductive hypothesis, there exists a quasi-NP basis $\{w_{i}\}_{i=1}^{r-1}$ of $W$. Extend $\{w_{i}\}_{i=1}^{r-1}$ to a basis $\{w_{i}\}_{i=1}^{r-1}\cup\{v\}$ of $V$. We construct a sequence $\{v^{(1)},\dots,v^{(l)}\}$ in $V$. We set the first term to be
	\begin{equation*}
		v^{(1)}=v-\sum_{i=1}^{d_{1}}A_{\boldsymbol{\mu}(w_{i})}(v)\cdot{}w_{i}\rlap{.}
	\end{equation*}
	For $h=1,\dots,l-2$, suppose $v^{(h)}$ is already given, we define the next term of the sequence by
	\begin{equation*}
		v^{(h+1)}=v^{(h)}-\sum_{i=d_{h}+1}^{d_{h+1}}A_{\boldsymbol{\mu}(w_{i})}(v^{(h)})\cdot{}w_{i}\rlap{.}
	\end{equation*}
	By this, we obtain $v^{(1)},\dots,v^{(l-1)}$. Now, we define the last term of the sequence to be
	\begin{equation*}
		v^{(l)}=v^{(l-1)}-\sum_{i=d_{l-1}+1}^{r-1}A_{\boldsymbol{\mu}(w_{i})}(v)\cdot{}w_{i}\rlap{.}
	\end{equation*}
	Let $v_{r}=A_{\boldsymbol{\mu}(v^{(l)})}(v^{(l)})^{-1}\cdot{}v^{(l)}$, then we have $\boldsymbol{\mu}(v_{r})=\boldsymbol{\mu}(v^{(l)})$ and $A_{\boldsymbol{\mu}(v_{r})}(v_{r})=1$. We now show that $A_{\boldsymbol{\mu}(w_{i})}(v_{r})=0$ for $i=1,\dots,r-1$. Since $\{w_{i}\}_{i=1}^{r-1}$ is a quasi-NP basis, for $i,j=1,\dots,d_{1}$, we have $A_{\boldsymbol{\mu}(w_{i})}(w_{j})=1$ if $i=j$, and $A_{\boldsymbol{\mu}(w_{i})}(w_{j})=0$ if $i\neq{}j$. Thus
	\begin{align*}
		A_{\boldsymbol{\mu}(w_{i})}(v^{(1)})={}&A_{\boldsymbol{\mu}(w_{i})}(v)-\sum_{j=1}^{d_{1}}A_{\boldsymbol{\mu}(w_{j})}(v)\cdot{}A_{\boldsymbol{\mu}(w_{i})}(w_{j})\\
		={}&A_{\boldsymbol{\mu}(w_{i})}(v)-A_{\boldsymbol{\mu}(w_{i})}(v)\cdot{}1=0\rlap{.}
	\end{align*}
	For $h=1,\dots,l-2$, suppose that $A_{\boldsymbol{\mu}(w_{i})}(v^{(h)})=0$ for $i=1,\dots,d_{h}$, we show $A_{\boldsymbol{\mu}(w_{i})}(v^{(h+1)})=0$ for $i=1,\dots,d_{h+1}$. Since $\{w_{i}\}_{i=1}^{r-1}$ is a quasi-NP basis, we have $A_{\boldsymbol{\mu}(w_{i})}(w_{j})=0$ for $i=1,\dots,d_{h}$ and $j=d_{h}+1,\dots,d_{h+1}$. Therefore, for $i=1,\dots,d_{h}$, we have
	\begin{align*}
		A_{\boldsymbol{\mu}(w_{i})}(v^{(h+1)})={}&A_{\boldsymbol{\mu}(w_{i})}(v^{(h)})-\sum_{j=d_{h}+1}^{d_{h+1}}A_{\boldsymbol{\mu}(w_{j})}(v^{(h)})\cdot{}A_{\boldsymbol{\mu}(w_{i})}(w_{j})\\
		={}&A_{\boldsymbol{\mu}(w_{i})}(v^{(h)})=0\rlap{.}
	\end{align*}
	For $i,j=d_{h}+1,\dots,d_{h+1}$, we have $A_{\boldsymbol{\mu}(w_{i})}(w_{j})=1$ if $i=j$, and $A_{\boldsymbol{\mu}(w_{i})}(w_{j})=0$ if $i\neq{}j$. Thus
	\begin{align*}
		A_{\boldsymbol{\mu}(w_{i})}(v^{(h+1)})={}&A_{\boldsymbol{\mu}(w_{i})}(v^{(h)})-\sum_{j=d_{h}+1}^{d_{h+1}}A_{\boldsymbol{\mu}(w_{j})}(v^{(h)})\cdot{}A_{\boldsymbol{\mu}(w_{i})}(w_{j})\\
		={}&A_{\boldsymbol{\mu}(w_{i})}(v^{(h)})-A_{\boldsymbol{\mu}(w_{i})}(v^{(h)})\cdot{}1=0\rlap{.}
	\end{align*}
	The argument above shows $A_{\boldsymbol{\mu}(w_{i})}(v^{(l-1)})=0$ for $i=1,\dots,d_{l-1}$. Since $\{w_{i}\}_{i=1}^{r-1}$ is a quasi-NP basis, we have $A_{\boldsymbol{\mu}(w_{i})}(w_{j})=0$ for $i=1,\dots,d_{l-1}$ and $j=d_{l-1}+1,\dots,r-1$. Therefore, for $i=1,\dots,d_{l-1}$, we have
	\begin{align*}
		A_{\boldsymbol{\mu}(w_{i})}(v^{(l)})={}&A_{\boldsymbol{\mu}(w_{i})}(v^{(l-1)})-\sum_{j=d_{l-1}+1}^{r-1}A_{\boldsymbol{\mu}(w_{j})}(v^{(l-1)})\cdot{}A_{\boldsymbol{\mu}(w_{i})}(w_{j})\\
		={}&A_{\boldsymbol{\mu}(w_{i})}(v^{(l-1)})=0\rlap{.}
	\end{align*}
	For $i,j=d_{l-1}+1,\dots,r-1$, we have $A_{\boldsymbol{\mu}(w_{i})}(w_{j})=1$ if $i=j$, and $A_{\boldsymbol{\mu}(w_{i})}(w_{j})=0$ if $i\neq{}j$. Thus
	\begin{align*}
		A_{\boldsymbol{\mu}(w_{i})}(v^{(l)})={}&A_{\boldsymbol{\mu}(w_{i})}(v^{(l-1)})-\sum_{j=d_{l-1}+1}^{r-1}A_{\boldsymbol{\mu}(w_{j})}(v^{(l-1)})\cdot{}A_{\boldsymbol{\mu}(w_{i})}(w_{j})\\
		={}&A_{\boldsymbol{\mu}(w_{i})}(v^{(l-1)})-A_{\boldsymbol{\mu}(w_{i})}(v^{(l-1)})\cdot{}1=0\rlap{.}
	\end{align*}
	Then, for $i=1,\dots,r-1$, we have $A_{\boldsymbol{\mu}(w_{i})}(v_{r})=A_{\boldsymbol{\mu}(v^{(l)})}(v^{(l)})^{-1}\cdot{}A_{\boldsymbol{\mu}(w_{i})}(v^{(l)})=0$. For $i=1,\dots,d_{l-1}$, set $v_{i}=w_{i}$. For $i=d_{l-1}+1,\dots,r-1$, let
	\begin{equation*}
		v_{i}=w_{i}-A_{\boldsymbol{\mu}(v_{r})}(w_{i})\cdot{}v_{r}\rlap{.}
	\end{equation*}
	Next, we show that $\{v_{i}\}_{i=1}^{r}$ is a quasi-NP-basis of $V$. For $i=d_{l-1}+1,\dots,r-1$, we have
	\begin{equation*}
		A_{\boldsymbol{\mu}(w_{i})}(v_{i})=A_{\boldsymbol{\mu}(w_{i})}(w_{i})-A_{\boldsymbol{\mu}(v_{r})}(w_{i})\cdot{}A_{\boldsymbol{\mu}(w_{i})}(v_{r})=A_{\boldsymbol{\mu}(w_{i})}(w_{i})=1\rlap{.}
	\end{equation*}
	For $i=d_{l-1}+1,\dots,r-1$ and $\boldsymbol{u}\in\mathrm{M}_{\mathrm{NP}}$, we show that $\operatorname{ord}A_{\boldsymbol{u}}(v_{i})>0$ if $\boldsymbol{\mu}(w_{i})\prec\boldsymbol{u}$. Supposing the contrary, then there exists $i'\in\{d_{l-1}+1,\dots,r-1\}$ and $\boldsymbol{u}'\in\mathrm{M}_{\mathrm{NP}}$, such that $\boldsymbol{\mu}(w_{i'})\prec\boldsymbol{u}'$ and $\operatorname{ord}A_{\boldsymbol{u}'}(v_{i'})=0$. The assumption $\boldsymbol{\mu}(w_{i'})\prec\boldsymbol{u}'$ implies $\operatorname{ord}A_{\boldsymbol{u}'}(w_{i'})>\operatorname{ord}A_{\boldsymbol{\mu}(w_{i'})}(w_{i'})=0$. Therefore, the assumption $\operatorname{ord}A_{\boldsymbol{u}'}(v_{i'})=0$ implies $\operatorname{ord}(A_{\boldsymbol{u'}}(w_{i'})-A_{\boldsymbol{u}'}(v_{i'}))=0$. Note that
	\begin{equation*}
		A_{\boldsymbol{u}'}(v_{i'})=A_{\boldsymbol{u}'}(w_{i'})-A_{\boldsymbol{\mu}(v_{r})}(w_{i'})\cdot{}A_{\boldsymbol{u}'}(v_{r})\rlap{.}
	\end{equation*}
	Thus $\operatorname{ord}(A_{\boldsymbol{\mu}(v_{r})}(w_{i'})\cdot{}A_{\boldsymbol{u}'}(v_{r}))=\operatorname{ord}(A_{\boldsymbol{u'}}(v_{i'})-A_{\boldsymbol{u}'}(w_{i'}))=0$, which implies $\operatorname{ord}A_{\boldsymbol{\mu}(v_{r})}(w_{i'})=\operatorname{ord}A_{\boldsymbol{u}'}(v_{r})=0$, since we have $\operatorname{ord}A_{\boldsymbol{\mu}(v_{r})}(w_{i'})\geq{}0$ and $\operatorname{ord}A_{\boldsymbol{u}'}(v_{r})\geq{}0$. Now $\operatorname{ord}A_{\boldsymbol{\mu}(v_{r})}(w_{i'})=0$ gives $\boldsymbol{\mu}(v_{r})\preceq\boldsymbol{\mu}(w_{i'})$, and $\operatorname{ord}A_{\boldsymbol{u}'}(v_{r})=0$ gives $\boldsymbol{u}'\preceq\boldsymbol{\mu}(v_{r})$. Thus, we obtain $\boldsymbol{u}'\preceq\boldsymbol{\mu}(w_{i'})$, which contradicts the premise $\boldsymbol{\mu}(w_{i'})\prec\boldsymbol{u}'$. Hence, for $i=d_{l-1}+1,\dots,r-1$ and $\boldsymbol{u}\in\mathrm{M}_{\mathrm{NP}}$, we have $\operatorname{ord}A_{\boldsymbol{u}}(v_{i})>0$ if $\boldsymbol{\mu}(w_{i})\prec\boldsymbol{u}$. This implies $\boldsymbol{\mu}(v_{i})\preceq\boldsymbol{\mu}(w_{i})$. Recall that we have already shown that $A_{\boldsymbol{\mu}(w_{i})}(v_{i})$, which implies $\boldsymbol{\mu}(w_{i})\preceq\boldsymbol{\mu}(v_{i})$. Thus $\boldsymbol{\mu}(w_{i})=\boldsymbol{\mu}(v_{i})$. Now, for $i=d_{l-1}+1,\dots,r-1$, we have $A_{\boldsymbol{\mu}(v_{i})}(v_{i})=A_{\boldsymbol{\mu}(w_{i})}(v_{i})=1$. For $i=d_{l-1}+1,\dots,r-1$ and $j=1,\dots,d_{l-1}$, we have
	\begin{equation*}
		A_{\boldsymbol{\mu}(v_{j})}(v_{i})=A_{\boldsymbol{\mu}(v_{j})}(w_{i})-A_{\boldsymbol{\mu}(v_{r})}(w_{i})\cdot{}A_{\boldsymbol{\mu}(v_{j})}(v_{r})=0\rlap{.}
	\end{equation*}
	For $i,j=d_{l-1}+1,\dots,r-1$ where $i\neq{}j$, we have
	\begin{align*}
		A_{\boldsymbol{\mu}(v_{j})}(v_{i})={}&A_{\boldsymbol{\mu}(v_{j})}(w_{i})-A_{\boldsymbol{\mu}(v_{r})}(w_{i})\cdot{}A_{\boldsymbol{\mu}(v_{j})}(v_{r})\\
		={}&A_{\boldsymbol{\mu}(w_{j})}(w_{i})=0\rlap{.}
	\end{align*}
	For $i=d_{l-1}+1,\dots,r-1$, we have $A_{\boldsymbol{\mu}(v_{i})}(v_{r})=A_{\boldsymbol{\mu}(w_{i})}(v_{r})=0$, and
	\begin{equation*}
		A_{\boldsymbol{\mu}(v_{r})}(v_{i})=A_{\boldsymbol{\mu}(v_{r})}(w_{i})-A_{\boldsymbol{\mu}(v_{r})}(w_{i})\cdot{}A_{\boldsymbol{\mu}(v_{r})}(v_{r})=0\rlap{.}
	\end{equation*}
	So far, we have verified all the requirements for $\{v_{i}\}_{i=1}^{r}$ to be a quasi-NP basis. In other words, we have shown that the assertion is true for $d=r$. By mathematical induction, we know that the assertion is true for $d\in\mathbb{Z}_{\geq{}1}$.
\end{proof}

\begin{dfn}\label{dfn:np-agreeable}
	\par Assume that $f$ is nondegenerate and $\dim\Delta(f)=n$. Let $\phi$ be a bijective $\sigma$-semilinear endomorphism on $V_{\mathrm{NP},m}$. We may write
	\begin{equation*}
		\phi(\pi^{pw(\boldsymbol{u}_{i})}\boldsymbol{x}^{\boldsymbol{u}_{i}})=\sum_{j=1}^{d_{\mathrm{NP}}}A_{\mathrm{NP},\phi}(i,j)\cdot\pi^{pw(\boldsymbol{u}_{j})}\boldsymbol{x}^{\boldsymbol{u}_{j}}
	\end{equation*}
	with unique $A_{\mathrm{NP},\phi}(i,j)\in{}K_{m}$. We say that $\phi$ is \emph{NP-agreeable} if for all $i,j=1,\dots,d_{\mathrm{NP}}$, we have
	\begin{equation*}
		\operatorname{ord}A_{\mathrm{NP},\phi}(i,j)\geq{}-w(\boldsymbol{u}_{j})\rlap{.}
	\end{equation*}
	Moreover, if $\phi$ is NP-agreeable, then we say that $(V_{\mathrm{NP},m},\phi,F^{\ast}_{\mathrm{NP}})\in\mathbf{M\overline{F}}{}^{\Phi}_{K_{m}/K_{m}}$ is \emph{NP-agreeable}.
\end{dfn}

\begin{lem}\label{lem:np-agreeable}
	\par Assume that $f$ is nondegenerate and $\dim\Delta(f)$. Let $\phi$ be a NP-agreeable bijective $\sigma$-semilinear endomorphism on $V_{\mathrm{NP},m}$. Let $V\subseteq{}V_{\mathrm{NP},m}$ be a subobject of $(V_{\mathrm{NP},m},\phi,F^{\ast}_{\mathrm{NP}})$ in $\mathbf{M\overline{F}}{}^{\Phi}_{K_{m}/K_{m}}$. Any quasi-NP basis of $V$ is an agreeable basis.
\end{lem}
\begin{proof}
	\par Let $d=\dim_{K_{m}}V$, and let $\{v_{i}\}_{i=1}^{d}$ be a quasi-NP basis of $V$. For $i=1,\dots,d$, we may write
	\begin{equation*}
		\phi(v_{i})=\sum_{j=1}^{d}A_{\phi}(v_{i},v_{j})\cdot{}v_{j}
	\end{equation*}
	with unique $A_{\phi}(v_{i},v_{j})\in{}K_{m}$. Then, for $\boldsymbol{u}\in\mathrm{M}_{\mathrm{NP}}$ and $i=1,\dots,d$, we have
	\begin{equation*}
		A_{\boldsymbol{u}}(\phi(v_{i}))=\sum_{j=1}^{d}A_{\phi}(v_{i},v_{j})\cdot{}A_{\boldsymbol{u}}(v_{j})\rlap{.}
	\end{equation*}
	At the same time, for $i=1,\dots,d$, since $v_{i}=\sum_{i'=1}^{d_{\mathrm{NP}}}A_{\boldsymbol{u}_{i'}}(v_{i})\cdot\pi^{pw(\boldsymbol{u}_{i'})}\boldsymbol{x}^{\boldsymbol{u}_{i'}}$, we have
	\begin{align*}
		\phi(v_{i})={}&\sum_{i'=1}^{d_{\mathrm{NP}}}\sigma(A_{\boldsymbol{u}_{i'}}(v_{i}))\cdot\phi(\pi^{pw(\boldsymbol{u}_{i'})})\boldsymbol{x}^{\boldsymbol{u}_{i'}}\\
		={}&\sum_{i'=1}^{d_{\mathrm{NP}}}\sum_{j'=1}^{d_{\mathrm{NP}}}\sigma(A_{\boldsymbol{u}_{i'}}(v_{i}))\cdot{}A_{\mathrm{NP},\phi}(i',j')\cdot\pi^{pw(\boldsymbol{u}_{j'})}\boldsymbol{x}^{\boldsymbol{u}_{j'}}\rlap{.}
	\end{align*}
	this gives $A_{\boldsymbol{u}_{j'}}(\phi(v_{i}))=\sum_{i'=1}^{d_{\mathrm{NP}}}\sigma(A_{\boldsymbol{u}_{i'}}(v_{i}))\cdot{}A_{\mathrm{NP},\phi}(i',j')$ for $j'=1,\dots,d_{\mathrm{NP}}$. Since $\phi$ is NP-agreeable, for $i',j'=1,\dots,d_{\mathrm{NP}}$, we have
	\begin{equation*}
		\operatorname{ord}A_{\mathrm{NP},\phi}(i',j')\geq{}-w(\boldsymbol{u}_{j'})\rlap{.}
	\end{equation*}
	Note that for $\boldsymbol{u}\in\mathrm{M}_{\mathrm{NP}}$ and $i=1,\dots,d$, we have $\operatorname{ord}A_{\boldsymbol{u}}(v_{i})\geq{}0$. Therefore, for $i=1,\dots,d$ and $j'=1,\dots,d_{\mathrm{NP}}$, we have
	\begin{align*}
		\operatorname{ord}A_{\boldsymbol{u}_{j'}}(\phi(v_{i}))\geq{}&\min\{\operatorname{ord}(\sigma(A_{\boldsymbol{u}_{i'}}(v_{i}))\cdot{}A_{\mathrm{NP},\phi}(i',j'))\mid{}i'=1,\dots,d_{\mathrm{NP}}\}\\
		\geq{}&\min\{\operatorname{ord}A_{\mathrm{NP},\phi}(i',j')\mid{}i'=1,\dots,d_{\mathrm{NP}}\}\geq{}-w(\boldsymbol{u}_{j'})\rlap{.}
	\end{align*}
	For $i=1,\dots,d$, let $\mathfrak{Bad}_{\phi}(i)=\{j\in\{1,\dots,d\}\mid\operatorname{ord}A_{\phi}(v_{i},v_{j})<-w(v_{j})\}$. To show that $\{v_{i}\}_{i=1}^{d}$ is agreeable, it suffices to show that $\mathfrak{Bad}_{\phi}(i)=\emptyset$ for all $i=1,\dots,d$. Supposing the contrary, then there exists $i_{0}\in\{1,\dots,d\}$, such that $\mathfrak{Bad}_{\phi}(i_{0})\neq\emptyset$. Let $j_{0}=\min\mathfrak{Bad}_{\phi}(i_{0})$. Since $\{v_{i}\}_{i=1}^{d}$ is a quasi-NP basis, for $j=1,\dots,d$, we have $A_{\boldsymbol{\mu}(v_{j_{0}})}(v_{j})=0$ if $w(v_{j})\geq{}w(v_{j_{0}})$. Thus
	\begin{align*}
		A_{\boldsymbol{\mu}(v_{j_{0}})}(\phi(v_{i_{0}}))={}&\sum_{j=1}^{d}A_{\phi}(v_{i_{0}},v_{j})\cdot{}A_{\boldsymbol{\mu}(v_{j_{0}})}(v_{j})\\
		={}&A_{\phi}(v_{i_{0}},v_{j_{0}})+\sum_{w(v_{j})<w(v_{j_{0}})}A_{\phi}(v_{i_{0}},v_{j})\cdot{}A_{\boldsymbol{\mu}(v_{j_{0}})}(v_{j})\rlap{.}
	\end{align*}
	We show that $\operatorname{ord}A_{\boldsymbol{\mu}(v_{j_{0}})}(\phi(v_{i_{0}}))=\operatorname{ord}A_{\phi}(v_{i_{0}},v_{j_{0}})$. If $w(v_{j_{0}})=\min\{w(v_{j})\mid{}j=1,\dots,d\}$, then $A_{\boldsymbol{\mu}(v_{j_{0}})}(\phi(v_{i_{0}}))=A_{\phi}(v_{i_{0}},v_{j_{0}})$. If $w(v_{j_{0}})>\min\{w(v_{j})\mid{}j=1,\dots,d\}$, then there exists $j_{1}\in\{1,\dots,d\}$, such that $w(j_{1})<w(j_{0})$. Since $\{v_{i}\}_{i=1}^{d}$ is filtration-generating, we have $j_{1}<j_{0}$, which implies $j_{1}\notin\mathfrak{Bad}_{\phi}(i_{0})$. Therefore, for all $j\in\{1,\dots,d\}$ such that $w(v_{j})<w(v_{j_{0}})$, we have
	\begin{equation*}
		\operatorname{ord}A_{\phi}(v_{i_{0}},v_{j})\geq{}-w(v_{j})>-w(v_{j_{0}})>\operatorname{ord}A_{\phi}(v_{i_{0}},v_{j_{0}})\rlap{,}
	\end{equation*}
	which implies $\operatorname{ord}A_{\boldsymbol{\mu}(v_{j_{0}})}(\phi(v_{i_{0}}))=\operatorname{ord}A_{\phi}(v_{i_{0}},v_{j_{0}})$. On the other hand, we have already shown
	\begin{equation*}
		\operatorname{ord}A_{\boldsymbol{\mu}(v_{j_{0}})}(\phi(v_{i_{0}}))\geq{}-w(\boldsymbol{\mu}(v_{j_{0}}))\rlap{,}
	\end{equation*}
	which implies $\operatorname{ord}A_{\boldsymbol{\mu}(v_{j_{0}})}(\phi(v_{i_{0}}))\geq{}w(v_{j_{0}})$, since $w(\boldsymbol{\mu}(v_{j_{0}}))\leq{}w(v_{j_{0}})$. However, this contradicts the premise that $\operatorname{ord}A_{\phi}(v_{i_{0}},v_{j_{0}})<-w(v_{j_{0}})$. Therefore, we have $\mathfrak{Bad}_{\phi}(i)=\emptyset$ for all $i=1,\dots,d$, and hence $\{v_{i}\}_{i=1}^{d}$ is agreeable.
\end{proof}

\begin{thm}\label{thm:weak-admissibility}
	\par If $f$ is nondegenerate and $\dim\Delta(f)=n$, then $(V_{\mathrm{NP}},\tilde{\phi}_{\mathrm{NP}},F^{\ast}_{\mathrm{NP}})\in\mathbf{M\overline{F}}{}^{\Phi}_{K_{1}/K_{1}}$ is weakly admissible.
\end{thm}
\begin{proof}
	\par Let $\tilde{\phi}_{\mathrm{NP},m}=\sigma\otimes\tilde{\phi}_{\mathrm{NP}}$, which is a bijective $\sigma$-semilinear endomorphism on $V_{\mathrm{NP},m}$. We first show that $\tilde{\phi}_{\mathrm{NP},m}$ is NP-agreeable. For $\boldsymbol{u}\in\mathrm{M}_{\mathrm{NP}}$, we have $\pi^{pw(\boldsymbol{u})}\boldsymbol{x}^{\boldsymbol{u}}\in{}L(\frac{p}{p-1},0)$. By the proof of \cite{dwork1962zeta}*{Lemma~4.1}, we have $\exp(\varphi(\widetilde{F}(\boldsymbol{x}))-\widetilde{F}(\boldsymbol{x}))\in{}L(\frac{1}{p-1},0)$, which implies that
	\begin{equation*}
		\exp(\varphi(\widetilde{F}(\boldsymbol{x}))-\widetilde{F}(\boldsymbol{x}))\cdot\pi^{pw(\boldsymbol{u})}\boldsymbol{x}^{\boldsymbol{u}}\in{}L\left(\frac{1}{p-1},0\right)
	\end{equation*}
	For $i=1,\dots,d_{\mathrm{NP}}$, there exist unique $\widetilde{A}_{\mathrm{NP}}(i,1),\dots,\widetilde{A}_{\mathrm{NP}}(i,d_{\mathrm{NP}})\in{}K_{m}$, such that
	\begin{equation*}
		\tilde{\phi}_{\mathrm{NP},m}(\pi^{pw(\boldsymbol{u}_{i})}\boldsymbol{x}^{\boldsymbol{u}_{i}})=\sum_{j=1}^{d_{\mathrm{NP}}}\widetilde{A}_{\mathrm{NP}}(i,j)\cdot\pi^{pw(\boldsymbol{u}_{j})}\boldsymbol{x}^{\boldsymbol{u}_{j}}\rlap{,}
	\end{equation*}
	where, by \cref{cor:decomposition}, we know that $\widetilde{A}_{\mathrm{NP}}(i,j),\dots,\widetilde{A}_{\mathrm{NP}}(i,d_{\mathrm{NP}})$ are characterized by
	\begin{equation*}
		\exp(\varphi(\widetilde{F}(\boldsymbol{x}))-\widetilde{F}(\boldsymbol{x}))\cdot\pi^{pw(\boldsymbol{u}_{i})}\boldsymbol{x}^{\boldsymbol{u}_{i}}\equiv\sum_{j=1}^{d_{\mathrm{NP}}}\widetilde{A}_{\mathrm{NP}}(i,j)\cdot\pi^{pw(\boldsymbol{u}_{j})}\boldsymbol{x}^{\boldsymbol{u}_{j}}\bmod\sum_{l=1}^{n}D_{\widetilde{F},i}L\left(\frac{1}{p-1}\right)\rlap{.}
	\end{equation*}
	Therefore, by \cref{prop:decomposition.growth}, for $i=1,\dots,d_{\mathrm{NP}}$, we have
	\begin{equation*}
		\sum_{j=1}^{d_{\mathrm{NP}}}\widetilde{A}_{\mathrm{NP}}(i,j)\cdot\pi^{pw(\boldsymbol{u}_{j})}\boldsymbol{x}^{\boldsymbol{u}_{j}}\in{}V\left(\frac{1}{p-1},0\right)\rlap{.}
	\end{equation*}
	This means $\operatorname{ord}(\widetilde{A}_{\mathrm{NP}}(i,j)\cdot\pi^{pw(\boldsymbol{u}_{j})})\geq{}\frac{1}{p-1}\cdot{}w(\boldsymbol{u}_{j})$ for all $i,j=1,\dots,d_{\mathrm{NP}}$, namely
	\begin{equation*}
		\operatorname{ord}\widetilde{A}_{\mathrm{NP}}(i,j)\geq\frac{1}{p-1}\cdot{}w(\boldsymbol{u}_{j})-\operatorname{ord}\pi^{pw(\boldsymbol{u}_{j})}=-w(\boldsymbol{u}_{j})\rlap{.}
	\end{equation*}
	This shows that $\tilde{\phi}_{\mathrm{NP},m}$ is NP-agreeable. Next, we show that $(V_{\mathrm{NP},m},\tilde{\phi}_{\mathrm{NP},m},F^{\ast}_{\mathrm{NP}})$ is weakly admissible. Let $V\subseteq{}V_{\mathrm{NP},m}$ be a subobject of $(V_{\mathrm{NP},m},\tilde{\phi}_{\mathrm{NP},m},F^{\ast}_{\mathrm{NP}})$ in $\mathbf{M\overline{F}}{}^{\Phi}_{K_{m}/K_{m}}$. By \cref{lem:np-basis}, we know that $V$ admits a quasi-NP basis with respect to $F^{\ast}_{\mathrm{NP}}$. By \cref{lem:np-agreeable}, we know that a quasi-NP basis of $V$ is agreeable with respect to $\tilde{\phi}_{\mathrm{NP},m}$ and $F^{\ast}_{\mathrm{NP}}$. Thus $V$ admits an agreeable basis, so by \cref{prop:agreeable}, we get $t_{\mathrm{N}}(V,\tilde{\phi}_{\mathrm{NP},m})\geq{}t_{\mathrm{H}}(V,F^{\ast}_{\mathrm{NP},m})$. By the proof of \cite{adolphson1993twisted}*{Theorem~3.17}, we get $t_{\mathrm{N}}(V_{\mathrm{NP},m},\tilde{\phi}_{\mathrm{NP},m})=t_{\mathrm{H}}(V_{\mathrm{NP},m},F^{\ast}_{\mathrm{NP}})$. Thus $(V_{\mathrm{NP},m},\tilde{\phi}_{\mathrm{NP},m},F^{\ast}_{\mathrm{NP}})$ is weakly admissible as an object in $\mathbf{M\overline{F}}{}^{\Phi}_{K_{m}/K_{m}}$. Finally, we show that $(V_{\mathrm{NP}},\tilde{\phi}_{\mathrm{NP}},F^{\ast}_{\mathrm{NP}})$ is weakly admissible. If $(W,\tilde{\phi}_{\mathrm{NP}},F^{\ast})$ is a subobject of $(V_{\mathrm{NP}},\tilde{\phi}_{\mathrm{NP}},F^{\ast}_{\mathrm{NP}})$ in $\mathbf{M\overline{F}}{}^{\Phi}_{K_{1}/K_{1}}$, then $(W_{m},\tilde{\phi}_{\mathrm{NP},m},F^{\ast}_{\mathrm{NP}})$ is a subobject of $(V_{\mathrm{NP},m},\tilde{\phi}_{\mathrm{NP},m},F^{\ast}_{\mathrm{NP}})$ in $\mathbf{M\overline{F}}{}^{\Phi}_{K_{m}/K_{m}}$, where $W_{m}=K_{m}\otimes_{K_{1}}W$. This implies that $t_{\mathrm{N}}(W,\tilde{\phi}_{\mathrm{NP}})=t_{\mathrm{N}}(W_{m},\tilde{\phi}_{\mathrm{NP},m})\geq{}t_{\mathrm{H}}(W_{m},F^{\ast}_{\mathrm{NP}})=t_{\mathrm{H}}(W,F^{\ast}_{\mathrm{NP}})$. On the other hand, we note that $t_{\mathrm{N}}(V_{\mathrm{NP}},\tilde{\phi}_{\mathrm{NP}})=t_{\mathrm{N}}(V_{\mathrm{NP},m},\tilde{\phi}_{\mathrm{NP},m})=t_{\mathrm{H}}(V_{\mathrm{NP},m},F^{\ast}_{\mathrm{NP}})=t_{\mathrm{H}}(V_{\mathrm{NP}},F^{\ast}_{\mathrm{NP}})$. Therefore, we obtain that $(V_{\mathrm{NP}},\tilde{\phi}_{\mathrm{NP}},F^{\ast}_{\mathrm{NP}})$ is weakly admissible.
\end{proof}

\begin{dfn}\label{dfn:np-dominating}
	\par Assume that $f$ is nondegenerate and $\dim\Delta(f)=n$. Let $T:V_{\mathrm{NP},m}\rightarrow{}V_{\mathrm{NP},m}$ be an automorphism in $\mathbf{Mod}_{K_{m}}$. For $i,i'=1,\dots,d_{\mathrm{NP}}$, we may write
	\begin{equation*}
		T(\pi^{w(\boldsymbol{u}_{i})}\boldsymbol{x}^{\boldsymbol{u}_{i}})=\sum_{i=1}^{d_{\mathrm{NP}}}T(i,j)\cdot\pi^{w(\boldsymbol{u}_{j})}\boldsymbol{x}^{\boldsymbol{u}_{j}}
	\end{equation*}
	with unique $T(i,j)\in{}K_{m}$. We say that $T$ is \emph{NP-dominating} if for $i,j,i',j'=1,\dots,d_{\mathrm{NP}}$, we have
	\begin{equation*}
		\operatorname{ord}T(i,j)+\operatorname{ord}T^{-1}(i',j')\geq{}w(\boldsymbol{u}_{j})-w(\boldsymbol{u}_{i})\rlap{.}
	\end{equation*}
\end{dfn}

\begin{lem}\label{lem:np-dominating}
	\par Assume that $f$ is nondegenerate and $\dim\Delta(f)=n$. Let $\phi$ and $\phi'$ be bijective $\sigma$-semilinear endomorphisms on $V_{\mathrm{NP},m}$. Let $T:(V_{\mathrm{NP},m},\phi')\rightarrow(V_{\mathrm{NP},m},\phi)$ be an isomorphism in $\mathbf{Mod}{}^{\Phi}_{K_{m}}$. If $\phi'$ is NP-agreeable and $T$ is NP-dominating, then $\phi$ is NP-agreeable.
\end{lem}
\begin{proof}
	\par We first fix some notations. For $i,i'=1,\dots,d_{\mathrm{NP}}$, we may write
	\begin{align*}
		\phi(\pi^{pw(\boldsymbol{u}_{i})}\boldsymbol{x}^{\boldsymbol{u}_{i}})={}&\sum_{j=1}^{d_{\mathrm{NP}}}A(i,j)\cdot\pi^{pw(\boldsymbol{u}_{j})}\boldsymbol{x}^{\boldsymbol{u}_{j}}\rlap{,}\\
		\phi'(\pi^{pw(\boldsymbol{u}_{i'})}\boldsymbol{x}^{\boldsymbol{u}_{i'}})={}&\sum_{j'=1}^{d_{\mathrm{NP}}}A'(i',j')\cdot\pi^{pw(\boldsymbol{u}_{j'})}\boldsymbol{x}^{\boldsymbol{u}_{j'}}\rlap{,}
	\end{align*}
	with unique $A(i,j)$ and $A(i',j')$ in $K_{m}$. Since $T$ is an isomorphism in $\mathbf{Mod}{}^{\Phi}_{K_{m}}$, it is Frobenius-compatible, which means $\phi=T^{-1}\comp\phi'\comp{}T$. Therefore, for $i=1,\dots,d_{\mathrm{NP}}$, we have
	\begin{align*}
		\phi(\pi^{pw(\boldsymbol{u}_{i})}\boldsymbol{x}^{\boldsymbol{u}_{i}})={}&T^{-1}(\phi'(T(\pi^{pw(\boldsymbol{u}_{i})})))\\
		={}&\sum_{i'=1}^{d_{\mathrm{NP}}}T^{-1}(\phi'((-p)^{w(\boldsymbol{u}_{i})}T(i,i')\cdot{}\pi^{w(\boldsymbol{u}_{i'})}\boldsymbol{x}^{\boldsymbol{u}_{i'}}))\\
		={}&\sum_{i',j'=1}^{d_{\mathrm{NP}}}(-p)^{w(\boldsymbol{u}_{i})}\sigma(T(i,i'))T^{-1}((-p)^{-w(\boldsymbol{u}_{i'})}A'(i',j')\cdot\pi^{pw(\boldsymbol{u}_{j'})}\boldsymbol{x}^{\boldsymbol{u}_{j'}})\\
		={}&\sum_{i',j',j=1}^{d_{\mathrm{NP}}}(-p)^{w(\boldsymbol{u}_{i})-w(\boldsymbol{u}_{i'})+w(\boldsymbol{u}_{j'})-w(\boldsymbol{u}_{j})}\sigma(T(i,i'))A'(i',j')T^{-1}(j',j)\cdot\pi^{pw(\boldsymbol{u}_{j})}\boldsymbol{x}^{\boldsymbol{u}_{j}}\rlap{.}
	\end{align*}
	Hence, for $i,j=1,\dots,d_{\mathrm{NP}}$, we have
	\begin{equation*}
		A(i,j)=\sum_{i',j'=1}^{d_{\mathrm{NP}}}(-p)^{w(\boldsymbol{u}_{i})-w(\boldsymbol{u}_{i'})+w(\boldsymbol{u}_{j'})-w(\boldsymbol{u}_{j})}\sigma(T(i,i'))A'(i',j')T^{-1}(j',j)
	\end{equation*}
	Since $\phi'$ is NP-agreeable, we have $\operatorname{ord}A'(i',j')\geq{}-w(\boldsymbol{u}_{j'})$ for all $i',j'=1,\dots,d_{\mathrm{NP}}$. Since $T$ is dominating, we have $\operatorname{ord}T(i,i')+\operatorname{ord}T^{-1}(j',j)\geq{}w(\boldsymbol{u}_{i'})-w(\boldsymbol{u}_{i})$ for all $i,j,i',j'=1,\dots,d_{\mathrm{NP}}$. Thus
	\begin{align*}
		&\operatorname{ord}\left\{(-p)^{w(\boldsymbol{u}_{i})-w(\boldsymbol{u}_{i'})+w(\boldsymbol{u}_{j'})-w(\boldsymbol{u}_{j})}\sigma(T(i,i'))A'(i',j')T^{-1}(j',j)\right\}\\
		={}&w(\boldsymbol{u}_{i})-w(\boldsymbol{u}_{i'})+w(\boldsymbol{u}_{j'})-w(\boldsymbol{u}_{j})+\operatorname{ord}A'(i',j')+\operatorname{ord}T(i,i')+\operatorname{ord}T^{-1}(j',j)\\
		\geq{}&w(\boldsymbol{u}_{i})-w(\boldsymbol{u}_{i'})+w(\boldsymbol{u}_{j'})-w(\boldsymbol{u}_{j})-w(\boldsymbol{u}_{j'})+w(\boldsymbol{u}_{i'})-w(\boldsymbol{u}_{i})=-w(\boldsymbol{u}_{j})\rlap{.}
	\end{align*}
	Therefore, for $i,j=1,\dots,d_{\mathrm{NP}}$, we have $\operatorname{ord}A_{\mathrm{NP}}(i,j)\geq{}-w(\boldsymbol{u}_{j})$, so that $\phi$ is NP-agreeable.
\end{proof}

\begin{lem}\label{lem:valuation}
	\par Assume that $p\neq{}2$. We have $\exp(\widetilde{F}(\boldsymbol{x})-\widehat{F}(\boldsymbol{x}))-1\in{}L(\frac{1}{p-1},p-2)$
\end{lem}
\begin{proof}
	\par Note that $\exp(\widetilde{F}(\boldsymbol{x})-\widehat{F}(\boldsymbol{x}))=\prod_{\boldsymbol{u}\in\operatorname{supp}f}\exp((\gamma-\pi)\hat{\alpha}_{\boldsymbol{u}}\boldsymbol{x}^{\boldsymbol{u}})\cdot\exp(\sum_{l=1}^{\infty}\gamma_{l}\hat{\alpha}_{\boldsymbol{u}}^{p^{l}}\boldsymbol{x}^{p^{l}\boldsymbol{u}})$. By the proof of \cref{lem:valuation.transformation}, for $\boldsymbol{u}\in\operatorname{supp}f$, we have $\exp((\gamma-\pi)\hat{\alpha}_{\boldsymbol{u}}\boldsymbol{x}^{\boldsymbol{u}})\in{}L(p-1,0)$, which implies
	\begin{equation*}
		\exp((\gamma-\pi)\hat{\alpha}_{\boldsymbol{u}}\boldsymbol{x}^{\boldsymbol{u}})-1\in{}L\left(\frac{1}{p-1},p-2\right)\rlap{,}
	\end{equation*}
	since $p-1\geq{}\frac{1}{p-1}+p-2\geq{}\frac{1}{p-1}\cdot{}w(\boldsymbol{u})+p-2$. For $l\geq{}1$ and $\boldsymbol{u}\in\operatorname{supp}f$, we have
	\begin{equation*}
		\exp(\gamma_{l}(\hat{\alpha}_{\boldsymbol{u}}\boldsymbol{x}^{\boldsymbol{u}})^{p^{l}})=\sum_{i=0}^{\infty}\frac{\gamma_{l}^{i}\hat{\alpha}_{\boldsymbol{u}}^{p^{l}i}}{i!}\cdot\boldsymbol{x}^{p^{l}i\cdot\boldsymbol{u}}\rlap{.}
	\end{equation*}
	Since $w(\boldsymbol{u})\leq{}1$ for $\boldsymbol{u}\in\operatorname{supp}f$, to show that $exp(\gamma_{l}(\hat{\alpha}_{\boldsymbol{u}}\boldsymbol{x}^{\boldsymbol{u}})^{p^{l}})-1\in{}L(\frac{1}{p-1},p-2)$, it suffices to show $\operatorname{ord}\frac{\gamma_{l}^{i}}{i!}\geq{}\frac{p^{l}i}{p-1}+p-2$ for $i\geq{}i$ and $l\geq{}1$. If $i=1$, this assertion follows from the observation
	\begin{equation*}
		\operatorname{ord}\gamma_{l}=\frac{p^{l+1}}{p-1}-l-1\geq{}\frac{p^{l}}{p-1}+p-2
	\end{equation*}
	for $l\geq{}1$. For $i\geq{}2$ and $l\geq{}1$, the assertion follows from the observation
	\begin{equation*}
		\operatorname{ord}\frac{\gamma_{l}^{i}}{i!}\geq\left(\frac{p^{l+1}}{p-1}-l-1-\frac{1}{p-1}\right)\cdot{}i\geq{}\frac{p^{l}i}{p-1}+p-2\rlap{.}\qedhere
	\end{equation*}
\end{proof}

\begin{proof}[Proof of \cref{thm:main}]
	\par By \cref{thm:object.cohomology}, we know that $\iota_{\widehat{F}}$ is an isomorphism, so that we obtain an object
	\begin{equation*}
		((H^{n}_{\mathrm{rig}}(\mathbb{T}^{n}_{k}/K_{1},f^{\ast}\mathcal{L}_{\pi}),\phi_{f}),(H^{n}_{\mathrm{dR}}(\mathbb{T}^{n}_{K_{1}},\nabla_{\widehat{F}}),F^{\ast}_{\mathrm{irr}}),\iota_{\widehat{F}})\in\mathbf{M\overline{F}}{}^{\Phi}_{K_{1}/K_{1}}\rlap{.}
	\end{equation*}
	By \cref{prop:np-module}, we get an isomorphism
	\begin{equation*}
		((H^{n}_{\mathrm{rig}}(\mathbb{T}^{n}_{k}/K_{1},f^{\ast}\mathcal{L}_{\pi}),\phi_{f}),(H^{n}_{\mathrm{dR}}(\mathbb{T}^{n}_{K_{1}},\nabla_{\widehat{F}}),F^{\ast}_{\mathrm{irr}}),\iota_{\widehat{F}})\cong(V_{\mathrm{NP}},\hat{\phi}_{\mathrm{NP}},F^{\ast}_{\mathrm{NP}})
	\end{equation*}
	in $\mathbf{M\overline{F}}{}^{\Phi}_{K_{1}/K_{1}}$. Therefore, in order to prove the theorem, it suffices to show that $(V_{\mathrm{NP}},\hat{\phi}_{\mathrm{NP}},F^{\ast}_{\mathrm{NP}})$ is weakly admissible. Let $\hat{\phi}_{\mathrm{NP},m}=\sigma\otimes\hat{\phi}_{\mathrm{NP}}$, which is a bijective $\sigma$-semilinear endomorphism on $V_{\mathrm{NP},m}$. We first show that $\hat{\phi}_{\mathrm{NP},m}$ is NP-agreeable. By the proof of \cref{thm:weak-admissibility}, we know that $\tilde{\phi}_{\mathrm{NP},m}$ is NP-agreeable. By \cref{prop:transformation}, the automorphism $T_{\mathrm{NP},m}=1\otimes{}T_{\mathrm{NP}}$ on $V_{\mathrm{NP},m}$ gives an isomorphism $T_{\mathrm{NP},m}:(V_{\mathrm{NP},m},\tilde{\phi}_{\mathrm{NP},m})\rightarrow(V_{\mathrm{NP},m},\hat{\phi}_{\mathrm{NP},m})$ in $\mathbf{Mod}{}^{\Phi}_{K_{1}}$. By \cref{lem:np-dominating}, we only need to show that $T_{\mathrm{NP},m}$ is NP-dominating. For $i,i'=1,\dots,d_{\mathrm{NP}}$, we may write
	\begin{align*}
		T_{\mathrm{NP},m}(\pi^{w(\boldsymbol{u}_{i})}\boldsymbol{x}^{\boldsymbol{u}_{i}})={}&\sum_{j=1}^{d_{\mathrm{NP}}}T_{\mathrm{NP},m}(i,j)\cdot\pi^{w(\boldsymbol{u}_{j})}\boldsymbol{x}^{\boldsymbol{u}_{j}}\rlap{,}\\
		T_{\mathrm{NP},m}^{-1}(\pi^{w(\boldsymbol{u}_{i'})}\boldsymbol{x}^{\boldsymbol{u}_{i'}})={}&\sum_{j'=1}^{d_{\mathrm{NP}}}T^{-1}_{\mathrm{NP},m}(i,j)\cdot\pi^{w(\boldsymbol{u}_{j'})}\boldsymbol{x}^{\boldsymbol{u}_{j'}}\rlap{,}
	\end{align*}
	with unique $T(i,j)$ and $T^{-1}(i,j)$ in $K_{m}$, which, by \cref{cor:decomposition}, are characterized by
	\begin{align*}
		\exp(\widetilde{F}(\boldsymbol{x})-\widehat{F}(\boldsymbol{x}))\cdot\pi^{w(\boldsymbol{u}_{i})}\boldsymbol{x}^{\boldsymbol{u}_{i}}\equiv{}&\sum_{j=1}^{d_{\mathrm{NP}}}T_{\mathrm{NP},m}(i,j)\cdot\pi^{w(\boldsymbol{u}_{j})}\boldsymbol{x}^{\boldsymbol{u}_{j}}\bmod\sum_{i=1}^{n}D_{\widehat{F},i}L\left(\frac{1}{p-1}\right)\rlap{,}\\
		\exp(\widehat{F}(\boldsymbol{x})-\widetilde{F}(\boldsymbol{x}))\cdot\pi^{w(\boldsymbol{u}_{i'})}\boldsymbol{x}^{\boldsymbol{u}_{i'}}\equiv{}&\sum_{j'=1}^{d_{\mathrm{NP}}}T^{-1}_{\mathrm{NP},m}(i,j)\cdot\pi^{w(\boldsymbol{u}_{j'})}\boldsymbol{x}^{\boldsymbol{u}_{j'}}\bmod\sum_{i=1}^{n}D_{\widehat{F},i}L\left(\frac{1}{p-1}\right)\rlap{.}
	\end{align*}
	Note that we have $T_{\mathrm{NP},m}(i,i)=1$ and $T_{\mathrm{NP},m}(j,j)=1$ for $i,j=1,\dots,d_{\mathrm{NP}}$. By \cref{lem:valuation}, we have $\exp(\widetilde{F}(\boldsymbol{x})-\widehat{F}(\boldsymbol{x}))-1\in{}L(\frac{1}{p-1},p-2)$, which implies $\exp(\widehat{F}(\boldsymbol{x})-\widetilde{F}(\boldsymbol{x}))\in{}L(\frac{1}{p-1},p-2)$. By \cref{prop:decomposition.growth}, for $i,j,i',j'=1,\dots,d_{\mathrm{NP}}$ where $i\neq{}i'$ and $j\neq{}j'$, we have
	\begin{align*}
		T_{\mathrm{NP},m}(i,j)\cdot\pi^{w(\boldsymbol{u}_{j})}\boldsymbol{x}^{\boldsymbol{u}_{j}}\in{}&L\left(\frac{1}{p-1},p-2\right)\rlap{,}\\
		T^{-1}_{\mathrm{NP},m}(i',j')\cdot\pi^{w(\boldsymbol{u}_{j'})}\boldsymbol{x}^{\boldsymbol{u}_{j'}}\in{}&L\left(\frac{1}{p-1},p-2\right)\rlap{,}
	\end{align*}
	which implies that $\operatorname{ord}T_{\mathrm{NP},m}(i,j)\geq{}p-2$ and $\operatorname{ord}T^{-1}_{\mathrm{NP},m}(i',j')\geq{}p-2$. As a summary, for $i,j,i',j'=1,\dots,d_{\mathrm{NP}}$, we have
	\begin{align*}
		\operatorname{ord}T_{\mathrm{NP},m}(i,j)={}&\begin{cases}
			=0
				&i=j\rlap{,}\\
			\geq{}p-2
				&i\neq{}j\rlap{,}
		\end{cases}\\
		\operatorname{ord}T^{-1}_{\mathrm{NP},m}(i',j')={}&\begin{cases}
			=0
				&i'=j'\rlap{,}\\
			\geq{}p-2
				&i'\neq{}j'\rlap{,}
		\end{cases}
	\end{align*}
	which gives $\operatorname{ord}T_{\mathrm{NP},m}(i,j)+\operatorname{ord}T^{-1}_{\mathrm{NP},m}(i',j')\geq{}p-2\geq\ell_{\mathrm{HT}}(V_{\mathrm{NP},m},F^{\ast}_{\mathrm{NP}})\geq{}w(\boldsymbol{u}_{j})-w(\boldsymbol{u}_{i})$. Here, we note that $\ell_{\mathrm{HT}}(V_{\mathrm{NP},m},F^{\ast}_{\mathrm{NP}})=\ell_{\mathrm{HT}}(V_{\mathrm{NP}},F^{\ast}_{\mathrm{NP}})$. This shows that $T_{\mathrm{NP},m}$ is NP-dominating, and hence $\hat{\phi}_{\mathrm{NP},m}$ is NP-agreeable. Next, we show that $(V_{\mathrm{NP},m},\hat{\phi}_{\mathrm{NP},m},F^{\ast}_{\mathrm{NP}})$ is weakly admissible as an object in $\mathbf{M\overline{F}}{}^{\Phi}_{K_{m}/K_{m}}$. Let $V\subseteq{}V_{\mathrm{NP},m}$ be a subobject of $(V_{\mathrm{NP},m},\hat{\phi}_{\mathrm{NP},m},F^{\ast}_{\mathrm{NP}})$ in $\mathbf{M\overline{F}}{}^{\Phi}_{K_{m}/K_{m}}$. By \cref{lem:np-basis}, we know that $V$ admits a quasi-NP basis with respect to $F^{\ast}_{\mathrm{NP}}$. By \cref{lem:np-agreeable}, a quasi-NP basis of $V$ is agreeable with respect to $\hat{\phi}_{\mathrm{NP},m}$ and $F^{\ast}_{\mathrm{NP}}$. Thus $V$ admits an agreeable basis, so by \cref{prop:agreeable}, we get $t_{\mathrm{N}}(V,\hat{\phi}_{\mathrm{NP},m})\geq{}t_{\mathrm{H}}(V,F^{\ast}_{\mathrm{NP},m})$. Since $(V_{\mathrm{NP},m},\tilde{\phi}_{\mathrm{NP},m})$ and $(V_{\mathrm{NP},m},\hat{\phi}_{\mathrm{NP},m})$ are isomorphic in $\mathbf{Mod}{}^{\Phi}_{K_{1}}$, their Newton numbers are equal. By \cref{thm:weak-admissibility}, we have $t_{\mathrm{N}}(V_{\mathrm{NP},m},\tilde{\phi}_{\mathrm{NP},m})=t_{\mathrm{H}}(V_{\mathrm{NP},m},F^{\ast}_{\mathrm{NP}})$. This implies $t_{\mathrm{N}}(V_{\mathrm{NP},m},\hat{\phi}_{\mathrm{NP},m})=t_{\mathrm{H}}(V_{\mathrm{NP},m},F^{\ast}_{\mathrm{NP}})$, and hence $(V_{\mathrm{NP},m},\hat{\phi}_{\mathrm{NP},m},F^{\ast}_{\mathrm{NP}})$ is weakly admissible. Finally, we show that $(V_{\mathrm{NP}},\hat{\phi}_{\mathrm{NP}},F^{\ast}_{\mathrm{NP}})$ is weakly admissible. Let $(W,\hat{\phi}_{\mathrm{NP}},F^{\ast})$ be a subobject of $(V_{\mathrm{NP}},\hat{\phi}_{\mathrm{NP}},F^{\ast}_{\mathrm{NP}})$ in $\mathbf{M\overline{F}}{}^{\Phi}_{K_{1}/K_{1}}$, and let $W_{m}=K_{m}\otimes_{K_{1}}W$. Note that $(W_{m},\hat{\phi}_{\mathrm{NP},m},F^{\ast}_{\mathrm{NP}})$ is a subobject of $(V_{\mathrm{NP},m},\hat{\phi}_{\mathrm{NP},m},F^{\ast}_{\mathrm{NP}})$ in $\mathbf{M\overline{F}}{}^{\Phi}_{K_{m}/K_{m}}$. This implies that $t_{\mathrm{N}}(W,\hat{\phi}_{\mathrm{NP}})=t_{\mathrm{N}}(W_{m},\hat{\phi}_{\mathrm{NP},m})\geq{}t_{\mathrm{H}}(W_{m},F^{\ast}_{\mathrm{NP}})=t_{\mathrm{H}}(W,F^{\ast}_{\mathrm{NP}})$. On the other hand, note that $t_{\mathrm{N}}(V_{\mathrm{NP}},\hat{\phi}_{\mathrm{NP}})=t_{\mathrm{N}}(V_{\mathrm{NP},m},\hat{\phi}_{\mathrm{NP},m})=t_{\mathrm{H}}(V_{\mathrm{NP},m},F^{\ast}_{\mathrm{NP}})=t_{\mathrm{H}}(V_{\mathrm{NP}},F^{\ast}_{\mathrm{NP}})$. Thus $(V_{\mathrm{NP}},\hat{\phi}_{\mathrm{NP}},F^{\ast}_{\mathrm{NP}})$ is weakly admissible.
\end{proof}

\begin{rmk}\label{rmk:weak-admissibility}
	\par Assume that $p\neq{}2$. Let $F\in{}A_{0}=O_{1}[x_{1},\dots,x_{n},(x_{1}\dots{}x_{n})^{-1}]$ be a nondegenerate Laurent polynomial such that $\pi^{-1}F\in{}A_{0}$ and the reduction of $\pi^{-1}F$ by $\pi$ coincides with $f$. If we assume in addition that $\Delta(F)=\Delta(f)$ and the $p$-adic distance between $\pi^{-1}F$ and $\hat{f}$ is small enough, then by a strategy similar to that in the proof of \cref{thm:main}, we can show that the filtered $\Phi$-module defined in \cref{rmk:object.cohomology} is weakly admissible.
\end{rmk}

\section{Examples and questions}
\label{section:example}

\begin{eg}\label{eg:curve}
	\par Assume that $p\neq{}2$. Let us consider the case where $n=1$ and $f:\mathbb{T}^{1}_{k}\rightarrow\mathbb{A}^{1}_{k}$ is defined by $t\mapsto{}f(x)\in{}k[x,x^{-1}]$. Here $\mathbb{T}^{1}_{k}=\operatorname{Spec}k[x,x^{-1}]$ is the $1$-dimensional torus over $k$. In this situation, if $f$ is nondegenerate, then so is $\hat{f}$. Assume that $f$ is nondegenerate, we show that
	\begin{equation*}
		\ell_{\mathrm{HT}}(H^{1}_{\mathrm{dR}}(\mathbb{T}^{1}_{K_{1}},\nabla_{\widehat{F}}),F^{\ast}_{\mathrm{irr}})\leq{}1\rlap{.}
	\end{equation*}
	In particular, this implies that \cref{conj:main} is true.
	
	\par If $f(x)$ is a polynomial, namely $f(x)\in{}k[x]$, then $\Delta(f)=[0,d]$, where $d=\deg{}f$. The point $d$ is the only face of $[0,d]$ not containing the origin, so $f$ is nondegenerate if and only if $p\nmid{}d$. Assume that $f$ is nondegenerate. Note that $\mathrm{M}(f)=\mathbb{Z}_{\geq{}0}$. For $u\in\mathbb{Z}_{\geq{}0}$, we have $w(u)=\frac{u}{d}$. Let $\bar{x}$ denote the image of $x$ in $\overline{R}$, then the image of $x\frac{d}{dx}f$ in $\overline{R}$ is $d\alpha_{d}\cdot\bar{x}^{d}$. This implies that $\mathrm{M}_{\mathrm{NP}}=\{0,\dots,d-1\}$. By \cref{rmk:np-module}, we get
	\begin{equation*}
		\ell_{\mathrm{HT}}(H^{1}_{\mathrm{dR}}(\mathbb{T}^{1}_{K_{1}},\nabla_{\widehat{F}}),F^{\ast}_{\mathrm{irr}})=w(d-1)=\frac{d-1}{d}<1\rlap{.}
	\end{equation*}
	If $f(x)\in{}k[x,x^{-1}]$ is a Laurent polynomial with positive degree $d_{1}$ and negative degree $d_{2}$, then we have $\Delta(f)=[-d_{2},d_{1}]$. The points $-d_{2}$ and $d_{1}$ are the only faces of $[-d_{2},d_{1}]$ not containing the origin, so $f$ is nondegenerate if and only if $p\nmid{}d_{1}d_{2}$. Assume that $f$ is nondegenerate. Note that $\mathrm{M}(f)=\mathbb{Z}$. For $u\in\mathbb{Z}$, we have
	\begin{equation*}
		w(u)=\begin{cases}
			\frac{u}{d_{1}}
				&u\geq{}0\rlap{,}\\
			-\frac{u}{d_{2}}
				&u<0\rlap{.}
		\end{cases}
	\end{equation*}
	Note that the image of $x\frac{d}{dx}f$ in $\overline{R}$ is equal to $d_{1}\alpha_{d_{1}}\cdot\bar{x}^{d_{1}}-d_{2}\alpha_{d_{2}}\cdot\bar{x}^{-d_{2}}$. This implies that we can set $\mathrm{M}_{\mathrm{NP}}=\{-d_{2}+1,\dots,0,\dots,d_{1}\}$. By \cref{rmk:np-module}, we get
	\begin{equation*}
		\ell_{\mathrm{HT}}(H^{1}_{\mathrm{dR}}(\mathbb{T}^{1}_{K_{1}},\nabla_{\widehat{F}}),F^{\ast}_{\mathrm{irr}})=w(d_{1})=1\rlap{.}\qedhere
	\end{equation*}
\end{eg}

\begin{eg}\label{eg:counterexample}
	\par Assume that $p\neq{}2$. We look at an example where $n=1$ and $f(x)=x^{2}+x$. As is mentioned in \cref{rmk:transformation}, we show that the automorphism $T_{\mathrm{NP}}$ on $V_{\mathrm{NP}}$ is not filtration-compatible with respect to $F^{\ast}_{\mathrm{NP}}$.
	
	\par Note that $\Delta(f)=[0,2]$, and $\mathrm{M}(f)=\mathbb{Z}_{\geq{}0}$. For $u\in\mathbb{Z}_{\geq{}0}$, we have $w(u)=\frac{u}{2}$. As is discussed in \cref{eg:curve}, we have $\mathrm{M}_{\mathrm{NP}}=\{0,1\}$, so that $V_{\mathrm{NP}}=\langle{}x^{0},x^{1}\rangle_{K_{1}}$. Note that $F^{0}_{\mathrm{NP}}V_{\mathrm{NP}}=\langle{}x^{0}\rangle_{K_{1}}$. Our goal is to show that $T_{\mathrm{NP}}(x^{0})\notin{}F^{0}_{\mathrm{NP}}V_{\mathrm{NP}}=\langle{}x^{0}\rangle_{K_{1}}$. We may write
	\begin{equation*}
		T_{\mathrm{NP},2}(x^{0})=T_{0}\cdot{}x^{0}+T_{1}\cdot\pi^{\frac{1}{2}}x^{1}
	\end{equation*}
	with unique $T_{0},T_{1}\in{}K_{2}$. It suffices to show that $T_{1}\neq{}0$. More precisely, we show $\operatorname{ord}T_{1}<\infty$. Note that $\widehat{F}(x)=\pi(x^{2}+x)$ and $\widetilde{F}(x)=\sum_{l=0}^{\infty}\gamma_{l}(x^{2p^{l}}+x^{p^{l}})$. Let $D_{\widehat{F}}=x\frac{d}{dx}f+2\pi{}x^{2}+\pi{}x$. By \cref{cor:decomposition}, we know that $T_{0}$ and $T_{1}$ are characterized by
	\begin{equation*}
		\exp(\widetilde{F}(x)-\widehat{F}(x))\equiv{}T_{0}+T_{1}\cdot\pi^{\frac{1}{2}}x\bmod{}D_{\widehat{F}}L\left(\frac{1}{p-1}\right)\rlap{.}
	\end{equation*}
	Note that $\exp(\widetilde{F}(x)-\widehat{F}(x))=\exp((\gamma-\pi)(x^{2}+x))\cdot\exp(\sum_{l=1}^{\infty}\gamma_{l}(x^{2p^{l}}+x^{p^{l}}))$. We may write
	\begin{align*}
		E(x)={}&\exp(\widetilde{F}(x)-\widehat{F}(x))=\sum_{u=0}^{\infty}A_{u}\cdot{}x^{u}\rlap{,}\\
		E_{1}(x)={}&\exp((\gamma-\pi)(x^{2}+x))=\sum_{u_{1}=0}^{\infty}A_{1,u_{1}}\cdot{}x^{u_{1}}\rlap{,}\\
		E_{2}(x)={}&\exp(\sum_{l=1}^{\infty}\gamma_{l}(x^{2p^{l}}+x^{p^{l}}))=\sum_{u_{2}=0}^{\infty}A_{2,u_{2}}\cdot{}x^{u_{2}}\rlap{.}
	\end{align*}
	Note that for $u\in\mathbb{Z}_{\geq{}0}$, we have $A_{u}=\sum_{u_{1}+u_{2}=u}A_{1,u_{1}}A_{2,u_{2}}$. By \cref{cor:decomposition}, for $u\in\mathbb{Z}_{\geq{}1}$, there exists unique $R_{u}\in{}K_{2}$, such that
	\begin{equation*}
		\pi^{\frac{u}{2}}x^{u}\equiv{}R_{u}\cdot\pi^{\frac{1}{2}}x\bmod{}D_{\widehat{F}}L\left(\frac{1}{p-1}\right)\rlap{.}
	\end{equation*}
	Note that $T_{1}=\sum_{u=1}^{\infty}\pi^{\frac{u}{2}}A_{u}R_{u}$. By the proof of \cref{lem:valuation}, we have $E_{2}(x)\in{}L(p-1,0)$, which implies $\operatorname{ord}A_{1,u_{1}}\geq(p-1)\cdot\frac{u_{1}}{2}$ for $u_{1}\in\mathbb{Z}_{\geq{}0}$. By the proof of \cite{li2022exponential}*{Proposition~2.1}, we have $E_{2}(x)\in{}L(\frac{p-1}{p},0)$, which implies $\operatorname{ord}A_{2,u_{2}}\geq{}\frac{p-1}{p}\cdot\frac{u_{2}}{2}$ for $u_{2}\in\mathbb{Z}_{\geq{}0}$. Therefore, for $u_{1},u_{2}\in\mathbb{Z}_{\geq{}0}$, we have
	\begin{equation*}
		\operatorname{ord}(\pi^{-\frac{u_{1}+u_{2}}{2}}A_{1,u_{1}}A_{2,u_{2}})\geq\left(p-1+\frac{1}{p-1}\right)\cdot\frac{u_{1}}{2}+\left(\frac{p-1}{p}-\frac{1}{p-1}\right)\cdot\frac{u_{2}}{2}\rlap{.}
	\end{equation*}
	Since $p\neq{}2$, we have $p\geq{}3$, which implies $\frac{p-1}{p}-\frac{1}{p-1}>0$ and $p-1-\frac{1}{p-1}>p-2+\frac{1}{p-1}\cdot\frac{1}{2}$. Therefore, for $u_{1}\geq{}2$ and $u_{2}\in\mathbb{Z}_{\geq{}0}$, we have
	\begin{equation*}
		\operatorname{ord}(\pi^{-\frac{u_{1}+u_{2}}{2}}A_{1,u_{1}}A_{2,u_{2}})\geq{}p-1-\frac{1}{p-1}>p-2+\frac{1}{p-1}\cdot\frac{1}{2}\rlap{.}
	\end{equation*}
	By \cref{lem:valuation.uniformizer}, we have $\operatorname{ord}A_{1,1}=\operatorname{ord}(\gamma-\pi)=p-1+\frac{1}{p-1}$. Therefore, for $u_{2}\in\mathbb{Z}_{\geq{}0}$, we have
	\begin{equation*}
		\operatorname{ord}(\pi^{-\frac{1+u_{2}}{2}}A_{1,1}A_{2,u_{2}})\geq{}p-1+\frac{1}{p-1}\cdot\frac{1}{2}>p-2+\frac{1}{p-1}\cdot\frac{1}{2}\rlap{.}
	\end{equation*}
	For $u_{2}\in\mathbb{Z}_{\geq{}0}$, we have $A_{2,u_{2}}=0$ if $p\nmid{}u_{2}$. Therefore, for $u\in\mathbb{Z}_{\geq{}1}$ such that $p\nmid{}u$, we have
	\begin{equation*}
		\operatorname{ord}(\pi^{-\frac{u}{2}}A_{u})>p-2+\frac{1}{p-1}\cdot\frac{1}{2}\rlap{.}
	\end{equation*}
	When $u=p$, we note that $A_{p}=A_{1,0}A_{2,p}+A_{1,p}A_{2,0}$, so that
	\begin{equation*}
		\operatorname{ord}(\pi^{-\frac{p}{2}}A_{p})=\operatorname{ord}(\pi^{-\frac{p}{2}}\gamma_{1})=p-\frac{3}{2}+\frac{1}{p-1}\cdot\frac{1}{2}>p-2+\frac{1}{p-1}\cdot\frac{1}{2}\rlap{.}
	\end{equation*}
	When $u=2p$, we note that $A_{2p}=A_{1,0}A_{2,2p}+A_{1,p}A_{2,p}+A_{1,2p}A_{2,0}$, so that
	\begin{equation*}
		\operatorname{ord}(\pi^{-p}A_{2p})=\operatorname{ord}(\pi^{-p}\gamma_{1})=p-2\rlap{.}
	\end{equation*}
	If $p\neq{}3$, namely $p\geq{}5$, then we have
	\begin{equation*}
		\frac{p-1}{p}-\frac{1}{p-1}>\frac{2}{3p}\cdot\left(p-2+\frac{1}{p-1}\cdot\frac{1}{2}\right)\rlap{.}
	\end{equation*}
	Therefore, for $u_{2}\in\mathbb{Z}_{\geq{}0}$ such that $u_{2}\geq{}3p$, we have $\frac{u_{2}}{2}\geq\frac{3p}{2}>(p-2+\frac{1}{p-1}\cdot\frac{1}{2})/(\frac{p-1}{p}-\frac{1}{p-1})$, so that
	\begin{equation*}
		\operatorname{ord}(\pi^{-\frac{u_{2}}{2}}A_{2,u_{2}})\geq\left(\frac{p-1}{p}-\frac{1}{p-1}\right)\cdot\frac{u_{2}}{2}>p-2+\frac{1}{p-1}\cdot\frac{1}{2}\rlap{.}
	\end{equation*}
	Hence, for $u\in{}p\mathbb{Z}_{\geq{}0}$ such that $u\geq{}3p$, we have
	\begin{equation*}
		\operatorname{ord}(\pi^{-\frac{u}{2}}A_{u})>p-2+\frac{1}{p-1}\cdot\frac{1}{2}\rlap{.}
	\end{equation*}
	If $p=3$, then $\frac{p-1}{p}-\frac{1}{p-1}=\frac{2}{15}(p-2+\frac{1}{p-1}\cdot\frac{21}{2})$. Thus, for $u_{2}>15$, we have $\frac{u_{2}}{2}>\frac{15}{2}=(p-2+\frac{1}{p-1}\cdot\frac{1}{2})/(\frac{p-1}{p}-\frac{1}{p-1})$, so that
	\begin{equation*}
		\operatorname{ord}(\pi^{-\frac{u_{2}}{2}}A_{2,u_{2}})\geq\left(\frac{p-1}{p}-\frac{1}{p-1}\right)\cdot\frac{u_{2}}{2}>p-2+\frac{1}{p-1}\cdot\frac{1}{2}\rlap{.}
	\end{equation*}
	Therefore, for $u\in{}p\mathbb{Z}_{\geq{}0}$ such that $u>15$, we have
	\begin{equation*}
		\operatorname{ord}(\pi^{-\frac{u}{2}}A_{u})>p-2+\frac{1}{p-1}\cdot\frac{1}{2}\rlap{.}
	\end{equation*}
	It is straightforward to verify that for $u=12,15$, we also have
	\begin{equation*}
		\operatorname{ord}(\pi^{-\frac{u}{2}}A_{u})>p-2+\frac{1}{p-1}\cdot\frac{1}{2}\rlap{.}
	\end{equation*}
	As a summary, when $p\neq{}2$, for $u\in\mathbb{Z}_{\geq{}1}\setminus\{2p\}$, we have $\operatorname{ord}(\pi^{-\frac{u}{2}}A_{u})>p-2+\frac{1}{p-1}\cdot\frac{1}{2}$. On the other hand, we have $R_{0}=0$ and $R_{1}=1$. For $u\geq{}2$, we have
	\begin{equation*}
		\operatorname{ord}R_{u}=\min\left\{\operatorname{ord}R_{u-2}+\operatorname{ord}(u-2),\operatorname{ord}R_{u-1}+\frac{1}{p-1}\cdot\frac{1}{2}\right\}\rlap{.}
	\end{equation*}
	Therefore, for $1\leq{}u\leq{}p-1$, we have $R_{u}=\frac{1}{p-1}\cdot\frac{1+(-1)^{u}}{4}$. Since $p\neq{}2$ so that $2\nmid{}p$, we have $\operatorname{ord}R_{p}=0$ and $\operatorname{ord}R_{p+1}=\frac{1}{p-1}\cdot\frac{1}{2}$. Thus $\operatorname{ord}R_{p+2}=\frac{1}{p-1}$. For $p+1\leq{}u\leq{}2p+1$, we have $\operatorname{ord}R_{u}=\frac{1}{p-1}\cdot\frac{3-(-1)^{u}}{4}$. Since $2\mid{}2p$, we have $\operatorname{ord}R_{2p}=\frac{1}{p-1}\cdot\frac{1}{2}$. Therefore, we have
	\begin{equation*}
		\operatorname{ord}(\pi^{-p}A_{p}R_{p})=p-2+\frac{1}{p-1}\cdot\frac{1}{2}\rlap{.}
	\end{equation*}
	By \cref{prop:decomposition.growth}, we know that $\operatorname{ord}R_{u}\geq{}0$ for $u\in\mapsto{Z}_{\geq{}0}$. Hence, for $u\in\mathbb{Z}_{\geq{}}\setminus\{2p\}$, we have
	\begin{equation*}
		\operatorname{ord}(\pi^{-\frac{u}{2}}A_{u}R_{u})\geq{}\operatorname{ord}(\pi^{-\frac{u}{2}}A_{u})>p-2+\frac{1}{p-1}\cdot\frac{1}{2}=\operatorname{ord}(\pi^{-p}A_{p}R_{p})\rlap{.}
	\end{equation*}
	This implies that $\operatorname{ord}T_{1}=\operatorname{ord}(\pi^{-p}A_{p}R_{p})=p-2+\frac{1}{p-1}\cdot\frac{1}{2}<\infty$, so that $T_{1}\neq{}0$.\qed
\end{eg}

\begin{eg}\label{eg:tensor}
	\par Assume that $p\neq{}2$. Let $f(\boldsymbol{x})=\sum_{i=1}^{n}f_{i}(x_{i})$ with $f_{i}(x_{i})\in{}k[x_{i},x_{i}^{-1}]$. Assume that $f$ is nondegenerate and $\dim\Delta(f)=n$, then $f_{1},\dots,f_{n}$ are nondegenerate. Since $\hat{f}(\boldsymbol{x})=\sum_{i=1}^{n}\hat{f}_{i}(x_{i})$, we know that $\hat{f}$ is nondegenerate. We show that \cref{conj:main} is true for $f$ without requiring $\ell_{\mathrm{HT}}(H^{n}_{\mathrm{dR}}(\mathbb{T}^{n}_{K_{1}},\nabla_{\widehat{F}}),F^{\ast}_{\mathrm{irr}})\leq{}p-2$. Note that
	\begin{equation*}
		f^{\ast}\mathcal{L}_{\pi}\cong{}f_{1}^{\ast}\mathcal{L}_{\pi}\boxtimes\dots\boxtimes{}f_{n}^{\ast}\mathcal{L}_{\pi}\rlap{,}
	\end{equation*}
	where $f_{i}:\mathbb{T}^{1}_{k}=\operatorname{Spec}k[x_{i},x_{i}^{-1}]\rightarrow\mathbb{A}^{1}_{k}$ is the morphism defined by $t\mapsto{}f_{i}(x_{i})$. This gives
	\begin{equation*}
		(H^{n}_{\mathrm{rig}}(\mathbb{T}^{n}_{k}/K_{1},f^{\ast}\mathcal{L}_{\pi}),\phi_{f})\cong\bigotimes_{i=1}^{n}(H^{1}_{\mathrm{rig}}(\mathbb{T}^{1}_{k}/K_{1},f_{i}^{\ast}\mathcal{L}_{\pi}),\phi_{f_{i}})\in\mathbf{Mod}{}^{\Phi}_{K_{1}}\rlap{.}
	\end{equation*}
	At the same time, note that $\nabla_{\widehat{F}}\cong\nabla_{\widehat{F}_{1}}\boxtimes\dots\boxtimes\nabla_{\widehat{F}_{n}}$. Here $\widehat{F}_{i}=\pi\hat{f}_{i}$ for $i=1,\dots,n$. By \cite{chen2018kunneth}*{Theorem~1}, this gives
	\begin{equation*}
		(H^{n}_{\mathrm{dR}}(\mathbb{T}^{n}_{K_{1}},\nabla_{\widehat{F}}),F^{\ast}_{\mathrm{irr}})\cong\bigotimes_{i=1}^{n}(H^{1}_{\mathrm{dR}}(\mathbb{T}^{1}_{K_{1}},\nabla_{\widehat{F}_{i}}),F^{\ast}_{\mathrm{irr}})\in\mathbf{M\overline{F}}_{K_{1}}\rlap{.}
	\end{equation*}
	Combining the argument above, we obtain an isomorphism
	\begin{equation*}
	\begin{tikzcd}
		((H^{n}_{\mathrm{rig}}(\mathbb{T}^{n}_{k}/K_{1},f^{\ast}\mathcal{L}_{\pi}),\phi_{f}),(H^{n}_{\mathrm{dR}}(\mathbb{T}^{n}_{K_{1}},\nabla_{\widehat{F}}),F^{\ast}_{\mathrm{irr}}),\iota_{\widehat{F}})
			\ar[d,phantom,sloped,"\cong"]
			\\
		\bigotimes_{i=1}^{n}((H^{1}_{\mathrm{rig}}(\mathbb{T}^{1}_{k}/K_{1},f_{i}^{\ast}\mathcal{L}_{\pi}),\phi_{f}),(H^{1}_{\mathrm{dR}}(\mathbb{T}^{1}_{K_{1}},\nabla_{\widehat{F}_{i}}),F^{\ast}_{\mathrm{irr}}),\iota_{\widehat{F}_{i}})
	\end{tikzcd}
	\end{equation*}
	in $\mathbf{M\overline{F}}{}^{\Phi}_{K_{1}/K_{1}}$. For $i=1,\dots,n$, let $(V_{i},\hat{\phi}_{i},F^{\ast}_{\mathrm{NP}})$ be the Newton polyhedron module associated with $f_{i}$, and let $\mathrm{M}_{i}\subseteq\mathbb{Z}$ be the finite subset such that $V_{1}=\langle{}x_{i}^{u}\mid{}u\in\mathrm{M}_{i}\rangle_{K_{1}}$. It is straightforward to verify that $\mathrm{M}_{\mathrm{NP}}=\mathrm{M}_{1}\times\dots\times\mathrm{M}_{n}$. Note that the map defined by
	\begin{equation*}
		x_{1}^{u_{1}}\otimes\dots\otimes{}x_{n}^{u_{n}}\mapsto{}x_{1}^{u_{1}}\dots{}x_{n}^{u_{n}}
	\end{equation*}
	gives an isomorphism $\bigotimes_{i=1}^{n}(V_{i},\hat{\phi}_{i},F^{\ast}_{\mathrm{NP}})\rightarrow(V_{\mathrm{NP}},\hat{\phi}_{\mathrm{NP}},F^{\ast}_{\mathrm{NP}})$ in $\mathbf{M\overline{F}}{}^{\Phi}_{K_{1}/K_{1}}$. By the proof of \cref{thm:main} and \cref{eg:curve}, for $i=1,\dots,n$, we know that $(V_{i,m},\hat{\phi}_{i,m},F^{\ast}_{\mathrm{NP}})$ is NP-agreeable, where $V_{i,m}=K_{m}\otimes_{K_{1}}V_{i}$ and $\hat{\phi}_{i,m}=1\otimes\hat{\phi}_{i}$. For $\boldsymbol{u}=(u_{1},\dots,u_{n})\in\mathrm{M}(f)$, we note that $w(\boldsymbol{u})=w_{1}(u_{1})+\dots+w_{n}(u_{n})$, where $w_{i}$ is the weight defined with respect to $f_{i}$. For $\boldsymbol{u}=(u_{1},\dots,u_{n})\in\mathrm{M}(f)$, we also note that $\hat{\phi}_{\mathrm{NP}}(\boldsymbol{x}^{\boldsymbol{u}})=\hat{\phi}_{1}(x_{1}^{u_{1}})\dots\hat{\phi}_{n}(x_{n}^{u_{n}})$. Then, it is straightforward to verify that $(V_{\mathrm{NP},m},\hat{\phi}_{\mathrm{NP},m},F^{\ast}_{\mathrm{NP}})$ is NP-agreeable. Note that for $i=1,\dots,n$, we have $t_{\mathrm{N}}(V_{i,m},\hat{\phi}_{i,m})=t_{\mathrm{H}}(V_{i,m},F^{\ast}_{\mathrm{NP}})$. This implies that $t_{\mathrm{N}}(V_{\mathrm{NP},m},\hat{\phi}_{\mathrm{NP},m})=t_{\mathrm{H}}(V_{\mathrm{NP},m},F^{\ast}_{\mathrm{NP}})$. Thus $(V_{\mathrm{NP},m},\hat{\phi}_{\mathrm{NP},m},F^{\ast}_{\mathrm{NP}})$ is weakly admissible, so that $(V_{\mathrm{NP}},\hat{\phi}_{\mathrm{NP}},F^{\ast}_{\mathrm{NP}})$ is weakly admissible.
	
	\par Now, we consider a specific example where $f(\boldsymbol{x})=\sum_{i=1}^{n}x_{i}^{d}$. In this situation, we have
	\begin{equation*}
		\ell_{\mathrm{HT}}(H^{n}_{\mathrm{dR}}(\mathbb{T}^{n}_{K_{1}},\nabla_{\widehat{F}}),F^{\ast}_{\mathrm{irr}})=\frac{d-1}{d}\cdot{}n\rlap{.}
	\end{equation*}
	If we assume in addition that $n>\frac{d}{d-1}\cdot(p-2)$, then $\ell_{\mathrm{HT}}(H^{n}_{\mathrm{dR}}(\mathbb{T}^{n}_{K_{1}},\nabla_{\widehat{F}}),F^{\ast}_{\mathrm{irr}})>p-2$. However, the argument above shows that \cref{conj:main} is true for $f$. This means that our constraints in \cref{thm:main} is not indispensable.
\end{eg}

\par Let us conclude this article by raising some questions. To be parallel to the story of $p$-adic Hodge theory, we change our perspective and start from a nondegenerate function $F:\mathbb{T}^{n}\rightarrow\mathbb{A}^{1}$ over $O_{1}$. Assume that $\pi^{-1}F$ also gives a regular function over $O_{1}$, and its reduction by $\pi$ gives a nondegenerate function $f:\mathbb{T}^{n}_{k}\rightarrow\mathbb{A}^{1}_{k}$, then it still makes sense to consider the specialization map
\begin{equation*}
	\iota_{F}:H^{i}_{\mathrm{dR}}(\mathbb{T}^{n}_{K_{1}},\nabla_{F})\rightarrow{}H^{i}_{\mathrm{rig}}(\mathbb{T}^{n}_{k}/K_{1},f^{\ast}\mathcal{L}_{\pi})\rlap{.}
\end{equation*}
However, this may not not always be an isomorphism. As is mentioned in \cref{rmk:object.cohomology}, if $\Delta(F)=\Delta(f)$ and the $p$-adic distance between $\pi^{-1}F$ and $\hat{f}$ is less than $p^{-\frac{1}{p-1}}$, then the specialization map $\iota_{F}$ is still an isomorphism. In this situation, we obtain an object
\begin{equation*}
	((H^{n}_{\mathrm{rig}}(\mathbb{T}^{n}_{k}/K_{1},f^{\ast}\mathcal{L}_{\pi}),\phi_{f}),(H^{n}_{\mathrm{dR}}(\mathbb{T}^{n}_{K_{1}},\nabla_{F}),F^{\ast}_{\mathrm{irr}}),\iota_{F})\in\mathbf{M\overline{F}}{}^{\Phi}_{K_{1}/K_{1}}\rlap{.}
\end{equation*}
It is then natural to ask whether this object is weakly admissible. As is mentioned in \cref{rmk:weak-admissibility}, when the $p$-adic distance between $\pi^{-1}F$ and $\hat{f}$ is small enough, we can show that this object is weakly admissible. Our question is that, what could be a sensible constraint for $F$, so that we can expect the associated object to be weakly admissible.

\par Moreover, we consider a global regular function $F$ on a smooth quasi-projective scheme $X$ over $O_{1}$ such that the pair $(X_{K_{1}},F)$ admits a good compactification in the sense of Yu in \cite{yu2014irregular}*{\S1}. Assume further that $\pi^{-1}F$ is a global regular function over $O_{1}$ and the reduction of $\pi^{-1}F$ by $\pi$ gives a regular function $f:X_{k}\rightarrow\mathbb{A}^{1}_{k}$. Then, it makes sense to consider the specialization map
\begin{equation*}
	\iota_{X,F}:H^{i}_{\mathrm{dR}}(X_{K_{1}},\nabla_{F})\rightarrow{}H^{i}_{\mathrm{rig}}(X_{k}/K_{1},f^{\ast}\mathcal{L}_{\pi})\rlap{.}
\end{equation*}
Our question is: what could be the condition for $(X,F)$ to be deserved to be called a \emph{good reduction}, namely $\iota_{X,F}$ is an isomorphism and the associated exponentially twisted cohomology gives a weakly admissible filtered $\Phi$-module. The result in \cite{li2024exponentially} provides a class of pairs $(X,F)$ where $\iota_{X,F}$ is an isomorphism. We can use their result to address these questions. There is also a result concerning the Newton-above-Hodge property given by Kramer-Miller (cf. \cite{kramer2021p}). Updating this result to weakly admissible filtered $\Phi$-modules is an interesting question.

\end{document}